\documentclass[11pt]{article}
\usepackage{newlfont,amsfonts,amssymb, oldgerm,amsmath,amsthm,amsgen,amscd,
}
\usepackage{epic}
\usepackage{eepic}
\usepackage{rotating}
\usepackage[dvips]{epsfig}
\begin{document}
\swapnumbers
\theoremstyle{definition}
\newtheorem{de}{Definition}[section]
\newtheorem{bem}[de]{Remark}
\newtheorem{bez}[de]{Notation}
\theoremstyle{plain}
\newtheorem{lem}[de]{Lemma}
\newtheorem{satz}[de]{Proposition}
\newtheorem{folg}[de]{Corollary}
\newtheorem{theo}[de]{Theorem}

\newcommand{\bde}{\begin{de}}
\newcommand{\ede}{\end{de}}
\newcommand{\ul}{\underline}
\newcommand{\ol}{\overline}
\newcommand{\tbf}{\textbf}
\newcommand{\mc}{\mathcal}
\newcommand{\mb}{\mathbb}
\newcommand{\mf}{\mathfrak}
\newcommand{\bs}{\begin{satz}}
\newcommand{\es}{\end{satz}}
\newcommand{\btheo}{\begin{theo}}
\newcommand{\etheo}{\end{theo}}
\newcommand{\bfolg}{\begin{folg}}
\newcommand{\efolg}{\end{folg}}
\newcommand{\blem}{\begin{lem}}
\newcommand{\elem}{\end{lem}}
\newcommand{\bprf}{\begin{proof}}
\newcommand{\eprf}{\end{proof}}
\newcommand{\bd}{\begin{displaymath}}
\newcommand{\ed}{\end{displaymath}}
\newcommand{\be}{\begin{eqnarray*}}
\newcommand{\ee}{\end{eqnarray*}}
\newcommand{\eeqa}{\end{eqnarray}}
\newcommand{\beqa}{\begin{eqnarray}}
\newcommand{\bi}{\begin{itemize}}
\newcommand{\ei}{\end{itemize}}
\newcommand{\bnum}{\begin{enumerate}}
\newcommand{\enum}{\end{enumerate}}
\newcommand{\la}{\langle}
\newcommand{\ra}{\rangle}
\newcommand{\ve}{\varepsilon}
\newcommand{\vp}{\varphi}
\newcommand{\lra}{\longrightarrow}
\newcommand{\Lra}{\Leftrightarrow}
\newcommand{\Ra}{\Rightarrow}
\newcommand{\sub}{\subset}
\newcommand{\ems}{\emptyset}
\newcommand{\sms}{\setminus}
\newcommand{\ints}{\int\limits}
\newcommand{\sums}{\sum\limits}
\newcommand{\lims}{\lim\limits}
\newcommand{\bcup}{\bigcup\limits}
\newcommand{\bcap}{\bigcap\limits}
\newcommand{\beq}{\begin{equation}}
\newcommand{\eeq}{\end{equation}}
\newcommand{\einhalb}{\frac{1}{2}}
\newcommand{\rr}{\mathbb{R}}
\newcommand{\rn}{\mathbb{R}^n}
\newcommand{\ccc}{\mathbb{C}}
\newcommand{\cn}{\mathbb{C}^n}
\newcommand{\M}{{\cal M}}
\newcommand{\drehgleich}{\mbox{\begin{rotate}{90}$=$  \end{rotate}}}
\newcommand{\turngleich}{\mbox{\begin{turn}{90}$=$  \end{turn}}}
\newcommand{\turnsimeq}{\mbox{\begin{turn}{270}$\simeq$  \end{turn}}}
\newcommand{\vf}{\varphi}
\newcommand{\earr}{\end{array}\]}
\newcommand{\barr}{\[\begin{array}}
\newcommand{\bvec}{\left(\begin{array}{c}}
\newcommand{\evec}{\end{array}\right)}
\newcommand{\sumk}{\sum_{k=1}^n}
\newcommand{\sumi}{\sum_{i=1}^n}
\newcommand{\suml}{\sum_{l=1}^n}
\newcommand{\sumj}{\sum_{j=1}^n}
\newcommand{\suminf}{\sum_{k=0}^\infty}
\newcommand{\inv}{\frac{1}}
\newcommand{\wzbw}{\hfill $\Box$\\[0.2cm]}
\newcommand{\lag}{\mathfrak{g}}
\newcommand{\+}{\oplus}
\newcommand{\x}{\times}
\newcommand{\lx}{\ltimes}
\newcommand{\rrn}{\mathbb{R}^n}
\newcommand{\laso}{\mathfrak{so}}
\newcommand{\lason}{\mathfrak{so}(n)}
\newcommand{\w}{\omega}
\newcommand{\pmh}{{\cal P}(M,h)}
\newcommand{\s}{\sigma}
\newcommand{\deri}{\frac{\partial}}
\newcommand{\ddx}{\frac{\partial}{\partial x}}
\newcommand{\ddz}{\frac{\partial}{\partial z}}
\newcommand{\ddi}{\frac{\partial}{\partial y_i}}
\newcommand{\ddk}{\frac{\partial}{\partial y_k}}
\newcommand{\xz}{^{(x,z)}}
\newcommand{\mh}{(M,h)}
\newcommand{\wxz}{W_{(x,z)}}
\newcommand{\qmh}{{\cal Q}(M,h)}
\newcommand{\bbem}{\begin{bem}}
\newcommand{\ebem}{\end{bem}}
\newcommand{\bbez}{\begin{bez}}
\newcommand{\ebez}{\end{bez}}
\newcommand{\pr}{pr_{\lason}}
\newcommand{\huts}{\hat{\s}}
\newcommand{\whut}{\w^{\huts}}
\newcommand{\bhg}{{\cal B}_H(\lag)}
\newcommand{\aaa}{\alpha}
\newcommand{\bb}{\beta}
\newcommand{\lam}{\lambda}
\newcommand{\LL}{\Lambda}
\newcommand{\D}{\Delta}
\newcommand{\ß}{\beta}
\newcommand{\ä}{\alpha}
\newcommand{\W}{\Omega}

\bibliographystyle{alpha}


\title{Towards a classification of Lorentzian holonomy groups}

\author{Thomas Leistner
}

\maketitle

\begin{abstract}
If the holonomy representation of an $(n+2)$--dimensional simply-connected Lorentzian manifold
$(M,h)$ admits
a degenerate invariant subspace its holonomy group is contained in
the parabolic group $( \mathbb{R} \times
SO(n) )\ltimes \mathbb{R}^n$. The main ingredient of such a holonomy group is
the $SO(n)$--projection $G:=pr_{SO(n)}(Hol_p(M,h))$ and one may ask whether it has to be
a Riemannian holonomy group. In this paper we show that this is the case if $G\subset U(n/2)$
or if the irreducible acting components of $G$ are simple.
\end{abstract}
\setcounter{tocdepth}{1}
\tableofcontents

\section*{Introduction}
The very first step in a classification of the holonomy groups of semi-Riemannian manifolds is the
decomposition theorem of de Rham and Wu (\cite{derham52} for Riemannian manifolds
and \cite{wu64} for general semi-Riemannian manifolds).
It asserts that every simply-connected, complete semi-Riemannian manifold is
isometric to a product of    simply-connected, complete {semi-Riemannian}
manifolds, of which one can be flat and all other are indecomposable (often called
``weakly-irreducible'', i.e. with
no non-degenerate invariant subspace under holonomy representation).  For a
Riemannian manifold  this theorem asserts that the  holonomy representation is
completely reducible, i.e. decomposes into factors which are trivial or irreducible,
and are again Riemannian holonomy representations.  For
pseudo-Riemannian manifolds  indecomposability is not the same as
irreducibility. We can have degenerate invariant subspaces under holonomy
representation.

On the other hand all irreducible factors are known by the
Berger classification of possible irreducible semi-Riemannian holonomy groups
(\cite{berger55}, \cite{simons62}, \cite {alekseevskii68},
\cite{brown-gray72}
and \cite{bryant87}).
This classification uses an algebraic condition which has to be satisfied by
every holonomy group of a torsionfree connection.
It follows from the first Bianchi identity and the Amrose-Singer holonomy theorem \cite{as}
and can be
formulated very easily: If $\mf{h}$ is the Lie algebra of the holonomy group of a torsionfree connection,
acting on the vector space $E\simeq T_p M$, then it obeys
$\mf{h}=\left\{
R(u,v)\ |\ u,v\in V, R\in{\cal K}(\mf{h})\right\}$, where
\[{\cal K}(\mf{h})\ :=\  \left\{R\in\wedge^2V^*\otimes \mf{h}\ |\
R(u,v)w+R(v,w)u+R(w,u)v=0\mbox{ for all }u,v,w\in V \right\}\]
is the space of curvature endomorphisms.

Lie algebras satisfying this conditions are called Berger algebras.
All irreducible Berger algebras are classified in \cite{schwachhoefer1} and \cite{schwachhoefer2}.

For non-irreducible, indecomposable holonomy representations (resp. Berger algebras) such a classification is
missing.

For a  Lorentzian manifold $(M,h)$ of dimension $m>2$ the de
Rham/Wu--decomposition yields the following two cases:
\begin{description}
\item{Completely reducible:} Here $(M,h)$ decomposes  into
irreducible or flat
Riemannian
manifolds and a manifold which is an irreducible or flat
Lorentzian manifold or $(\mathbb{R},-dt)$.
The irreducible Riemannian holonomies are known, as well as the irreducible
Lorentzian one, which has to be the whole $SO(1,m-1)$. (The latter follows from the Berger list but
was directly proved by \cite{olmos-discala01}.)
\item{Not completely reducible:} This is equivalent to the existence of a
degenerate invariant subspace and entails the
existence of exactly one holonomy invariant lightlike subspace. The Lorentzian manifold decomposes into
 irreducible or
flat Riemannian manifolds and a Lorentzian manifold with indecomposable, but non-irreducible
holonomy representation, i.e. with  invariant lightlike (i.e.
{one-dimensional})  subspace.
\end{description}
Thus in order to classify holonomy groups of simply-connected Lorentzian manifolds one
has to find the possible holonomy groups of indecomposable, but non-irreducible
Lorentzian manifolds.

The holonomy algebra of such a manifold of dimension $m:=n+2>2$  is contained in $
( \mathbb{R} \oplus \mathfrak{so}(n) )\ltimes
\mathbb{R}^n$.  L. Berard-Bergery and A. Ikemakhen studied in \cite{bb-ike93}
the projections of
such a holonomy algebra and achieved two important results. The first gives a
classification into four types based on the possible
projections on $\mathbb{R}$ and $\mathbb{R}^n$. For two of these types  the projections are coupled
and for the remaining two uncoupled to the $\lason$--component.

The second result is a
decomposition property for the $\mathfrak{so}(n)$--projection
(see theorem \ref{theoI}), i.e. there is a decomposition of the representation space into irreducible components
and of the Lie algebra into ideals which act irreducible on the components.

The relation between the $\lason$--part and the $\rr$-- and $\rrn$--parts is understood quite well
(\cite{boubel00}, or very recently \cite{galaev}): If one has a simply-connected, indecomposable, non-irreducible
Lorentzian manifold with holonomy of uncoupled type, then, under certain conditions,
one can construct a Lorentzian
manifold with coupled type holonomy.

Now one may ask: Which algebras can occur as
$\mathfrak{so}(n)$-projection of an indecomposable,
but not-irreducible Lorentzian manifold? Of course it has to satisfy the
decomposition property.
Riemannian holonomy algebras are the first examples, because there is a  method to construct
from a given Riemannian manifold an indecomposable Lorentzian manifold with holonomy of
uncoupled type for which the $\lason$--projection equals to the Riemannian holonomy.
Furthermore one can show that the Lorentzian manifold is a $pp$-wave if and only if
the $\rr$-- and the $SO(n)$--component vanish \cite{leistner01}.

In \cite{leistner02} we derived an algebraic criterion on the $\lason$--component of an indecomposable,
non-irreducible, simply-connected Lorentzian manifold $(M,h)$, in analogy to the well known Berger criterion
for holonomy algebras. If $\lag$ is the $\lason$-component of an indecomposable,
non-irreducible, simply-connected Lorentzian manifold, acting on an $n$--dimensional Euclidean vector space
$(E,h)$ then it obeys
$\mf{g}=\left\{
Q(u) | Q\in {\cal B}_{h}(\lag), u\in E\right\}$ where ${\cal B}_{h}(\lag)$ is defined as follows
\[
{\cal B}_{h}(\lag) =
 \left\{ Q\in E^*\otimes \mf{g}\ |\
h( Q(u)v,w)+h( Q(v)w,u)+h( Q(w)u,v)=0,\ \forall\   u,v,w\in E\right\}.\]
Since orthogonal Berger algebras do satisfy this criterion we called these algebras weak-Berger algebras.
Furthermore we showed that every irreducible weak-Berger algebra, which is contained in $\mf{u}(n/2)$ is
a Berger algebra, in particular a Riemannian holonomy algebra. This, together with the
decomposition property implies that $\lag:=pr_{\lason}\mf{hol}_p(M,h)$ is
a Riemannian holonomy algebra if it is contained in $\mf{u}(n/2)$.

In the present paper we prove the  following: If $\lag$ is a simple weak-Berger algebra, not contained in
$\mf{u}(n/2)$, which acts irreducible on $\rrn$, then it is a Berger algebra, and in particular a Riemannian holonomy algebra.
This of course applies to the
irreducible components of the $\lason$--projection of $\mf{hol}_p(M,h)$. In the proof we
proceed analogously to \cite{schwachhoefer2}, where the holonomy groups of
torsion free connections are classified.
This will be the main part of this paper and is contained in section \ref{sectionzwei}.

In the first section we recall the results of \cite{bb-ike93} and our results from \cite{leistner02}
introducing the
notion of weak-Berger algebras.
The third section presents again for sake of completeness
the proof of the fact that weak-Berger algebras in $\mf{u}(n/2)$
are Berger algebras. In the appendix we recall facts about representations of real Lie algebras.

These results leave open the question: Are there semisimple, non simple,
irreducible acting Lie algebras, not contained in $\mf{u}(n/2)$, which are weak-Berger, but not Berger?
We guess that this is not the case, and we intent to apply the
methods of the present paper also in the semisimple case.
Up to dimension eleven this was proved very recently by \cite{galaev} also for algebras not contained in
$\mf{u}(n/2)$.
In his paper he studied the space of curvature endomorphisms for subalgebras in
$( \mathbb{R} \oplus \mathfrak{so}(n) )\ltimes \mathbb{R}^n$ which are of the types found in
\cite{bb-ike93}.
Reducing everything to one uncoupled type he proved the other direction of our result:
a subalgebra  of one of these types is a Berger algebra,
if its $\lason$--projection is a weak-Berger algebra (in our terms).

We are aware that the proofs we will present here
are a cumbersome case-by-case analysis using the methods of representation theory.
It is very desirable to get a direct and more geometric proof of the proposition
that every $SO(n)$--projection of an indecomposable, non-irreducible Lorentzian holonomy group
is a Riemannian holonomy group, which includes the remaining semisimple case of course.

We want to remark that the starting point of this investigation was the question
for the existence of parallel spinors on Lorentzian manifolds.
Such a spinor defines a parallel vector field which can be light like. Hence the manifold
has an indecomposable, non-irreducible factor. But the existence of parallel spinors on
indecomposable Lorentzian manifolds with parallel lightlike vector field
depends only on the $SO(n)$--projection.
Thus a complete list of the latter would answer this question. In the physically important dimensions
below twelve the question for the maximal indecomposable Lorentzian
holonomy groups admitting parallel spinors is answered \cite{bryant00}, \cite{farrill99mw}.

\section{Indecomposable Lorentzian holonomy and weak-Berger algebras}
\subsection{Basic properties}

Let $(M,h)$ be  an indecomposable, non-irreducible Lorentzian manifold with $dim
\ M
=n+2>2$. The holonomy group in a point $p\in M$ acting on $T_p M$ --- defined
as the group of parallel displacements along loops starting at $p$ --- then
has a lightlike, one-dimensional invariant subspace $\Xi_p$
which is the fibre of a parallel distribution $\Xi$. This is equivalent to the existence of
a recurrent lightlike vector field. The subspace $\Xi_p^\bot$ also
is
holonomy invariant and the fibre of a parallel distribution $\Xi^\bot$.
(We call a distribution parallel if it is closed under $\nabla_U$ for every $U\in TM$.)

With respect to a basis
\begin{equation}
\label{basis}
\begin{array}{l}
(X, E_1, \ldots E_n, Z) \mbox{ adapted to $\Xi_p\subset
\Xi_p^\bot$, i.e. $X\in \Xi_p, E_i\in  \Xi_p^\bot$}\\
\mbox{with $h(E_i, E_j)= \delta_{ij},
h(Z,Z)=h(Z,E_i)=h(X,E_i)=0$ and $h(X,Z)=1$}
\end{array}
\end{equation}
 the holonomy algebra is contained in
the following Lie algebra
\begin{equation}
\label{holform}
\mathfrak{hol}_p(M,h)\subset
\left\{ \left( \left.
\begin{array}{ccc}
a	 & 	u^t&0	\\
0	 & A &-u	\\
0	 &0^t&-a
\end{array} \right)
\right| a\in \mathbb{R}, u\in \mathbb{R}^n, A \in \mathfrak{so}(n)
\right\}\ = \ (\rr\+\lason )\ltimes \rrn.
\end{equation}

Choosing a different basis of type (\ref{basis})
corresponds to conjugation with an element
in $O(1,n+1)$ which respects the form (\ref{basis}).
Hence the $\mathfrak{so}(n)$--component is uniquely defined
with respect to conjugation in
$O(n)$.

The projections of $ \mathfrak{hol}_p(M,h)$ on the $\mathbb{R}$-- and on the
$ \mathbb{R}^n$--component are well understood. With respect to these projections
there exist four different types (see \cite{ike90}, \cite{bb-ike93} and \cite{ike96}).
For the types $I$ and $II$ the holonomy is equal to $(\rr \oplus \mathfrak{g}) \ltimes \rr^n$ resp.
$\mathfrak{g}\ltimes\rr^n$.
In case of types $II$ and $IV$ the projection on $\rr$ is zero, which implies the existence
not only of a recurrent lightlike vector field but also of a parallel one.
In case of types $III$ and $IV$ the $\rr$-- respectively the $\rrn$-- components are coupled
to the $\lason$--component, or more precisely to its center.

In the following shall be ${\frak g}:= pr_{{\frak
s}{\frak o}(n)} (\mathfrak{hol}_p(M,h) )$.
About ${\frak g}$ the in \cite{bb-ike93} is proved

\begin{theo}
\label{theoI}
\cite{bb-ike93}
Let ${\frak g}:= pr_{\mathfrak{so}(n)} \left(\mathfrak{hol}_p(M,h) \right)$ be the
projection
of the holonomy algebra of an indecomposable, non-irreducible, $n+2$--dimensional Lorentzian
manifold onto  ${\frak s}{\frak o}(n)$.
Then $ {\frak g}$ satisfies the following decomposition property: There
exists a decomposition of $\mathbb{R}^n$ into orthogonal subspaces and of
${\frak g}$ into ideals
\[
\mathbb{R}^n = E_0 \oplus E_1 \oplus \ldots \oplus E_r\ \mbox{ and }\
{\frak g} = {\frak g}_1 \oplus \ldots \oplus {\frak g}_r,\]
such that ${\frak g}$ acts trivial on $E_0$, ${\frak g}_i$ acts irreducible on $E_i$
and trivial on $E_j$ for $i\not=j$.
\end{theo}

This theorem has two important consequences making a further algebraic
investigation of ${\frak g}$
possible.

Irreducible acting, connected subgroups of $SO(n)$ are are closed and therefore compact.
Now by the theorem the group  $G:=pr_{SO(n)}Hol^0_p(M,h)$
decomposes in such irreducible
acting subgroups. Thus we have as first consequence that  $G$ is compact, although the
whole holonomy group must
not be compact (for such examples see also \cite{bb-ike93}).

The second is, that it suffices to study irreducible acting groups or algebras
${\frak g}$, a fact
which is necessary for trying a classification. We will see this in detail in
the following section.

We will describe the local situation briefly.
Locally there are $n$-dimensional Riemannian submanifolds defined via
special coordinates respecting the foliation $\Xi\subset
\Xi^{\bot}$, denoted by
$(x, y_1, \ldots , y_n, z)$ with $ \frac{\partial}{\partial x}\in \Xi$,
$\frac{\partial}{\partial y_i}\in
\Xi^\bot$. The restriction of $h$ to these submanifolds defined by $y_1, \ldots , y_n$
gives  a family of
Riemannian metrics $g_z$ on it, depending only on the coordinate $z$ (since
$ \frac{\partial}{\partial x}\in \Xi$, see
\cite{brinkmann25}, also  \cite{ike96}).

Although these coordinates are unique under certain conditions (see
\cite{boubel00}) it is not clear how the Lie algebra ${\frak g}$ can be obtained
by the holonomies $\mathfrak{hol}_{p(z)}(g_z)$ of the family of metrics $g_z$.
The only known point is, that all these $\mathfrak{hol}_{p(z)}(h_z)$ are
contained in ${\frak g}$ \cite{ike96}.

If the dependence on $z$ is trivial  --- i.e. $g_z\equiv g$ or $g_z\equiv f(z)
g$ --- then ${\frak g}$ is equal to the holonomy of the Riemannian metric $g$.

In particular this gives a way to construct indecomposable, non-irreducible
Lorentzian manifolds with holonomy equal to $ (\mathbb{R} \oplus $Riemannian
holonomy$)\ltimes \mathbb{R}^n$: Let $(N,g)$ be an $n$-dimensional Riemannian
manifold, $\theta $ a closed form on $N$ and $q$ a function on $N\times
\mathbb{R}^2$, the latter sufficiently general. Then
\[(M=N\times \mathbb{R}^2, h=dxdz + q dz^2 + \theta dz+ f (z) g)\]
is a Lorentzian manifold with holonomy
\[ \mathfrak{hol}_{(x,z,p)}(M,h)=(\mathbb{R} \oplus
\mathfrak{hol}_p(N,g))\ltimes \mathbb{R}^n.\]
In case of Riemannian K\"ahler- and hyper-K\"ahler manifolds $(N,g)$ these conditions
can be weakened \cite{leistner01}.

Furthermore there is a  method to construct manifolds with coupled holonomy from manifolds with uncoupled
holonomy \cite{boubel00}: If $(M,h)$ is a simply-connected, indecomposable, non-irreducible Lorentzian manifold
with uncoupled holonomy $\lag\ltimes\rrn$ or $(\rr\+\lag)\ltimes\rrn$ such that
$\lag$ has non-trivial center (and further conditions),
then there is a metric $\tilde{h}$ on $M$ such that $(M,\tilde{h})$ has holonomy of coupled type and
with $\lason$--projection $\lag$.

In the following we will go an algebraic way, in order to classify the possible algebras ${\frak g}$.
This algebraic way uses the Bianchi--identity,
restricted to $\mathbb{R}^n$ as representation space of ${\frak g}$. This is the aim of the next sections.

\subsection{Berger and weak-Berger algebras}

Here we will introduce the notion of
weak-Berger algebras in comparison to Berger algebras. We present some basic
properties,
in particular a decomposition property and the behavior under complexification.
For the details to this section see \cite{leistner02}


Let $E$ be a vector space over the field $\mathbb{K}$ and let
${\frak g}\subset \mathfrak{gl}(E) $ be a Lie algebra. Then one defines
\begin{align*}
{\cal K}({\frak g})&\ :=\  \{ R\in \Lambda^2 E^* \otimes {\frak g} \ |\ R(x,y)z +
R(y,z)x + R(z,x)y=0\}\\
\underline{{\frak g}}&\ :=\ span \{ R(x,y)\ |\  x,y\in E, R\in {\cal K}({\frak g})\},\\
\intertext{and for ${\frak g}\subset \mathfrak{so}(E,h)$:}
{\cal B}_h({\frak g})&\ :=\  \{ Q\in E^* \otimes {\frak g} \ |\ h(Q(x)y,z) + h(Q(y)z,x)
+ h(Q(z)x,y)=0\}\\
{\frak g}_h&\ :=\  span\{ Q(x)\ |\  x\in E, Q\in {\cal B}_h({\frak g})\}.
\end{align*}

Then we have the following basic properties.
\begin{satz}\label{module}
$   {\cal K}({\frak g})\subset \Lambda^2E^* \otimes
{\frak g}$ and ${\cal B}_h({\frak g})
\subset E^*
\otimes {\frak g}$ are $ {\frak g}$-modules. $ \underline{{\frak g}}$ and
${\frak
g}_h$ are ideals in $\lag$.

\end{satz}

The representation of ${\frak g}$ on $ {\cal B}_h({\frak g})$ and $ {\cal
K}({\frak
g})$ is given by the standard and the adjoint representation
\begin{eqnarray}
(A \cdot Q) (x) &=& - Q(Ax) + [A, Q(x)]\label{weakaction}\\
(A \cdot R) (x,y) &=& - R(Ax,y) - R(x,Ay) + [A, R(x,y)].
\end{eqnarray}

\begin{de}\label{weakdef}
Let ${\frak g}\subset \mathfrak{gl}( E)$ be a Lie algebra. Then ${\frak g}$ is
is called {\bf Berger algebra} if $\underline{\frak g} = {\frak g}$. If
${\frak g}\subset \mathfrak{so}(E,h)$ is an orthogonal Lie algebra with $ {\frak g}_h={\frak g}$,
then we
call it {\bf weak-Berger algebra}.
\end{de}

Equivalent to the (weak-)Berger property is the fact that there is no ideal
$\mf{h}$ in $\lag$ such that ${\cal K}(\mf{h})={\cal K}(\mf{\lag})$ (resp. ${\cal B}_h(\mf{h})= {\cal B}_h(\mf{g})$).

The notion ``weak-Berger'' is satisfied by the following
\begin{satz}
Every Berger algebra which is orthogonal is a weak-Berger algebra.
\end{satz}
This proposition has a

 \begin{folg}\label{corollar}
Let ${\frak g}\subset \mathfrak{so}(E,h)$ be an orthogonal Lie algebra. Then
\begin{equation}
\label{folgerung}
span \{R(x,y)+Q(z)|R\in {\cal K}({\frak g}), Q\in {\cal B}_h({\frak g}),x,y,z\in E \}
\subset
 {\frak g}_h.
\end{equation}
\end{folg}

Concerning the decomposition
of Berger and weak-Berger algebras the following proposition holds.

\begin{satz}
If ${\frak g}_1\subset  \mathfrak{gl}(V_1)  $ , ${\frak g}_2 \subset
\mathfrak{gl}(V_2)$ and ${\frak g}:= {\frak g}_1 \oplus {\frak g}_2 \subset
\mathfrak{gl}(V:=V_1 \oplus V_2)$, then it holds:
\begin{enumerate}
\item If ${\frak g}_1$ and $ {\frak g}_2$ are Berger algebras, then ${\frak g}$
is a Berger algebra.
\item
If in addition ${\frak g}_1\subset  \mathfrak{so}(V_1, h_1)  $, ${\frak g}_2
\subset \mathfrak{so}(V_2,h_2)$ and ${\frak g}:= {\frak g}_1 \oplus {\frak g}_2
\subset \mathfrak{so}(V:=V_1 \oplus V_2, h:= h_1 \oplus h_2)$, then holds:

${\frak g}_1$ and ${\frak g}_2$ are weak-Berger algebras if and only if
${\frak g}$ is a weak-Berger algebra.
\end{enumerate}
\label{zerlegung}
\end{satz}


The Ambrose-Singer holonomy theorem \cite{as} then implies that holonomy algebras of
torsion free connections --- in particular of a
Levi-Civita-connection --- are Berger algebras.  The list of all
irreducible Berger algebras is known (\cite{berger55} for orthogonal, non-symmetric Berger algebras,
\cite{berger57} for orthogonal symmetric ones, and \cite{schwachhoefer1} in the
general affine case).

We should mention that in our notation Berger algebras are not only non-symmetric
Berger algebras, as it is sometimes defined. For us only the possibility of
being
the holonomy algebra of a Riemannian manifold is of interest,  symmetric or non symmetric.

The $ {\frak s}{\frak o}(n)$--projection of an indecomposable, non-irreducible Lorentzian manifold
is no holonomy algebra, and therefore not necessarily a Berger algebra. But the following statement,
which we proved in \cite{leistner02},
asserts that it is a
weak-Berger algebra.

\begin{theo}\label{weak-Berger}
Let $(M,h)$ be an indecomposable, but non-irreducible, simply connected Lorentzian manifold
and  ${\frak g}=pr_{{\frak s}{\frak o}(n)}(\mathfrak{hol}_p(M,h))$. Then $ {\frak g}$ is a weak-Berger algebra.
 \end{theo}

From point two of proposition \ref{zerlegung} we get the following

\begin{folg}
Let $(M,h)$ be an indecomposable, but non-irreducible Lorentzian manifold
and  ${\frak g}=pr_{{\frak s}{\frak o}(n)}(\mathfrak{hol}_p(M,h)) \subset E^*
\otimes E$ and $ {\frak g} = {\frak g}_1 \oplus \ldots \oplus {\frak g}_r$ with
$ {\frak g}_i\in \mathfrak{so}(E_i,h_i)$ the decomposition in irreducible acting ideals
from theorem
\ref{theoI}.
Then these $ {\frak g}_i$  are irreducible weak-Berger algebras.
\end{folg}

This corollary ensures that we are at a similar point as in the Riemannian situation,
but reaching it by a different way. This is shown schematically in the following
diagram:

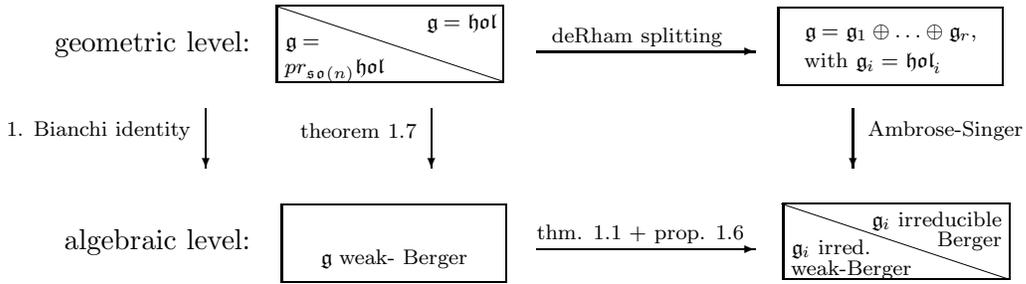
\[ \begin{array}{rccc}
\mbox{ geometric level:}&
\setlength{\unitlength}{1cm}
\begin{picture}(3.1,0.7)
\put(0,-0.4){\framebox(3,1)[tr]{
\put(-0.1,-0.1){\makebox(0,0)[tr]{{\scriptsize $\mathfrak{g}=\mathfrak{hol}$}}}}}
\put(0,-0.4){\framebox(3,1)[bl]{
\put(0.1,0.5){\makebox(0,0)[l]{{\scriptsize $\mathfrak{g}=$}}}
\put(0.1,0){\makebox(0,0)[bl]{{\scriptsize $pr_{\mathfrak{so}(n)}\mathfrak{hol}$}}}}}
\put(0,0.6){\line(3,-1){3}}
\end{picture}
&
\setlength{\unitlength}{1cm}
\begin{picture}(2.8,1)
\put(0,0){\vector(1,0){2.9}}
\put(0.2,0.2){\makebox(0,0)[l]{{\scriptsize deRham splitting}}}
\end{picture}
&\framebox[3cm]{{\scriptsize $\begin{array}{l}
{\frak g}={\frak g}_1 \oplus \ldots \oplus {\frak g}_r,\\
\mbox{with }{\frak g}_i=\mathfrak{hol}_i
\end{array}$}}
\\
\setlength{\unitlength}{1cm}
\begin{picture}(1,1)
\put(0.2,0.4){\makebox(0,0)[r]{{\scriptsize 1. Bianchi identity}}}
\put(0.4,0.7){\vector(0,-1){0.8}}
\end{picture}
&
\setlength{\unitlength}{1cm}
\begin{picture}(1,1)
\put(1,0.7){\vector(0,-1){0.8}}
\put(0.8,0.4){\makebox(0,0)[r]{{\scriptsize theorem \ref{weak-Berger}}}}
\end{picture}
&&
\setlength{\unitlength}{1cm}
\begin{picture}(1,1)
\put(0,0.7){\vector(0,-1){0.8}}
\put(0.2,0.4){\makebox(0,0)[l]{{\scriptsize Ambrose-Singer}}}
\end{picture}
\\
\mbox{algebraic level:}&\framebox[3cm]{{\scriptsize $\begin{array}{c}\\{\frak g}\mbox{ weak-
Berger}\\
\end{array}$}}&
\setlength{\unitlength}{1cm}
\begin{picture}(2.8,1)
\put(0,0){\vector(1,0){2.9}}
\put(0,0.2){\makebox(0,0)[l]{{\scriptsize thm. \ref{theoI} $+$ prop. \ref{zerlegung}}}}
\end{picture}
&
\setlength{\unitlength}{1cm}
\begin{picture}(3.1,0.6)
\put(0,-0.4){
\framebox(3,1)[tr]{
\put(-0.1,-0.1){\makebox(0,0)[tr]{{\scriptsize $\mathfrak{g}_i$ irreducible}}}
\put(-0.1,-0.5){\makebox(0,0)[r]{{\scriptsize
Berger}}}}}
\put(0,-0.4){
\framebox(3,1)[bl]{
\put(0.1,0.4){\makebox(0,0)[l]{{\scriptsize $\mathfrak{g}_i$ irred.}}}
\put(0.1,0){\makebox(0,0)[bl]{{\scriptsize weak-Berger}}}}}
\put(0.1,0.6){\line(3,-1){3}}
\end{picture}\\
&&&
\end{array}\]

The proof of the theorem gives another

\begin{folg}
Let $(M,h)$ be an indecomposable, non-irreducible Lorentzian manifold and
${\frak g}= pr_{ \mathfrak{so}(n)} \mathfrak{hol}_p(M,h)$. If there exists
coordinates $(x, y_1, \ldots , y_n,z)$ of the above form
(i.e. respecting the foliation $\Xi\subset \Xi^\perp$), with the property that
everywhere holds
${\cal R}(\frac{\partial}{\partial z}, \frac{\partial}{\partial y_i}, \frac{\partial}{\partial y_j},
\frac{\partial}{\partial y_k})=0$, then ${\frak g}$ is a
Berger-algebra.
\end{folg}
The aim of the following sections will be to classify all weak-Berger algebras. Before we do this
we have to say a word about real and complex (weak-) Berger algebras.

\subsection{Real and complex weak-Berger algebras}

Because of the above result we have to classify the real weak-Berger algebras. Since we will use the
representation theory of complex semisimple Lie algebras we have to describe the
transition of a real weak-Berger algebra to its complexification.

First we note that the spaces $ {\cal K}({\frak g})$ and ${\cal B}_h({\frak g})$ for
${\frak g}\subset
\mathfrak{so}(E,h)$ can be described by the following exact sequences:
\[
\begin{array}{rcccccl}
0 &\rightarrow & {\cal K}({\frak g}) & \hookrightarrow& \wedge^2 E^* \otimes
{\frak g} &\stackrel{\lambda}{\twoheadrightarrow} &\wedge^3 E^* \otimes E\\
0 &\rightarrow & {\cal B}_h({\frak g}) & \hookrightarrow &  E^* \otimes {\frak
g}
&\stackrel{\lambda_h}{\twoheadrightarrow}& \wedge^3 E^*,
\end{array}\]
where the  map $ \lambda$ is the skew-symmetrization and $\lambda_h$ the
dualization by $h$ and the skew-symmetrization.

If we now consider a real Lie algebra ${\frak g}$ acting orthogonal on a real vector
space
$E$, i.e. ${\frak g}\subset \mathfrak{so}(E,h)$, then $h$ extends by
complexification (linear in both components) to a non-degenerate complex-bilinear form
$h^\mathbb{C}$ which is invariant under ${\frak g}^\mathbb{C}$, i.e. ${\frak
g}^\mathbb{C}\subset \mathfrak{so}( E^\mathbb{C}, h^\mathbb{C})$. Then the
complexification of the above exact sequences gives
\begin{eqnarray}
\label{complex1}
{\cal K}({\frak g})^\mathbb{C}&=& {\cal K}({\frak g^\mathbb{C}})\\
\label{complex2}
\left({\cal B}_h({\frak g})\right)^\mathbb{C}&=& {\cal B}_{h^\mathbb{C}}({\frak
g^\mathbb{C}}).
\end{eqnarray}
 This implies
\begin{satz} \label{real-complex}
${\frak g}\subset \mathfrak{so}(E,h)$ is a (weak-) Berger algebra if and only if
${\frak
g}^\mathbb{C}\subset \mathfrak{so}(E^\mathbb{C}, h^\mathbb{C})$ is a (weak-) Berger algebra.
\end{satz}

I.e. complexification preserves the weak-Berger as well as the Berger property.

Because of proposition \ref{zerlegung} it suffices to classify the  real weak-Berger algebras which are irreducible.
Now irreducibility is a property which is not preserved under complexification. We have to deal with this problem.
At a first step one recalls the following definition, distinguishing two cases for a module of a
real Lie algebras.

\bde
Let $\lag$ be a real Lie algebra.
Irreducible real $\lag$-modules $E$ for which $E^\ccc$ is an irreducible $\lag$-module and
irreducible complex modules $V$ for which $V_\rr$ is a reducible $\lag$-module are
called of {\bf real type}.
Irreducible real $\lag$-modules $E$ for which $E^\ccc$ is a reducible $\lag$-module and
irreducible complex modules $V$ for which $V_\rr$ is a irreducible $\lag$-module are
called of {\bf non-real type}.
\ede

This notation corresponds to the distinction of complex irreducible $\lag$-modules
into real, complex and quaternionic ones. It makes sense because the complexification
of a module of real type is of real type --- recall that
$(E^\ccc)_\rr$ is a reducible $\lag$-module --- and the reellification of a module of non-real type is
of non-real type. These relations are described in the appendix
\ref{real representations}.

In the original papers of Cartan \cite{cartan1914} and Iwahori \cite{iwahori59},  see also \cite {goto78},
where these distinction is introduced,
a representation of real type is called as representation of {\bf first type} and a representation of
non-real type is called of {\bf second type}.

If one now complexifies the Lie algebra $ {\frak g}$ too, then $E^\ccc$ becomes a
$\lag^\ccc$--module. This transition preserves irreducibility.
\begin{lem}
Let $ {\frak g}^\mathbb{C}\subset \mathfrak{gl}(V)$ be the complexification of
${\frak g}\subset \frak{gl}(V)$ with a complex $\lag$-module $V$. Then it holds:
\begin{enumerate}
\item ${\frak g}$ is irreducible if and only if ${\frak g}^\mathbb{C}$ is
irreducible.
\item ${\frak g}\subset \mathfrak{so}(V,H)$ if and only if ${\frak
g}^\mathbb{C}\subset \mathfrak{so}(V, H)$, where $H$ is a symmetric bilinear form.
\end{enumerate}
\end{lem}
In the following sections
we will describe the weak-Berger property for real and non-real modules of a real Lie algebra $\lag$.

\section{Weak-Berger algebras of real type}
\label{sectionzwei}
\label{complex weak berger}
In this section we will make efforts to classify weak-Berger algebras of real type, at least the
simple ones.
The argumentation in this section is analogously to the reasoning in \cite{schwachhoefer1}.

\bigskip

$\lag_0$ shall be a real Lie algebra and $E$ a real irreducible module of real type.
Furthermore we suppose $\lag_0\in \laso(E,h)$ with $h$ positive definite. Then $E^\ccc$
is an irreducible $\lag_0$-module (also of real type). If $h^\ccc$ denotes the complexification of $h$,
bilinear in both components we have that $\lag_0\subset \laso(E^\ccc,h^\ccc)$.

Now we can extend $h$ also sesqui-linear on $E^\ccc$ and get a hermitian form $\theta^h$ on $V$
which is invariant under $\lag_0$. Thus we have $\lag_0\subset \mf{u}(V,\theta^h).$ $\theta^h$ has the same index as
$h$ (see appendix \ref{real representations}).

Since the bilinear form $h$ we start with is positive definite we can make another simplification.
Subalgebras of $\laso(E,h)$ with positive definite $h$ are compact and therefore reductive.
I.e. its Levi-decomposition is $\lag_0=\mf{z}_0\+\mf{d}_0$, with center $\mf{z}_0$ and semisimple derived
algebra $\mf{d}_0$.
Thus $\lag_0^\ccc=\mf{z}\+ \mf{d}$ is also reductive. But since it is irreducible by assumption,
the Schur lemma implies that the center $\mf{z}$ is $\ccc\  Id$ or zero. But $\ccc \ Id$ is not contained in
$\laso(V,H)$. Thus the center has to be zero and $\lag$ is semisimple.
Proposition \ref{real-complex} gives the following.

\bs
If ${\frak g}_0\subset \mathfrak{so}(E,h)$ is a weak-Berger algebra of real type
then, $ {\frak g}_0^\mathbb{C}\subset \mathfrak{so}(E^\mathbb{C}, h^\mathbb{C})$ is an irreducible weak-Berger
algebra. $E^\ccc$ is a $\lag_0$-module of real type and if $h$ is positive definite then
$\lag_0^\ccc $ is semisimple.

If $ {\frak g}\subset \mathfrak{so}(V, H)$ is an irreducible complex weak-Berger which is semisimple. Then
$\lag$ has a compact real form $\lag_0$ and if
$V$ is a $\lag_0$-module of real type, then
$V=E^\ccc$, $\lag_0$ is unitary with respect to a hermitian form $\theta$ and
${\frak g}_0\subset \mathfrak{so}(E,h)$ is a weak-Berger algebra of real type. The indices of $h$ and $\theta$
are equal.
\es

\bprf
The first direction follows obviously from proposition \ref{real-complex}. That $E^\ccc$ is a module of real type
holds because of $(E^\ccc)_\rr$ is reducible (see appendix proposition \ref{realtype1}).

Since $\lag$ is semisimple it has a compact real form $\lag_0$.
If  $V$ is a $\lag_0$-module of real type, then $\lag_0$ is unitary since it is orthogonal
(see proposition \ref{dualtype}) and it is $V=E^\ccc$ (proposition \ref{realtype1}).
By proposition \ref{realtypeorthogonal} follows that $\lag_0$ is orthogonal w.r.t.
$h$ which has the same index as $\theta$.
Then the proposition follows by proposition \ref{real-complex}.
\eprf

The main point of this proposition is the implication that if $\lag_0\subset\laso(E,h)$ is weak-Berger
of real type, then $\lag_0^\ccc\subset\laso(E^\ccc, h^\ccc)$ is  an irreducible acting, complex
semisimple weak-Berger algebra. These we have to classify.

\bbem
Before we start we have to make a remark about definition of holonomy up to conjugation.
The $SO(n)$--component of an indecomposable, non-irreducible Lorentzian manifold was defined
modulo conjugation in $O(n)$.
Hence we shall not distinct between subalgebras of $\mf{gl}(n,\ccc)$ which
are isomorphic under $Ad_\vf$ where $\vf$ is an element from $O(n,\ccc)$ and
$Ad$ the adjoint action in of $Gl(n,\ccc)$ on $\mf{gl}(n,\ccc)$.
We say that
an orthogonal representation  $\kappa_1$ of a complex semisimple Lie algebra $\lag$
 is {\bf congruent} to an orthogonal representation $\kappa_2$ if there is an element
 $\vf\in O(n,\ccc)$ such that the following equivalence of
$\lag$--representations is valid: $\kappa_1\sim Ad_\vf\circ\kappa_2$.
Hence we have to classify semisimple, orthogonal, irreducible acting,
complex weak-Berger algebras of real type up to this congruence of representations.

If the automorphism $Ad_\vf$ is inner, then the representations are equivalent, if it is outer
then only congruent.

For semisimple Lie algebras it holds that $Out(\lag):=Aut(\lag)/Inn(\lag)$ counts the connection
components of $Aut(\lag)$ and (see for example \cite{onishchik-vinberg3})
$Out(\lag)$ is isomorphic to the automorphism of the fundamental system, i.e. symmetries of the
Dynkin diagram.
\\[.2cm]
\begin{minipage}[b]{10.5cm}{
For us this becomes relevant in case of $\laso(8,\ccc)$.
In the picture one sees that the symmetries of the Dynkin diagram
generate the symmetric group ${\cal S}_3$, i.e.
$Out(\laso(8,\ccc))={\cal S}_3$ and it contains
the so-called ``triality
automorphism'' which interchanges vector and spin representations of
$ \laso(8,\ccc)$ without fixing one.
}
\end{minipage}
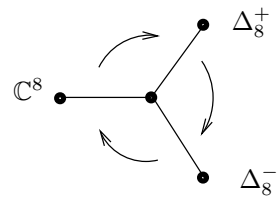
\hfill
\begin{minipage}[b]{3.5cm}{
\begin{picture}(0,0)%
\includegraphics{trias.pstex}%
\end{picture}%
\setlength{\unitlength}{4144sp}%
\begingroup\makeatletter\ifx\SetFigFont\undefined%
\gdef\SetFigFont#1#2#3#4#5{%
  \reset@font\fontsize{#1}{#2pt}%
  \fontfamily{#3}\fontseries{#4}\fontshape{#5}%
  \selectfont}%
\fi\endgroup%
\begin{picture}(1355,1110)(86,-556)
\put(1396,389){\makebox(0,0)[lb]{\smash{\SetFigFont{10}{12.0}{\rmdefault}{\mddefault}{\updefault}$\Delta_8^+$}}}
\put(1441,-556){\makebox(0,0)[lb]{\smash{\SetFigFont{10}{12.0}{\rmdefault}{\mddefault}{\updefault}$\Delta_8^-$}}}
\put( 86,-26){\makebox(0,0)[lb]{\smash{\SetFigFont{10}{12.0}{\rmdefault}{\mddefault}{\updefault}$\mathbb{C}^8$}}}
\end{picture}
}
\end{minipage}\\[.1cm]
We will use that the automorphism
which interchanges  the
vector representation with one spinor representation and
fixes the second spinor representation resp. interchanges the spinor representations and
fixes the vector representation comes from $Ad_\vf$ with $\vf\in O(n,\ccc)$.
Hence the vector and the spinor representations of
$\laso(8,\ccc)$ are congruent to each other.

Finally we should remark that  compact real forms equivalent to a given one correspond to inner
automorphism of $\lag$. Hence the corresponding representations are equivalent
\ebem

\subsection{Irreducible, complex, orthogonal, semisimple Lie algebras}

In the following $V$ will be a complex vector space equipped with a non-degenerate symmetric bilinear
2--form $H$. $\lag$ shall be an irreducible acting, complex, semisimple subalgebra of $\laso(V,H)$.

Thus all the tools of root space decomposition and representation theory will apply.
Let $\mf{t}$ be the Cartan subalgebra of $\lag$.
We denote by $\Delta\subset \mf{t}^*$ the roots of $\lag$ and we set $\Delta_0:= \Delta\cup \{0\}$. Then
$\lag $ decomposes into its root spaces
$\lag_\alpha:=\{ A\in \lag| [T, A]=\alpha(T) \cdot A\mbox{ for all } T\in \mf{t}\}\not=\{0\}$. It is
\[\lag =\bigoplus_{\alpha\in \Delta_0} \lag_\alpha
\;\;\;\mbox{ where $\lag_0=\mf{t}$.}\]

By $\Omega\subset \mf{t}^*$ we denote  the weights of $\lag\subset \laso (V,H)$. Then $V$ decomposes into
the weight spaces
$V_\mu:=\{ v\in V| T(v)=\mu(T)\cdot v\mbox{ for all } T\in \mf{t}\}\not=\{0\}$, i.e.
\[V=\bigoplus_{\mu\in \Omega} V_\mu.\]

Now  the following holds.

\begin{satz}
\label{ortogonal}
Let $ {\frak g}\subset \laso (V,H)$ be a complex, semisimple Lie algebra with weight space decomposition.
Then
\[ V( \mu)\bot V(\lambda)\mbox{ if and only if } \lambda\not= -
\mu.\]
In particular if $\mu$ is a weight, then $-\mu$ too.
\end{satz}

\bprf
For any $T\in \mf{t}$, $u\in V_\mu $ and $v\in V_\lambda$ we have
\[0=
H(Tu,v)+H(u,Tv)=\left(\mu(T) + \lambda(T)\right) H(u,v).\]
Now if $\lambda\not=-\mu$ there is a $T$ such that $\mu(T) + \lambda(T)\not=0$. But this implies
$V_\lambda \bot V_\mu$.

On the other hand $V_\mu\bot V_{-\mu}$ would imply $V_\mu\bot V$ which contradicts the non-degeneracy of $H$.

Its non-degeneracy also implies that $\mu\in \mf{t}^*$ is a weight if and only if $- \mu$ is a weight.
\eprf

\subsection{Irreducible complex weak-Berger algebras}

We will now draw consequences from the weak-Berger property.  Therefore we consider the space
${\cal B}_H(\lag)$ defined by the Bianchi identity. If $\lag$ is weak-Berger it has to be non-zero, i.e.
by proposition \ref{module} it is a non-zero $\lag$--module. If we denote by $\Pi$ all its weights then it
decomposes into weight spaces
\[ \bhg\ = \ \bigoplus_{\phi\in \Pi}{\cal B}_\phi.\]
If  $\Omega$  are the weights of $V$ then we define the following set
\[\Gamma := \left\{\ \mu+\phi\ \left|
\begin{array}{l}
\mu\in \Omega,\ \phi \in \Pi\
\mbox{ and there is an $u\in V_\mu$}\\
\mbox{and a $Q\in {\cal B}_\phi$ such that $Q(u)\not=0$}
\end{array}\right.\right\}\subset \mf{t}^*.\]
Then one proves a      

\blem
$\Gamma\subset \Delta_0$.
\elem

\bprf
We have to show that every $\mu+\phi\in \Gamma$ is a root of $\lag$. Therefore we consider weight elements
$Q_\phi\in {\cal B}_\phi$ and $u_\mu\in V_\mu$ with $0\not=Q_\phi(u_\mu)$.
Then  for every $T\in \mf{t}$ holds (because of the definition of the $\lag$-module $\bhg$):
\be
\left[ T, Q_\phi(u_\mu)\right]&=& (TQ_\phi)(u_\mu)+ Q_\phi(T(u_\mu))\\
&=&\left(\phi(T) + \mu(T)\right) Q_\phi(u_\mu)
\ee
I.e. $\phi+\mu$ is a root or zero.
\eprf

For weak-Berger algebras now the other inclusion is true.

\bs
If $\lag\subset \laso(V,h)$ is irreducible, semisimple Lie algebra. If it is weak-Berger then
$\Gamma =\Delta_0$. If $\Gamma =\Delta_0$ and $span\{Q_{\-\mu}(u_\mu)\ |\ \mu\in\Omega\}=\mf{t}$
then it is weak-Berger.
\es

\bprf
The decomposition of $\bhg$ and $V$ into weight spaces and the fact that
$Q_\phi(u_\mu)\in \lag_{\phi+\mu}$ imply the following inclusion
\[\lag_H\ =\ span\{Q_\phi(u_\mu)\ |\ \phi+ \mu\in \Gamma\}\subset
\bigoplus_{\beta\in \Gamma}\lag_\beta.\]
But if $\lag=\bigoplus_{\alpha\in \Delta_0}\lag_\alpha$ is weak-Berger it holds that $\lag \subset \lag_H$
and thus
\[\bigoplus_{\alpha\in \Delta_0}\lag_\alpha\subset
\bigoplus_{\beta\in \Gamma}\lag_\beta
\subset
\bigoplus_{\alpha\in \Delta_0}\lag_\alpha.
\]
But this implies $\Gamma = \Delta_0$.

If now $\Gamma = \Delta_0$ and $span\{Q_{\-\mu}(u_\mu)\ |\ \mu\in\Omega\}=\mf{t}$ we have that
\[\lag_H\ =\ span\{Q_\phi(u_\mu)\ |\ \phi+ \mu\in \Gamma\}=
\bigoplus_{\beta\in \Gamma}\lag_\beta= \mf{t}\oplus\bigoplus_{\beta\in \Delta}\lag_\beta =\lag .\]
This completes the proof.\eprf

To derive  necessary conditions for the weak Berger property we have to fix a notation.
Let $\alpha\in \Delta$ be a root. Then we denote by $\Omega_\alpha$ the following subset of $\Omega$:
\[\Omega_\alpha:= \left\{ \lambda\in \Omega\ |\ \lambda+\alpha\in \Omega\right\}.\]
Then of course $\alpha+\Omega_\alpha$ are the weights of $\lag_\alpha V$.

\bs\label{musatz}
Let $\lag$ be a semisimple
Lie algebra with roots $\Delta$ and $\Delta_0=\Delta\cup\{0\}$. Let $\lag\subset
\laso(V,H)$ irreducible, weak-Berger with weights $\Omega$. Then the following properties are satisfied:
\begin{description}
\item[(PI)] There is a $\mu\in \Omega$ and a hyperplane $U\subset \mf{t}^*$ such that
\beq \label{mu}
\Omega\subset \left\{ \mu+\beta\ |\ \beta\in \Delta_0\right\}\cup U\cup
\left\{- \mu+\beta\ |\ \beta\in \Delta_0\right\}.
\eeq
\item[(PII)] For every $\alpha\in \Delta$ there is a $\mu_\alpha\in \Omega$ such that
\beq \label{mua}
\Omega_\alpha\subset \left\{ \mu_\alpha-\alpha+\beta\ |\ \beta\in \Delta_0\right\}\cup
\left\{ -\mu_\alpha+\beta\ |\ \beta\in \Delta_0\right\}.
\eeq
\end{description}
\es

\bprf
If $\lag$ is weak-Berger we have $\Gamma=\Delta_0$. We will use this property for $0\in\Delta_0$ as well as
for every $\alpha\in \Delta$.

{\bf (PI)} $\Gamma=\Delta_0$ implies that there is
$\phi\in \Pi$ and $\mu\in \Omega$ such that $0=\phi+\mu$ with  $Q\in {\cal B}_\phi$ and $
u\in V_\mu$ such that $0\not= Q(u)\in \mf{t}$, i.e. $\phi=-\mu \in \Pi$.
We fix such  $u, Q $ and $\mu$.
For  arbitrary $\lambda\in \Omega$ then occur the following cases:
\begin{description}
\item{{\em Case 1:}} There is a $v_+\in V_\lambda$ such that $Q(v_+)\not=0$ or a $v_-\in V_{-\lambda}$  such that
$Q(v_-)\not=0$. This implies $-\mu+\lambda \in \Delta_0$ or $-\mu-\lambda \in \Delta_0$, i.e.
$\lambda\in \left\{\mu+\beta\ |\ \beta\in\Delta_0\right\}\cup \left\{-\mu+\beta\ |\ \beta\in\Delta_0\right\}$.
\item{{\em Case 2:}} For all $v\in V_\lambda \+ V_{-\lambda}$ holds $Q(v)=0$. Then the Bianchi identity implies for
$v_+\in V_\lambda$ and $v_-\in V_{-\lambda}$
that $0=\lambda( Q(u))H(v_+,v_-).$
Now one can choose $v_+$ and $v_-$ such that $H(v_+,v_-)\not=0$. This implies 
$\lambda\in Q(u)^\bot=:U$ and we get (PI).
\end{description}

{\bf (PII)} Let $\alpha\in \Delta$. $\Gamma=\Delta_0$ implies the existence of
$\phi\in \Pi$ and $\mu_\alpha\in \Omega$ such that $\alpha=\phi+\mu_\alpha$ with  $Q\in {\cal B}_\phi$ and  $
u\in V_{\mu_\alpha}$ such that $0\not = Q(u)\in \lag_\alpha$. We fix $Q$ and $u$ for $\alpha$
and have that $\alpha-\mu_\alpha=\phi\in \Pi$ a weight of ${\cal B}_H$.

Let now $\lambda$ be a weight in $ \Omega_\alpha$, i.e. $\lambda+\alpha$ is also a weight. Thus
$-\lambda-\alpha$ is a weight. If $v\in V_\lambda$ then $Q(u)v\in V_{\lambda+\alpha}$. Since $H$ is non-degenerate
there is a $w\in
V_{-\lambda-\alpha}$ such that $H(Q(u)v,w )\not=0$. Since $Q\in {\cal B}_H(\lag)$ the Bianchi identity then gives
\[0=H(Q(u)v,w )+H(Q(v)w,u)+H(Q(w)u,v),\]
i.e. at least one of $Q(v)$ or $Q(w)$ has to be non-zero.
Hence we have two cases for $\lambda\in \Omega_\alpha$:
\begin{description}
\item{{\em Case 1:}} $Q(v)\not=0$.
This implies $-\mu_\alpha+\alpha+\lambda \in \Delta_0$, i.e.
$\lambda\in \left\{\mu_\alpha-\alpha+\beta\ |\ \beta\in\Delta_0\right\}$.
\item{{\em Case 2:}} $Q(w)\not=0$.
This implies $-\mu_\alpha+\alpha-\lambda-\alpha=-\mu_\alpha-\lambda \in \Delta_0$, i.e.
$\lambda\in \left\{-\mu_\alpha+\beta\ |\ \beta\in\Delta_0\right\}$.
\end{description}
But this is (PII).
\eprf

Of course it is desirable to find weights $\mu$ and $\mu_\alpha$ which are extremal in order to handle
criteria (PI) and (PII). To show in which sense this is possible we need a

\blem \label{nulllemma}
Let $\lag\subset \laso(V,H)$ an irreducible acting, complex semisimple Lie algebra.
For an extremal weight vector $u\in V_\Lambda$ there is a weight element $Q\in {\cal B}_H(\lag)$
such that $Q(u)\not=0$.
\elem

\bprf Let $u\in V_\Lambda$ be extremal with $Q(u)=0$ for every weight element $Q$.
Since $\bhg\ = \ \bigoplus_{\phi\in \Pi}{\cal B}_\phi$ the assumption implies 
$Q(u)=0$ for all $Q\in \bhg$. But this gives for every $A\in \lag$ and every weight element $Q$ that
\[Q(A u)\ = \ \left[A,Q(u)\right]- \underbrace{(A\cdot Q)}_{\in \bhg}  (u) \ =\ 0.\]
On the other hand $V$ is irreducible and thats why generated as vector space by elements of the form
$A_1\cdot \ldots \cdot A_k\cdot u$ with $A_i\in\lag$ and $k\in\mathbb{N}$ (see for example \cite{serre87}).
By successive application of $\lag$ to $u$ we get that $Q(v)=0$ for every weight element $Q$ and every weight
vector  $v$. But this gives $Q(v)=0$ for all $Q\in \bhg$ and every $v\in V$, hence $\bhg=0$.
\eprf

\bs\label{gewichtssatz}
Let $\lag$ be a semisimple
Lie algebra with roots $\Delta$ and $\Delta_0=\Delta\cup\{0\}$. Let $\lag\subset
\laso(V,H)$ irreducible, weak-Berger with weights $\Omega$.
Then there is an ordering of $\Delta$ such that the following holds:
If $\Lambda $ is the highest weight of
$\lag\subset\laso(V,H)$ with respect to that ordering, then the following properties are
satisfied:
\begin{description}
\item[(QI)] There is a $\delta\in \Delta_+\cup\{0\}$ and a hyperplane $U\subset \mf{t}^*$ such that
\beq \label{highest}
\Omega\subset \left\{ \Lambda-\delta+\beta\ |\ \beta\in \Delta_0\right\}\cup U\cup
\left\{ -\Lambda+\delta+\beta\ |\ \beta\in \Delta_0\right\}.
\eeq
\end{description}
If $\delta$ can not be chosen to be zero, then holds
\begin{description}
\item[(QII)] There is an $\alpha\in \Delta$  such that
\beq \label{qii}
\Omega_\alpha\subset \left\{ \Lambda-\alpha+\beta\ |\ \beta\in \Delta_0\right\}\cup
\left\{ -\Lambda+\beta\ |\ \beta\in \Delta_0\right\}.
\eeq
\end{description}
  \es

\bprf
First we consider the extremal weights of the representation, i.e. the images of the highest weight
under the Weyl group.
These do not lie in one hyper plane (because this would imply that all roots lie in one hyperplane).
Thus by proposition \ref{musatz} --- fixing $\mu\in \Omega$ --- there is an extremal weight $\Lambda$ with
$\Lambda+\mu\in \Delta_0$ or $\Lambda-\mu\in \Delta_0$. This one we fix.

Since the Weyl group acts transitively on the extremal weights we can find a fundamental root system,
i.e. an ordering on the roots, such that  $\Lambda$ is the highest weight. With respect to this
fundamental root system the roots splits into  positive and negative roots $\Delta=\Delta_+
\cup\Delta_-$. This implies
\beq\label{lambda-mu}
\mu=\varepsilon (\Lambda - \delta)
\eeq
with $\delta\in \Delta_+$ and
$\varepsilon=\pm 1$.

For an arbitrary $\lambda\in \Omega$ then holds $\lambda\in U=Q(u)^\bot$ or
$\lambda+\mu\in \Delta_0 $ or $\lambda-\mu\in \Delta_0$.
But with (\ref{lambda-mu}) this implies that we find an $\beta\in \Delta_0$
such that $\lambda = \pm(\Lambda-\delta) + \beta$ with $\beta \in \Delta_0$.
This is (QI). Note that we are still free to choose $\Lambda$ or $-\Lambda$ as highest weight.

Now we suppose that $\delta$ can not be chosen to be zero. Let $v\in V_\Lambda$ be a highest weight vector
or $v\in V_{-\Lambda}$.
Looking at the proof of proposition \ref{musatz}
one has that for all weight elements $Q\in {\cal B}_h(\lag)$
holds $Q(v)\in \lag_\alpha$ for a $\alpha\in \Delta$.
Since $\lag$ is weak-Berger $\bhg$ is non-zero. Thus we get
by lemma \ref{nulllemma}  that there is a weight element $Q$ such that $0\not=Q(v)\in \lag_\alpha$
and we are done (possibly by making $-\Lambda$ to the
highest weight).
\eprf

\paragraph{Representations of $\mf{sl}(2,\ccc)$}
To illustrate how these criteria shall work we apply them to irreducible representations of $\mf{sl}(2,\ccc)$.

\bs
Let $V$ be an irreducible, complex, orthogonal  $\mf{sl}(2,\ccc)$--module of highest weight $\Lambda$.
If it is weak-Berger then
$\Lambda \in \{2,4\}$.
\es
\bprf
Let $\mf{sl}(2,\ccc)\subset \laso (N, \ccc)$ be an irreducible representation of highest weight $\Lambda$.
I.e. $\Lambda(H)=l\in \mathbb{N}$ for $\mf{sl}(2,\ccc)=span (H,X,Y)$ where $X$ has the root $\alpha$.
Since the representation is orthogonal,
$l$ must be even 
and $0$ is a weight.
The hypersurface $U$ is the point
$0$.
Now property (\ref{mu}) ensures that $l\in \{2,4,6\}$. If $\mu=\Lambda$ we obtain $l\in\{2,4\}$.
If $\mu\not=\Lambda$ we can apply (QII): We have that $\Omega_\alpha=\Omega\setminus \{\Lambda\}$ and
$\Omega_{-\alpha}=\Omega\setminus \{-\Lambda\}$. Then (QII) implies that $l\in\{2,4\}$.
\eprf

So we get the first result.
\bfolg
Let $\mf{su}(2)\subset \laso(E,h)$ be a real irreducible weak-Berger algebra of real type. Then
it is a Berger algebra. In particular it is equivalent to the Riemannian holonomy representations of
$\laso(3,\rr)$ on $\rr^3$ or of the symmetric space of type $AI$, i.e.
$\mf{su}(3)/\laso(3,\rr)$ in the compact case or $\mf{sl}(3,\rr)/\laso(3,\rr)$ in the non-compact case.
\efolg

\subsection{Berger algebras, weak Berger algebras and spanning triples}

In this section we will describe a result of \cite{schwachhoefer1} and \cite{schwachhoefer2},
where holonomy groups of torsionfree connections, i.e. Berger algebras, are classified.
We will describe our results in their language such that we can use a partial result
of \cite{schwachhoefer2}.

For a Berger algebra holds that for every $\alpha \in \Delta_0$
there is a weight element $R\in {\cal K}(\lag)$ and weight vectors $u_1\in V_{\mu_1}$ and $u_2\in V_{\mu_2}$ such
that $0\not= R(u_1,u_2)\in \lag_\alpha$. The Bianchi identity then gives for an arbitrary $v\in V$
\[R(u_1,u_2)v\ =\ R(v,u_2)u_1 + R(u_1,v)u_2.\]
Choosing now $u_1,u_2$ such that $0\not=R(u_1,u_2)\in \mf{t}$ one gets for any $\lambda\in \Omega$ and
$v\in V_\lambda$ that
\[\lambda(R(u_1,u_2))v\ =\ R(v,u_2)u_1 + R(u_1,v)u_2.\]
This implies  $\lambda\in \left(R(u_1,u_2)\right)^\bot\subset \mf{t}^*$ or $V_\lambda\subset
\lag V_{\mu_1}\+ \lag V_{\mu_2}$.
This gives property
\begin{description}
\item[(RI)] There are weights $\mu_1,\mu_2\in \Omega$ such that
\[\Omega\subset \{\mu_1+\beta\ |\ \beta\in \Delta_0\}\cup U\cup \{\mu_2+\beta\ |\ \beta\in \Delta_0\}.\]
\end{description}
If one chooses $u_1,u_2$ such that $0\not=R(u_1,u_2)=A_\alpha\in \lag_\alpha$ with $\alpha\in \Delta$ then one gets
for  $\lambda\in \Omega$ that $A_\alpha V_\lambda\subset \lag V_{\mu_1}\+\lag V_{\mu_2}$. This means that
the weights of $A_\alpha V_\lambda$  are contained in $\{\mu_1+\beta|\beta\in \Delta_0\}\cup
\{\mu_2+\beta|\beta\in \Delta_0\}$. But this is property
\begin{description}
\item[(RII)] For every $\alpha\in \Delta$ there are weights $\mu_1,\mu_2\in \Omega$ such that
\[\Omega_\alpha\subset \{\mu_1-\alpha+\beta\ |\ \beta\in \Delta_0\}
\cup \{\mu_2-\alpha+\beta\ |\ \beta\in \Delta_0\}.\]
\end{description}
Of course our (PI) is a special case of (RI) with $\mu_1=-\mu_2$. (PII) is not a special case of (RII)
since $\mu_\ä+\ä$ is not a weight apriori.

To describe this situation further in \cite{schwachhoefer2} the following definitions are made.
We point out that here $\Omega_\alpha$ does not denote the weights of $\lag_\alpha V$ but the
weights $\lambda$ of $V$ such that $\lambda+\alpha$ is a weight.

\bde
Let $\lag\subset End(V) $ be an irreducible acting complex Lie algebra, $\Delta_0$ be the roots
and zero of the semisimple part of $\lag$, $\Omega$ the weights of $\lag$ and $\Omega_\alpha$ as above.
\bnum
\item A triple $(\mu_1,\mu_2,\alpha)\in \Omega\times\Omega\times\Delta$ is called {\bf spanning triple} if
\[\Omega_\alpha \subset \left\{ \mu_1-\alpha+\beta\ |\ \beta\in \Delta_0\right\}\cup
\left\{ \mu_2-\alpha+\beta\ |\ \beta\in \Delta_0\right\}.\]
\item A spanning triple $(\mu_1,\mu_2,\alpha)$ is called {\bf extremal} if $\mu_1$ and $\mu_2$ are
extremal.
\item A triple $(\mu_1,\mu_2, U)$ with $\mu_1,\mu_2$ extremal weights and  $U$ an affine hyperplane
in $\mf{t}^*$
is called {\bf planar spanning triple} if every extremal weight different from $\mu_1$ and $\mu_2$ is
contained in $U$ and
$\Omega\subset
 \left\{ \mu_1+\beta\ |\ \beta\in \Delta_0\right\}\cup U
 \cup
\left\{ \mu_2+\beta\ |\ \beta\in \Delta_0\right\}$.
\enum
\ede

From (RI) and (RII) in \cite{schwachhoefer2} the following proposition is deduced.
\bs \cite{schwachhoefer2}
Let $\lag\subset End(V)$ be an irreducible complex Berger algebra.
Then for every root $\alpha\in \Delta$ there is a spanning triple. Furthermore there is an
extremal spanning triple or a planar spanning triple.
\es

If we return to the weak-Berger case we can reformulate proposition \ref{gewichtssatz} as follows.
\bs\label{triplesatz}
Let $\lag\subset \laso(V,H)$ be an irreducible complex weak-Berger algebra. Then there is an extremal weight
$\Lambda$ such that
one of the following properties is satisfied:
\begin{description}
\item[(SI)] There is a planar spanning triple of the form $(\Lambda,-\Lambda, U)$.
\item[(SII)] There is an $\ä\in \Delta$ such that
$
\Omega_\alpha\subset \left\{ \Lambda-\alpha+\beta\ |\ \beta\in \Delta_0\right\}\cup
\left\{ -\Lambda+\beta\ |\ \beta\in \Delta_0\right\}.
$
\end{description}
There is a fundamental system such that the extremal weight in (SI) and (SII) is the highest weight.
\es

\bprf The proof is analogous the the one of proposition \ref{gewichtssatz}. If there is an $\alpha\in \Delta$
such that the corresponding $\mu_\alpha$ is extremal we are done.
If not, then for every extremal weight vector $u\in V_\Lambda$
and every weight element $Q\in {\cal B}_\phi$ holds that $Q(u)\in\mf{t}^*$. Then by lemma \ref{nulllemma}
there is a $Q$ such that $0\not=Q(u)\in \mf{t}^*$.
As before this implies
\[\Omega\subset \{\Lambda+\beta\ |\ \beta\in \Delta_0\}\cup U\cup \{-\Lambda+\beta\ |\ \beta\in \Delta_0\}.\]
To ensure that $(\Lambda,-\Lambda,U)$ is a planar spanning triple we have to show that every
extremal weight $\lambda$ different from $\Lambda$ and $-\Lambda$ is contained in $U=Q(u)^\bot$.
Let $\lambda$ be extremal and different from $\Lambda$ and $- \Lambda$,
$v_\pm\in V_{\pm\lambda}$ and $u\in V_\Lambda$. Since $Q(v_\pm)\in \mf{t}$ the Bianchi identity gives
\be
0 & =& H(Q(u)v_+,v_- )+H(Q(v_+)v_-,u)+H(Q(v_-)u,v_+)
\\&=&\lambda\left(Q(u)\right)\underbrace{H(v_+,v_- )}_{\not=0}-
\underbrace{\lambda\left(Q(v_+)\right)H(v_-,u)+\Lambda\left((Q(v_-)\right)H(u,v_+)}_{\mbox{$=0$ since $u$ is neither in
$V_\lambda$ nor in $V_{-\lambda}$}}.
\ee
Hence $\lambda\in U$.
\eprf
Obviously we are in a slightly different situation as in the Berger case since
$- \Lambda+\alpha$ is not necessarily a weight and in case it is a weight, it is not extremal in general.

\subsection{Properties of root systems}

In this section we will recall the properties of abstract root systems. Let $(E,\la.,.\ra)$ be a euclidian
vector space. A finite set of vectors $\Delta$ is called root system if it satisfies the following properties
\bnum
\item $\Delta$ spans $E$.
\item For every $\alpha \in \Delta$ the reflection on the hyperplane perpendicular to $\alpha$ defined by
\[s_\alpha (\vf):= \vf - \frac{2\la\vf, \alpha\ra}{\|\alpha\|^2}\alpha\] maps $\Delta$ onto itself.
\item For $\alpha,\beta\in \Delta$ the number $\frac{2\la\beta, \alpha\ra}{\|\alpha\|^2}$ is an integer.
\enum
A root system is called indecomposable if it does not split into orthogonal subsets. It is called reduced
if $2\alpha$ is not a root if $\alpha$ is a root.

The indecomposable, reduced root systems corresponds to the roots of simple Lie algebras.
They are classified in a finite
list: $A_n,\  B_n,\  C_n,\  D_n,\ E_6,\ E_7,\ E_8,\ F_4$ and $G_2$. The index designates the dimension of $E$.

We will cite some basic properties of root systems, which can be found for example in \cite{knapp96}.

\bs(See for example \cite{knapp96}, pp. 149) Let $\Delta$
be  an abstract reduced root system in $(E,\la.,.\ra)$.\label{knappsatz}
\bnum
\item \label{ks1} If $\alpha\in \Delta$, then the only root which is proportional to $\alpha$ is $-\alpha$.
\item \label{ks2} If $\alpha,\beta\in \Delta$, then $\frac{2\la\beta, \alpha\ra}{\|\alpha\|^2}\in \{0,\pm1,\pm 2,\pm 3\}$.
If $\Delta$ is one of the indecomposable root systems
$\pm 3$ occurs only for the root system $G_2$. If both roots are non proportional then $\pm 2$
only occurs for $B_n, C_n, F_4$ or $G_2$.
\item \label{ks3} If $\alpha$ and $\beta$ are nonproportional in $\Delta$ and $\|\beta\|\le\|\alpha\|$, then
$\frac{2\la\beta, \alpha\ra}{\|\alpha\|^2}\in \{0,\pm1\}$.
\item \label{ks4} Let be $\alpha,\beta\in \Delta$. If $\la\alpha,\beta\ra>0$, then $\alpha-\beta\in \Delta$.
If $\la\alpha,\beta\ra<0$, then $\alpha+\beta\in \Delta$. I.e. if neither $\alpha-\beta\in \Delta$ nor
$\alpha+\beta\in \Delta$, then $\la \alpha,\beta\ra=0$.
\item \label{ks5} The subset of $\Delta$ defined by $\left\{\beta+ k\alpha\in \Delta\cup\{0\}|k\in \mathbb{Z}\right\}$
is called
$\alpha$--string through $\beta$. It has no gaps, i.e. $\beta + k\alpha \in \Delta$ for $
-p\le k\le q$ with $p,q\ge 0$ and it holds $p-q=\frac{2\la\beta, \alpha\ra}{\|\alpha\|^2}$.
The maximal length of such string is given by $\max_{\alpha,\beta\in\Delta}
\frac{2\la\beta, \alpha\ra}{\|\alpha\|^2}+1$, i.e.
it contains at most four roots.
 \enum
 \es
As a consequence of that proposition we get the following lemmata.
In these we will refer to long and short roots. This notion is evident because in the indecomposable
reduced root systems of type $B_n,\ C_n,\ F_4$ and $G_2$ the roots have two different lengths.

\blem
\label{wurzellemma0}
Let $\D$ be an indecomposable, reduced root system. Then it holds:
\bnum
\item
If $a\ä+\bb\in \D$ for $a\in \mathbb{N}$ and $a>1$, then $\la\ä,\bb\ra<0$ and $\ä$ is a short root.
\item If $\Delta$ is a root system, where the roots have equal length or if $\ä$ is a long root,
then
$\ä+\bb\in\D$ implies $\la\ä,\bb\ra<0$.
\item
Let $\ä$ and $\bb$ be two short roots. If $\ä+\bb$ is a long root then $\la\ä,\bb\ra=0$,
if it is a short one then $\la\ä,\bb\ra<0$.
The sum of a short and a long root is a short one
\item If $\bb$ is a long root in $\D\not=G_2$, then there are orthogonal roots $\ä$ and $\gamma$ such
that
$\bb=\ä+\gamma$.
\enum
\elem
\bprf The proof follows directly from proposition \ref{knappsatz}.\eprf

\blem\label{wurzellemma1}
Let $\alpha$ and $\beta$ be two nonproportional roots and $a\in \mathbb{N}$. If $a(\alpha+\beta)\in \Delta$
then $a=1$.
\elem

\bprf
If $a>1$ then $\alpha+\beta$ is not a root. This implies $\la\alpha,\beta\ra\ge 0$
and yields for $a(\alpha+\beta)=\gamma\in \Delta$:
\barr{rcccl}
0&<& a\left(\|\alpha\|^2 +\la \alpha,\beta\ra \right)&=&\la\alpha,\gamma\ra\\
0&<& a\left(\la \alpha,\beta\ra+\|\beta\|^2  \right)&=&\la\gamma,\beta\ra.
\earr
On the other hand we have
\[\|\gamma\|^2=a\left(\la \alpha,\gamma\ra + \la \beta,\gamma\ra \right).\]
But this gives
\[1=\frac{a}{2} \Big(\underbrace{
\underbrace{\frac{2\la\gamma, \alpha\ra}{\|\gamma\|^2}}_{\mbox{$>0$ in $\mathbb{N}$}}
+ \underbrace{\frac{2\la\gamma, \beta\ra}{\|\gamma\|^2}}_{\mbox{$>0$ in $\mathbb{N}$ }}}_{\mbox{$\ge 2$
 in $\mathbb{N}$}}\Big)
.\]
This is a contradiction. Hence $a=1$.
\eprf
The next lemma is a little more general.

\blem\label{wurzellemma2}
Let $\alpha$ and $\beta$ be two non-proportional roots in an indecomposable root system
and $a,b\in \mathbb{N}$ with $a\le b$ such that
$a\alpha+b\beta\in \Delta$.
\bnum
\item
If $\Delta$ is not $G_2$ then $a=1$. If $\Delta= A_n, D_n, E_6,E_7, E_8$ then $b=1$ too.
If $\Delta=B_n, C_n , F_4$ then $b\le 2$.
\item If $\Delta=G_2$ then $a\le 2$ and $b\le 3$.
\enum
\elem

\bprf We suppose $a\alpha + b \beta =\gamma\in \Delta$.

First we consider the case $\la \alpha, \beta\ra\ge 0$. This gives
\barr{rcccl}
0&<& a\|\alpha\|^2 +b\la \alpha,\beta\ra &=&\la\alpha,\gamma\ra\\
0&<& a(\la \alpha,\beta\ra+b\|\beta\|^2  &=&\la\gamma,\beta\ra.
\earr
On the other hand we have $\|\gamma\|^2=a\la \alpha,\gamma\ra + b\la \beta,\gamma\ra $ and thus
\be
1&=&\frac{a}{2}\overbrace{\frac{2\la\gamma, \alpha\ra}{\|\gamma\|^2}}^{>0}
+\frac{b}{2}\overbrace{\frac{2\la\gamma, \beta\ra}{\|\gamma\|^2}}^{>0}
\\&\ge&
\frac{a}{2} \Big(\underbrace{
\frac{2\la\gamma, \alpha\ra}{\|\gamma\|^2}
+ \frac{2\la\gamma, \beta\ra}{\|\gamma\|^2}}_{\mbox{$\ge 2$
 in $\mathbb{N}$}}\Big)
.\ee
Hence $a=1$.

Let now be $\la\alpha,\beta\ra<0$. This implies, that $\alpha+\beta=:\delta$ is a root with the property
$\delta-\beta=\alpha\in \Delta$

Although the above proposition does not assert that this implies $\la \delta,\beta\ra\ge 0$ we can show this.
Suppose that $\la\delta,\beta\ra<0$. Hence $\delta+\beta=\alpha+2\beta$ is a root. If we exclude the root system
$G_2$ point \ref{ks5} of proposition \ref{knappsatz}  implies
$\frac{2\la\alpha, \beta\ra}{\|\beta\|^2}=-2$, i.e. $\la\alpha,\beta\ra=-\|\beta\|^2$ and finally
$\la\delta,\beta\ra=0$, which was excluded.

Thus we have that $\la\delta,\beta\ra\ge 0$. Analogously to the first case we get
\[
1=\frac{a}{2}\ \frac{2\la\gamma, \delta\ra}{\|\gamma\|^2}
+\frac{(b-a)}{2}\ \frac{2\la\gamma, \beta\ra}{\|\gamma\|^2}.
\]
In case that $a\le b-a$ we get again that $a=1$.
Otherwise we get $b-a=1$, i.e. $a\delta+\beta=\gamma$. Again by
point \ref{ks5} of proposition \ref{knappsatz} we get
$p-q=\frac{2\la\gamma, \delta\ra}{\|\delta\|^2}\ge 0$. But this implies $a\le 1$.

The possible values for $b$ follow also by proposition \ref{knappsatz}.

For $G_2$ the possible values of $a$ and $b$ can be calculated analogously.
\eprf

\blem\label{wurzellemma3}
Let $\eta$ be a long root of an indecomposable root system.
\bnum
\item Let $a,b\in\mathbb{N}$ and $\alpha\in \Delta$ not proportional to
$\eta$ such that $a\eta+b\alpha\in \Delta$. Then $a\le b$, i.e. $a=1$ if $\Delta$ not equal to $G_2$ and
$a\le 2$ otherwise.
\item
Let $\alpha,\beta$ in $\Delta$ not proportional to $\eta$ and $a\in \mathbb{N}$ such  that
$a\eta +\alpha +\beta\in \Delta$. Then $a\le2$.
\enum
\elem

\bprf
1.) First we exclude $G_2$ and suppose  that $b=1$, i.e.
$a\eta+\alpha=\gamma\in¸\Delta$. Hence $-p\le a\le q$ and
\[
\left| p-q \right|=\frac{2|\la\eta, \alpha\ra|}{\|\eta\|^2} < 2\frac{\|\eta\|\cdot \|\eta\|}{\|\eta\|^2}
\le 2,
\]
i.e. $|p-q|\le 1$. But since we have excluded $G_2$ we have that $a=1$.

For $G_2$ a long root $\eta$ is given by $2e_3-e_1-e_2$ with the notations of appendix C of \cite{knapp96}.
For this we get the wanted result.

2.)
Let $a\eta+ \alpha +\beta =\gamma$.

First we consider the case that $\alpha+\beta$ or $\alpha-\gamma$ or $\beta-\gamma$ is a root.
If this root is not proportional to $\eta$ we have by the first point that $a\le 1$.
If it is proportional to $\eta$ we get that $a\le 2$ and we are done.

Now we suppose that neither $\alpha+\beta$ nor $\alpha-\gamma$ nor $\beta-\gamma$ is a root. This
implies $\la \alpha,\beta\ra\ge 0$, $\la \alpha,\gamma\ra\le 0$ and $\la \beta,\gamma\ra\le 0$.
We consider the equations
\barr{rcccccccl}
a\la\eta,\alpha\ra &+& \|\alpha\|^2 &+& \underbrace{\la\alpha,\beta\ra}_{\ge 0}&=&\la\gamma,\alpha\ra&\le&0\\
a\la\eta,\beta\ra &+& \|\beta\|^2 &+& \underbrace{\la\alpha,\beta\ra}_{\ge 0}&=&\la\gamma,\beta\ra&\le&0\\
a\la\eta,\gamma\ra &+& \underbrace{\la\alpha,\gamma\ra}_{\le 0}
&+& \underbrace{\la\beta,\gamma\ra}_{\le 0}&=&\|\gamma\|^2&>&0.
\earr
Hence we have that $\la\eta,\alpha\ra<0$, $\la\eta,\beta\ra<0$ and $\la\eta,\gamma\ra>0$. But since $\eta$ is
long, not proportional neither to $\alpha$ nor to $\beta$ we have that
\be
\|\eta\|^2\ \ge\  \la \gamma,\eta\ra
&=&
 a \|\eta\|^2 +\la \alpha,\eta\ra+\la\beta,\eta\ra\\
 &=&
 a \|\eta\|^2 - \underbrace{|\la \alpha,\eta\ra|}_{<\|\alpha\|\cdot\|\eta\|\le \|\eta\|^2}
-\underbrace{|\la\beta,\eta\ra|}_{<\|\beta\|\cdot\|\eta\|\le \|\eta\|^2}\\
&>&(a-2)\|\eta\|^2.
\ee
This gives $a-2<1$ which is the proposition.
\eprf

\blem\label{wurzellemma5}
Let $\alpha$ be a long root and $\eta$ be a short one with $\la\aaa,\eta\ra>0$, i.e.
$\frac{2\la\alpha,\eta\ra}{\|\eta\|^2}\ge 2$.
Then there is a short root $\beta$ with $\bb\not\sim\eta$,
$ \la\beta,\aaa\ra <0$ and $\la\bb,\eta\ra\le 0$. If the rank of the root system is
greater than $2$ or if $\frac{2\la\alpha,\eta\ra}{\|\eta\|^2}= 3$ (which can only occur for $G_2$),
$\bb$ can be chosen such that $\la\bb,\eta\ra<0$.
\elem
\bprf
$\la\aaa,\eta\ra>0$ implies that $\eta-\alpha $ is a root, in particular a short one.
For the inner product we get
\[
\la\aaa,\eta-\alpha\ra=\la\aaa,\eta\ra-\|\aaa\|^2<\|\aaa\|\|\eta\|-\|\aaa\|^2<\|\aaa\|^2-\|\aaa\|^2=0
\]
and
\[
\la\eta,\eta-\alpha\ra=\|\eta\|^2-\la\aaa,\eta\ra\le 0.\]
In case of $\frac{2\la\alpha,\eta\ra}{\|\eta\|^2}= 3$ the last $\le$ is a $<$,
and we are done with the second point in case of $G_2$.

If the rank of the root system is greater than $2$ this can be seen with the help of the
definitions of the reduced, indecomposable  root systems (see appendix C of \cite{knapp96}).
\eprf

\section{Simple weak-Berger algebras of real type}

In this section we will apply the result of proposition
\ref{triplesatz} to simple complex irreducible acting Lie algebras.

We will do this step by step under the following special conditions:
\bnum
\item The highest weight of the representation is a root.
\item The representation satisfies (SI), i.e. admits a planar spanning triple $(\LL,-\LL,U)$.
\item The representation satisfies (SII) and has weight zero.
\item The representation satisfies (SII) and does not have weight zero.
\enum
 Throughout this section
the considered Lie algebra is supposed to be different from $\mf{sl}(2,\ccc)$.

\subsection{Representations with roots as highest weight}

\bs
\label{wurzelgewicht}
Let $\lag \subset \laso(N, \ccc ) $
be an irreducible representation of real type of a  complex simple Lie algebra different from
$\mf{sl}(2,\ccc)$ and satisfying  (SI) or (SII).
If we suppose in addition that there is an extremal weight $\Lambda$ with
$\Lambda = a \eta$ for a root $\eta\in \Delta$ and $a>0$, then holds the following:
\bnum
\item If $\eta $ is a long root, then $a=1$ and the representation is the adjoint one.
\item If $\eta $ is a short root, then holds the following for $a$:
\bnum
\item If $\D=B_n$ or $G_2$ then $a=1,2$.
\item If $\D=C_n$ or $F_4$ then $a=1$.
\enum
\enum
\es

\bprf
Let $\Lambda = a\eta$ with $ \eta\in \Delta$, $a\in \mathbb{N}$. W.l.o.g. we may suppose that
$\Lambda $ is the extremal weight in the properties (SI) and (SII).
(If not then there is an element of the Weyl group $\sigma$
mapping $\Lambda$ to the extremal weight of (SI) and (SII) $\Lambda^\prime$.
Then $\Lambda^\prime = a \sigma \eta$ and $\sigma \eta \in \Delta$.)

First we show that $a\in \mathbb{N}$. If we chose an fundamental system $(\pi_1,\ldots, \pi_n)$ such that
$\Lambda=a\eta$ is the highest weight we get that $\la\LL,\pi_i\ra=a\la\eta,\pi_i\ra\in \mathbb{N}$ for all $i$.
$a\not\in\mathbb{N}$ would imply that $\la\eta,\pi_i\ra\ge 2$ for all $i$ with
$\la\eta,\pi_i\ra\not=0$. This holds only for the root system $C_n$ where $\Lambda=\w_1= \einhalb \eta$.
But this representation is symplectic but not orthogonal. (For an explicit
formulation of this criterion see \cite{onishchik-vinberg3}.) So we get $a\in \mathbb{N}$.

Now we consider two cases.
\begin{description}
\item{{\em Case 1: $\eta$ is a long root:}} In this case the root system of long roots, denoted by $\Delta_l$
is the orbit of $\eta$ under the Weyl group. Hence $a\cdot \D_l$ are the extremal weights and $\D\subset\Omega$.
This implies $0\in \Omega_\alpha$ for every $\alpha\in \Delta$.

Furthermore  for all roots holds that $a\cdot \Delta\subset \Omega$. This is  true because we can find a
short root such that $\la\eta,\beta\ra>0$. This implies $\eta-\beta\in \D_s$. On the other hand it is
$\frac{2\la a\eta,\beta\ra}{\|\beta\|^2}\ge a$. Hence $a \eta-a\beta=a(\eta-\beta)\in \Omega$. Applying the
Weyl group to this weight we get the property for all short roots.

\begin{description}
\item{(SI)} Let $\Lambda$ satisfy (SI), i.e. $\Lambda$ and $-\Lambda$ define a planar spanning triple
$(\Lambda,-\Lambda,U)$.
This would imply that every long root different from $\eta$ lies in the hyperplane $U$. This is only possible
for the the root system $C_n$, because all other root systems have an indecomposable system of long roots.
For $C_n$ holds that $\D_l=A_1\times\ldots\times A_1$.
But we have still a root $\beta$ --- possibly a short one --- such that
$\beta\not\in U$ and $\beta$ not proportional to $\eta$.
This implies  $\Omega\ni a\beta=\Lambda+\gamma=a\eta+\gamma$ or
$\Omega\ni a\beta=-\Lambda+\gamma=-a\eta+\gamma$ with $\gamma\in \Delta_0$. Then Lemma \ref{wurzellemma1}
implies $a=1$.
\item{(SII)}
Lets suppose that $\Lambda$ satisfies (SII), i.e. there is an $\alpha\in \Delta$ such that
$\Omega_\alpha\subset \{\Lambda-\alpha+\beta|\beta\in \Delta_0\}\cup\{-\Lambda+\beta|\beta\in \Delta_0\}$.
$0\in \Omega_\alpha$ implies $0=\Lambda-\alpha+\beta=a\eta-\alpha+\beta$ or
$0=  -\Lambda+\beta=-a\eta+\beta$ with $\beta\in \Delta_0$. The second is not possible and the first implies by
lemma \ref{wurzellemma3} that $a=1$ or $a=2$ and $\eta=\ä$. In the second case we find a root
$\gamma\not\sim\ä$ such that $\la\gamma,\ä\ra<0$, hence $2\gamma\in\W_\ä$. Since $2\gamma-2\ä\not\in\D$
it has to be $2\gamma=\ä+\bb$, but this is prevented by $\la\gamma,\ä\ra<0$ and lemma \ref{wurzellemma0}.

\end{description}
Of course if $\eta $ is a long root the representation is the adjoint one.
\item{{\em Case 2: $\eta$ is a short root:}}
Lets denote by $\D_s$ the root system of short roots. It equals to the orbit of $\eta $ under the Weyl group.
It is a root system of the same rank as $\D$ and all roots have the same length.
Clearly  $\D_s\subset \Omega$ and $a\cdot\D_s$ are the extremal weights in $\Omega$.
 For the root system  $B_n$ the root system of short roots
$\D_s$ equals to $A_1\times \ldots \times A_1$, otherwise it is indecomposable.

Furthermore holds the following: If  $a\ge2$ then $\D\subset \Omega$. To verify this, we consider a long root
$\beta\in \Delta_l$ with the property that $\la\bb,\eta\ra >0$. Such a $ \bb$ always exists.
Then we have
$\frac{2\la \eta,\beta\ra}{\|\eta\|^2}> \frac{2\la \eta,\beta\ra}{\|\beta\|^2}\ge1$. This implies
$2\eta-\beta\in \D$ (see proposition \ref{knappsatz}).  On the other hand $a\ge 2$ ensures
that $\Omega\ni s_\bb(2\eta)=2\left(\eta-\frac{2\la \eta,\beta\ra}{\|\beta\|^2}\bb\right)$. This implies
that the long root $2\eta-\bb $ is a weight. Now applying the Weyl group to $\bb$
shows that every long root is a weight.

\begin{description}
\item{(SI)}
We suppose that there is a planar spanning triple $(\LL,-\LL,U)$.
This implies that $a\beta$ lies in the hyperplane $U$ if $\beta$ is a short root.
But this is only possible for $B_n$ because the short roots of all other root systems
are  indecomposable.

In case of $B_n$ we can at least find a long root $\alpha$
which is not in $U$. Since the long roots are weights,
we have $\alpha=a\eta +\gamma$ or $\alpha=-a\eta+\gamma$ with $\gamma\in \D_0$. But this implies
for $B_n$ that $a\le 2$.

\item{(SII)}
Suppose that there is an $\ä\in \Delta$ such that
$
\Omega_\alpha\subset \left\{ \Lambda-\alpha+\beta\ |\ \beta\in \Delta_0\right\}\cup
\left\{ -\Lambda+\beta\ |\ \beta\in \Delta_0\right\}$.
 $\D\subset \Omega$
implies  $0\in \Omega_\alpha$ for all $\aaa$.
$0=-a\eta+\gamma$ with $\gamma\in \D_0$ implies $a=1$. Hence if we suppose $a\ge 2$ we must have
\beq
\label{nulleqn}
0=a\eta-\alpha+\gamma
\eeq
Thus we have to deal with the following cases:
\bnum
\item[(a)] $\alpha=\eta$ and $a=2$.
\item[(b)] $\alpha\not\sim\eta$ and by \ref{ks5} of proposition \ref{knappsatz}
$a\le \frac{2\la \eta,\alpha\ra}{\|\eta\|^2}\le 3$. I.e. if $a\ge2$, $\alpha$ is a long root.
\enum
We exclude the first case for any root system different from $B_n$. Set $a=2$ and $\alpha=\eta$.
If $\D\not=B_n$ the short roots are indecomposable, i.e. there is a short root $\bb$ such that
$\bb\not\sim \eta$ and $\la\bb,\eta\ra<0$. Hence $2\bb\in \Omega_{\eta}$ and $\bb+\eta\in \Delta$.

The existence of a spanning triple implies then $2\beta=\eta+\gamma$
or $2\beta=-2\eta+\gamma$ with
$\gamma\in \Delta_0$. The second case is impossible because of lemma \ref{wurzellemma1}. The first implies
$2\beta-\eta\in\Delta$. Again this is not possible by \ref{wurzellemma0} and $\la\bb,\eta\ra<0$.
Hence the case (a) is excluded.

Now  we consider the case (b). First we show that $a=3$ is not possible. Set $a=3$.
We notice that
$\la\eta,\aaa\ra>0 $ implies  $\frac{2\la \eta,\alpha\ra}{\|\aaa\|^2}\ge 1$ and hence
$3\eta-3\aaa\in\Omega_\alpha$. Thus we have the alternative
$3\eta-3\aaa=3\eta-\aaa+\gamma$ or $3\eta-3\aaa=-3\eta+\gamma$ with $\gamma\in \D_0$. The first implies
$2\alpha\in \Delta$ and the second $6\eta-3\aaa\in \Delta$. Both are not true, hence $a=3$ is impossible.

We continue with case (b) and have that
$\alpha $ is a long root with
\[
\frac{2\la \eta,\alpha\ra}{\|\eta\|^2}\ge 2\ ,
\;\;\mbox{
i.e. $2\eta-\alpha\in\D$.}\]
From now on we suppose, that the root system is different from $G_2$. Then we have
\beq\label{alphaeta}
\frac{2\la \eta,\alpha\ra}{\|\eta\|^2}= 2.
\eeq
In a next step we will show that under these conditions
there is no short root $\beta$ with
\beq
\label{f4}
\bb\in\Delta_s\mbox{ with }\la\aaa,\beta\ra<0\ ,\; \la\aaa,\eta\ra<0 \mbox{ and }\; \bb\not\sim\eta.
\eeq
Suppose that there is such a $\beta$. Then the first condition implies
that $2\beta\in \Omega_\alpha$ and hence
$2\bb=2\eta-\aaa +\gamma$ or $2\bb=-2\eta+\gamma$ with $\gamma\in\Delta_0$. The latter is not possible.
The second implies the following using (\ref{alphaeta}):
\[
-2\ \ge\   2\cdot\frac{2\la \beta,\eta\ra}{\|\eta\|^2}
\ =\  \frac{2\la2\eta-\aaa,\eta\ra}{\|\eta\|^2}+  \frac{2\la \gamma,\eta\ra}{\|\eta\|^2}
\ =\  2+ \frac{2\la \gamma,\eta\ra}{\|\eta\|^2}.
\]
Hence $-4\ge \frac{2\la \gamma,\eta\ra}{\|\eta\|^2}$ which is impossible.

Now by the lemma \ref{wurzellemma5} there is such a $\beta$. Hence for any remaining root systems different from
$G_2$ and different from $B_n$ we have that $a=1$.
\end{description}
\end{description}
All in all we have shown, that for a long root holds $a=1$
and for a short root $a=2$ implies $\Delta=B_n$ or $G_2$.
\eprf

\bfolg\label{wurzelgewichtfolge}
Let $\lag \subset \laso(N, \ccc ) $
be an irreducible complex simple weak Berger algebra different from
$\mf{sl}(2,\ccc)$ and with the additional property that the highest weight
is of the form $\Lambda = a \eta$ for a root $\eta\in \Delta$. Then $\lag$ is the complexification
of a holonomy algebra of a
Riemannian manifold or the representation of $G_2$ with highest weight $2\w_1$.
\efolg

\bprf
Clearly if $\eta$ is a long root the representation is the adjoint one, i.e. the complexification of a
holonomy representation of
a Lie group with positive definite bi-invariant metric.
For a short root $\eta$ we get the following:
\begin{description}
\item[$B_n,\ a=1:$] This is the  representation of highest weight $\w_1$, i.e. the standard representation
of $\laso(2n+1,\ccc)$ on $\ccc^{2n+1}$. Of course this is
the complexification of the generic Riemannian holonomy representation.
\item[$B_n,\  a=2:$] This is the  representation of highest weight $2\w_1$. A further analysis shows that
this is the complexified representation of the Riemannian symmetric space of type $AI$, i.e.
of the symmetric spaces $SU(2n+1)/SO(2n+1,\rr)$, respectively $SL(2n+1,\rr)/SO(2n+1,\rr)$.
\item[$C_n,\  a=1:$] (for $n\ge 3$) This is the  representation of highest weight $\w_2$.
It is the complexified representation of the Riemannian symmetric space of type $AII$, i.e.
of the symmetric spaces $SU(2n)/Sp(n)$, respectively $SL(2n,\rr)/Sp(n)$.
\item[$F_4,\  a=1:$]  This is the  representation of highest weight $\w_1$.
It is the complexified representation of the Riemannian symmetric space of type $EIV$, i.e.
of the symmetric spaces $E_6/F_4$, respectively $E_{6(-26)}/F_4$.
\item[$G_2,\  a=1:$] This is the  representation of highest weight $\w_1$. It is the representation
of $G_2$ on $\ccc^7$, i.e. the complexification of the holonomy of a Riemannian $G_2$--manifold.
\item[$G_2,\  a=2:$] This is the representation $2\w_1$ of $G_2$.
It  is a $27$-dimensional representation of $G_2$ isomorphic to
$Sym_0^2\ccc^7$, where $\ccc^7$ denotes the standard module of $G_2$ and $Sym^2_0\ccc^7$ its
symmetric, trace free $(2,0)$--tensors.
This is the exception, because there is no Riemannian manifold with
this complexified holonomy  representation.
\end{description}
\eprf

\subsection{Representations with planar spanning triples}

Now we consider representations of a simple Lie algebra under the condition that there is a
planar spanning triple. The proof of this proposition is a copy of the proof in \cite{schwachhoefer2}
adding the additional properties of our planar spanning triple.

\bs\label{planarsatz}
Let $\lag \subset \laso(N, \ccc ) $
be an irreducible representation of real type of a  complex simple Lie algebra different from
$\mf{sl}(2,\ccc)$ and satisfying (SI), i.e. with a planar spanning triple
of the form $(\Lambda,-\Lambda,U)$. If there is no root $\alpha$ such that
$\Lambda = a \alpha$ then
$\lag$ is of type $D_n$ with $n\ge3$  and the representation is congruent to the one with highest weight
$\w_1$ or $2\w_1$.
\es

\bprf
The condition $\Lambda \not= a \alpha$ implies that there is no root such that
$-\LL=s_\alpha (\LL)$. The existence of a planar spanning triple then gives that
for any $\alpha\in \D$ with $\la\LL,\alpha\ra\not=0$ the image of the reflection
lies in $U$. If we set $U=T^\bot$ this gives
\beq\label{teqn}
\mbox{For $\alpha\in \Delta$ with $\la\alpha,\LL\ra\not=0$ holds }\;\;\;
\la \alpha,T\ra\ =\  \frac{\|\alpha\|^2}{2\la\LL,\alpha\ra}\la\LL,T\ra\ \not=\ 0.
\eeq
In the following we prove various claims to get the wanted result. We follow completely
the lines of reasoning in \cite{schwachhoefer2}.
\begin{description}
\item{Claim 1:}
{\em For any non-proportional $ \alpha,\beta\in \Delta$ with $\la\LL,\alpha\ra\not=0$ and
$\la\LL,\bb\ra\not=0$ holds that $\la \alpha,\bb\ra=0 $ or both have the same length.}

To show this we prove that for two such roots hold that they are orthogonal or that
$\la\LL,s_\aaa\bb\ra=\la\LL,s_\bb\aaa\ra=0$.

Suppose that $\la\LL,s_\aaa\bb\ra\not=0$. Then (\ref{teqn}) gives the following
\be
\|\bb\|^2&=&\|s_\aaa \bb\|^2
\\
&=&
\frac{2}{\la\LL,T\ra}\cdot
\la\LL,s_\aaa\bb\ra\cdot \la s_\aaa\bb,T\ra
\\&=&
\frac{2}{\la\LL,T\ra}\cdot
\left(
\la\LL,\bb\ra-\frac{2 \la \aaa,\bb\ra}{\|\aaa\|^2}\la\LL,\aaa\ra
\right)\cdot
\left(
\la\bb,T\ra -\frac{2 \la \aaa,\bb\ra}{\|\aaa\|^2}\la\aaa,T\ra
\right)
\\&=&
2\cdot
\left(
\la\LL,\bb\ra-\frac{2 \la \aaa,\bb\ra}{\|\aaa\|^2}\la\LL,\aaa\ra
\right)\cdot
\left(
\frac{\|\bb\|^2}{2\la\LL,\bb\ra} -\frac{2 \la \aaa,\bb\ra}{\la\LL,\aaa\ra}
\right)
\\
&=&
2\cdot
\left(
\frac{\|\bb\|^2}{2}-
2  \la \aaa,\bb\ra\
\frac{\la\LL,\bb\ra}{\la\LL,\aaa\ra}
-
2 \la \aaa,\bb\ra\
\frac{\la\LL,\aaa\ra}{\la\LL,\bb\ra}\frac{\|\bb\|^2}{\|\aaa\|^2}
+
4\ \frac{2 \la \aaa,\bb\ra^2}{\|\aaa\|^2}
\right).
 \ee
Subtracting $\|\bb\|^2$ and multiplying by the denominators gives
\[
0\ =\ \la\aaa,\bb\ra \left(
\|\bb\|^2\la\LL,\aaa\ra^2 + \|\aaa\|^2\la\LL,\bb\ra^2 -2\la\bb,\aaa\ra\la\LL,\aaa\ra\la\LL,\bb\ra\right).
\]
But this gives the following pair of equations
\be
0& =& \la\aaa,\bb\ra \Big(
\underbrace{\left(\|\bb\|\la\LL,\aaa\ra + \|\aaa\|\la\LL,\bb\ra\right)^2}_{>0}
-2\underbrace{\left(
\|\aaa\|\|\bb\|+\la\bb,\aaa\ra\right)}_{>0}\la\LL,\aaa\ra\la\LL,\bb\ra\Big)\\
0& =& \la\aaa,\bb\ra \Big(
\underbrace{\left(\|\bb\|\la\LL,\aaa\ra - \|\aaa\|\la\LL,\bb\ra\right)^2}_{>0}
+2\underbrace{\left(\|\aaa\|\|\bb\|-\la\bb,\aaa\ra\right)}_{>0}\la\LL,\aaa\ra\la\LL,\bb\ra\Big).
\ee
This implies $\la\aaa,\bb\ra=0$ or $\la\LL,\aaa\ra\la\LL,\bb\ra=0$, but this was excluded.
This argument is symmetric in $\aaa$ and $\bb$ hence we get the same result for $s_\bb\aaa$. Thus
we have proved that $\la\LL,s_\aaa\bb\ra=\la\LL,s_\bb \aaa\ra=0$ or $\la\aaa,\bb\ra=0$.

Now $\la\LL,s_\aaa\bb\ra=\la\LL,s_\bb \aaa\ra=0$ implies 
$\la\LL,\aaa\ra=\frac{2\la\aaa,\bb\ra}{\|\aaa\|^2}\cdot \frac{2\la\aaa,\bb\ra}{\|\bb\|^2} \cdot
\la\LL,\aaa\ra$. Since $\la\LL,\aaa\ra $ was supposed to be non zero we have that
$\frac{2\la\aaa,\bb\ra}{\|\aaa\|^2}\cdot \frac{2\la\aaa,\bb\ra}{\|\bb\|^2}=1$ which implies --- since
both factors are in $\mathbb{Z}$ --- that $\|\aaa\|^2=\|\bb\|^2$. This holds if $\la\aaa,\bb\ra\not=0$.

\item{Claim 2:} {\em All roots in $\Delta$ have the same length.}

Suppose we have short and long roots. Then we can write a
long root $\alpha$ as the sum of two short ones, lets say
$\aaa=\beta+\gamma$. This implies   $\la\aaa,\bb\ra\not=0$ and
$\la\aaa,\gamma\ra\not=0$. Since $\alpha$ is long and $\beta$ and $\gamma$ are short we have by the first
claim that $\la\LL,\aaa\ra\cdot \la\LL,\beta\ra=0$ and $\la\LL,\aaa\ra\cdot \la\LL,\gamma\ra=0$. Now
$\la\LL,\aaa\ra=\la\LL,\beta\ra+\la\LL,\gamma\ra$ gives that $\la\LL,\aaa\ra=0$ for every long root.
But this is impossible.
 Hence all roots have the same length and in particular holds for non-proportional roots
\beq\frac{2\la\aaa,\bb\ra}{\|\aaa\|^2}=\pm1.
\label{gleichlang}\eeq
\item{Claim 3:}
{\em There is an $a\in\mathbb{N}$ such that for every root $\aaa$ holds $\la\LL,\aaa\ra\in\{0,\pm a\}$.
Furthermore $a$ is less or equal than the length of the roots.}

We consider $\aaa\in \D$ with $\la\LL,\aaa\ra\not=0$ and set $a:=\la\LL,\aaa\ra$. Then we
define the vector space $A:= span \{\beta\in\D\ |\ \la\LL,\beta\ra=\pm a\}\subset \mf{t}^*$. We show that
$A=\mf{t}^*$ and that every root $\gamma$ with $\la\LL,\gamma\ra\not\in\{0,\pm a\}$ is orthogonal to $A$.

To verify $A=\mf{t}^*$ we show that every root is either in $A$ or in $A^\bot$. First consider
$\gamma\in\D$ with $\la\LL,\gamma\ra=0$. If it is not in $A^\bot$ then there is a root
$\beta\in A$ and a $\delta\not\in A$ such that $\gamma=\beta+\delta$.
But this implies $0=\la\LL,\gamma\ra=\la\LL,\beta\ra+\la\LL,\delta\ra=\pm a+\la\LL,\delta\ra$.
Hence $\delta\in A$ and therefore $\gamma\in A$ which is a contradiction. Thus $\gamma\in A^\bot$.

Now we consider a root $\gamma$ with $\la\LL,\gamma\ra\not\in \{0,\pm a\}$. For any
$\bb$ with $\la\LL,\bb\ra=\pm a$  then we have because of
(\ref{gleichlang}) that $\la\LL,s_\beta\gamma\ra=\la\LL,\gamma\ra\pm a\not=0$. Because of the proof of
claim 1 this gives $\la\beta,\gamma\ra=0$. Hence $\gamma\in A^\bot$. Since the root system is indecomposable
we have that $A=\mf{t}^*$. Furthermore we have shown that any root with $\la\LL,\gamma\ra\not\in \{0,\pm a\}$
is orthogonal to $A=\mf{t}^*$. Thus the first part of claim 3 is proved.

Now we suppose that $a > c$ where $c$ denotes the length of the roots.
We consider an $\ä\in\D$ with $\la\LL,\ä\ra=a$.
$s_\ä(\LL)=\LL-\frac{2a}{c}\ä$ is an extremal weight in $ U$. Then
$a>c$ implies 
$\LL-2\ä\in \Omega$ but not in $U$. Then the existence of the planar spanning triple $(\LL,-\LL,U)$ implies 
$\LL-2\ä=-\LL+\beta$ for a $\beta\in \D$. Hence
\[\frac{2\la\LL,\gamma\ra}{c}\ =\ 1+\frac{2\la\ä,\gamma\ra}{c}\ =\ 2\]
and therefore $\la\LL,\gamma\ra=a$ and $a=c$ which is a contradiction.
\end{description}

Now we consider for  any $\aaa\in \D$ the set $\D_\aaa^\bot:=\{\beta\in \D\ |\ \la\aaa,\bb\ra=0\}\subset \D$.
This set is a root system, reduced but not necessarily  indecomposable. But we can make the following claim.

\begin{description}
\item{Claim 4:} {\em Let $\alpha\in\D$ with $\la\LL,\ä\ra\not=0$. Then  one of the following cases holds:
\bnum
\item $\D_\aaa^\bot$ is orthogonal to $\LL$ or
\item there is a unique
$\beta\in\D_\aaa^\bot$ with $\la\LL,\bb\ra\not=0$ such that
\bnum
\item $\LL=\pm \frac{a}{c}(\aaa+\bb)$ where $c$ is the lengths of the roots, and
\item $\D_\aaa^\bot$ is decomposable with a direct summand $A_1=\{\pm\beta\}$.
\enum
\enum
}

Suppose that there is a $\beta\in\D_\aaa^\bot$ with $\la\LL,\beta\ra\not=0$. W.l.o.g. we can suppose that
$\la\LL,\beta\ra=\la\LL,\aaa\ra=\pm a$. $\la\aaa,\bb\ra$ implies then
\[s_\aaa s_\bb(\LL)\ =\ \LL\mp \frac{2a}{c}(\aaa+\beta).\]
Now we show with the help of (\ref{teqn}) that $s_\aaa  s_\bb (\LL) $ is not in $U$:
\be
\la s_\aaa  s_\bb (\LL),T\ra
&=&
\la\LL,T\ra - \frac{2\la\LL,\aaa\ra  }{\|\aaa\|^2}\la\aaa,T\ra
-\frac{2\la\LL,\bb\ra  }{\|\bb\|^2}\la\bb,T\ra
\\
&=&-\la\LL,T\ra\not=0.
\ee
But this implies
$-\LL=s_\aaa s_\bb (\LL)=\LL\pm \frac{2a}{c}(\aaa+\bb)$. By this equation $\aaa$ determines $\bb$ uniquely.

We still have to show that such  $\beta$ is orthogonal to all other roots in $\D_\alpha$.
For $\gamma\not\sim\beta$ in $\D_\alpha$ we have
\[\la\LL,s_\beta\gamma\ra= \underbrace{\la \LL,\gamma\ra}_{=0} - \frac{2\la\bb,\gamma\ra}{\|\bb\|^2}\la\LL,\bb\ra.\]
The uniqueness of $\bb$ implies that $\bb$ is orthogonal to $\D_\alpha$.

\item{Claim 5:} {\em The root system of $\lag$ is of type $A_n$ or $D_n$.}

The only root system with roots of equal length where the root system $\D_\ä^\bot$ is decomposable for a root $\ä$ is $D_n$.
Hence for every root system different
from $D_n$ we have that $\Delta_\ä^\bot \bot\  \LL$ by claim 4. Any root system different from $A_n$
satisfies that $span(\D_\ä^\bot)=\ä^\bot$. Both together imply that for any root system different
from $D_n$ and $A_n$ we have that $\ä=\LL$ but this was excluded.
\end{description}
To find the representations of  $A_n$ and $D_n$ which obey the above claims we introduce a fundamental system
$\Pi=(\pi_1,\ldots, \pi_n)$ which makes $\LL$ to the highest weight of the representation. Then we have that
$\LL=\sumk m_k\w_k$ with $m_k\in\mathbb{N}\cup\{0\}$ and $\w_k$ the fundamental representations.
$\la\w_i,\pi_j\ra=\delta_{ij}$ implies
$m_i=\la\LL,\pi_i\ra\in\{0,a\}$. Then we get
\begin{description}
\item{Claim 6.} {\em The root system is of type $D_n$ and the
 representation  is  the $a$-th power of a fundamental representation,
i.e. $\LL=a\w_i$.}

Applying $\LL$ to the root $\sumk \pi_k$ gives $\sumk m_k=a$. Applying $\LL$ to any of the $\pi_i$ gives that
$\sumk m_k=m_i$ for one $i$.

Now we consider the root system $A_n$. $n=1$ was excluded from the beginning. Recalling $A_3\simeq D_3$ we can
also exclude $A_3$. Now we impose the condition that the representation is orthogonal.
This forces $n$ to be odd and $\Lambda=a \w_{\frac{n+1}{2}}$ where $a $ has to be $2$ when
$\frac{n+1}{2}$ is odd. Thus we can suppose that $n>3$.
Using the usual notation  we consider now the root
$\sumk \pi_k=e_1-e_{n+1}$ for which holds that $\la\LL,\eta\ra=a$. Hence
by claim 4 we have that $\D_\eta^\bot $ is orthogonal to $\LL$. On the other hand
$\D_\eta^\bot=\{\pm (e_i-e_j)\ |\ 2\le i<j\le n\}$ with $n>3$ is not orthogonal
to $a \w_{\frac{n+1}{2}}=a\left(e_1+\ldots +e_{\frac{n+1}{2}}\right)$.
This yields a contradiction.
\end{description}

Finally we show that only the representations of $D_n$ given in the proposition satisfy the derived properties.
The fundamental representations of $D_n$ are given by
$\w_i=e_1+\ldots +e_i$ for $i=1\ldots n-2$ and $\w_i=\einhalb (e_1+\ldots +e_{n-1}\pm e_n)$ for $i=n-1, n$.
Then $\la a\w_i,\pi_i\ra=a$. On the other hand for the largest root $\eta=e_1+e_2$ holds
\[
\la a\w_i,\eta\ra=\left\{
\begin{array}{rcl}
a&:&i=1,n-1,n\\
2a&:&2\le i\le n-2.
\end{array}
\right.\]
Hence the representation of $a\w_i$ with $2\le i\le n-2$ does not satisfy claim 3.
Now we consider for $n>4$ the representations $\LL=\einhalb(e_1+\ldots +e_{n-1}\pm e_n)$.
For the root $\ä=e_{n-1}\pm e_n$ holds that $\la\LL,\ä\ra=a\not=0$.
The roots $\bb_1:=e_1-e_2$ and $\gamma:=e_1+e_3$ both satisfy
$\la\LL,\bb\ra=\la\LL,\gamma\ra=a$ and $\la \ä,\bb\ra=\la\ä,\gamma\ra=0$. But this is a violation of
the uniqueness property in claim 4. Hence $n=4$.

For $D_4$ it holds that, $\w_3$ and $w_4$ are congruent to $\w_1$, i.e. there is an
involutive automorphism of the Dynkin diagram which interchanges $\w_1$ with $\w_3$ respectively $\w_1$ with
$\w_4$. For $D_3\simeq A_3$
only the representations $\w_2$ and $2\w_2$ are orthogonal.
\eprf

Again we get a
\bfolg\label{planarfolge}
Every representation of a Lie algebra which satisfies the conditions of proposition
\ref{planarsatz} is the complexification of a Riemannian holonomy representation.
\efolg

\bprf
The representation with highest weight $\w_1$ of $D_n$ is the standard representation of
$\laso(2n,\ccc)$ in $\ccc^{2n}$. Hence it is the holonomy representation of a generic Riemannian
manifold.

The representation with highest weight $2\w_1$ is the complexified holonomy representation
of a symmetric space of type $AI$ for even dimensions,
 i.e. of $SU(2n)/SO(2n,\rr)$ respectively $Sl(2n,\rr)/SO(2n,\rr)$.
\eprf

\subsection{Representations with the property (SII) and weight zero}

Now we will study the property (SII) for representation for which zero is a weight.
For this we need a lemma.

\blem
Let $\lag\subset \laso(N,H)$ the irreducible representation of a simple Lie algebra with weights $\Omega$.
If $0\in \Omega$ then
\bnum
\item $ \D\subset \Omega$ or
\item the extremal weights are short roots or
\item $\D=C_n$ and the
representation is a fundamental one with highest weight $\w_{2k}$ for $k\ge 2$.
\enum
\elem

\bprf
$0\in\Omega $ implies that there is a $\lam\in \Omega$ and an $\eta\in \D$ such that $0=\lam-\eta$, i.e
$\lam=\eta$. Now we consider two cases.
\begin{description}
\item{{\em Case 1: $\eta$ is a long root.}}

Of course we have that the root system of long roots is contained in $\W$. We have to show that the short roots are in
$\W$. This is the case if one short root is in $\W$. For this we write
$\eta=\ä+\bb$ where $\ä$ and $ \bb$ are short roots.
If $\D\not= G_2$ we have that $\la\ä,\bb\ra=0$. In this case we have that
$\frac{2\la \eta,\ä\ra}{\|\ä\|^2}=2$, i.e. $\eta-\ä=\bb\in \W$.
For $\D=G_2$ we have that $\frac{2\la\ä,\bb\ra}{\|\ä\|^2}=\frac{2\la\ä,\bb\ra}{\|\bb\|^2}=1$ and therefore
$\frac{2\la \eta,\ä\ra}{\|\ä\|^2}=3$, i.e. $\eta-\ä=\bb\in \W$ too. Hence also the short roots are weights and we have
$\D\subset \W$.

\item{{\em Case 2: $\eta$ is a short root.}}

Again the short roots are weights. We have to show that one long root is a weight if $\eta$ is not extremal
or that we are in the case of the $C_n$ with the above representations.
If $\eta $ is not extremal then exists an $\alpha\in \D$ such that $\eta+\alpha\in \W$ and $¸\eta-\ä\in \W$.
This $\ä$ we fix and consider the following cases.

\begin{description}
\item{{\em Case A: $\alpha=\eta$, i.e. $2\eta\in\W$.}} If $\D\not= G_2$ we find a long root $\bb$ such that
$\frac{2\la \eta,\bb\ra}{\|\eta\|^2}=-2$.  This implies that $\bb+2\eta $
is a long root but also a weight. In case of $G_2$ we find a short root $\bb$ with
$\la\eta,\bb\ra<0$ and such
that $2\eta+\bb\in \D$ a long root. This long root is also in $\W$ since
$\la\eta,\bb\ra<0$.

\item{{\em Case B:  $\ä\not\sim\eta$ and $\la\ä,\eta\ra\not= 0$.}}
First we consider the case where $\ä$ is a short root. W.l.o.g. let be $\la\ä,\eta\ra< 0$
Then  $\ä+\eta$ is a root and a weight. If $\D$ is different from $C_n$ it is a long root and we are ready.
For $C_n$ we have to analyze the situation in detail (see the appendix of \cite{knapp96}):
Let $\eta=e_i+e_j$ and $\ä=e_k-e_j$ with $i\not=k$ be the two short roots. Since $\W
\ni \eta-\ä=2e_j+e_i-e_k$ we have that
$\frac{2\la \eta-\ä,e_i-e_k\ra}{\|e_i-e_k\|^2}=2$. Hence $\eta-\ä-(e_i-e_k)=2e_j\in\W$.
But $2e_j$ is a long root of $C_n$ and we are ready.

If $\ä$ is a long root we proceed as follows. For $G_2$ one of $\eta\pm\ä$ is a short root, lets say
$\eta-\ä$. Then we have that $\la \eta+\ä,\eta-\ä\ra<0$ hence $2\eta$ is a weight and we may argue as in the
first case A. If $\D$ is different from $G_2$ we  write $\ä=\ä_1+\ä_2$ with two
orthogonal short roots $\ä_1$ and $\ä_2$. For one of these is $\la\eta,\ä_i\ra\not=0$ and hence
$\eta\pm \ä_i$ a long root, but also a weight.

\item{{\em Case C: $\la\ä,\eta\ra=0$ and $\Delta\not=C_n$.}} For $G_2$ this case implies that $\ä$ is a long root and that
$\eta+\ä$ is two times a short root. Hence for $G_2$ we can proceed as above to get the result.

If $\D$ is different from $G_2$ we consider the root system $\D_\eta^\bot$ of roots orthogonal to $\eta$,
which contains $\ä$. In case of $C_n$ this root system is equal to
$A_1\times C_{n-2}$ and in the remaining cases --- $B_n$ and $F_4$ --- equal to $B_{n-1}$ resp. $B_3$.
Now we show that there is a short root $\ä_1$ in $\D_\eta^\bot$ such that
$\eta+\ä_1\in\W$. If $\ä$ is short this is trivial and if $\ä$ is long we write
$\ä=\ä_1+\ä_2$ with two orthogonal short roots from $\D_\eta^\bot$. Then $\la\eta+\ä,\ä_2\ra>0$ and thus
$\eta+\ä_1\in\Omega$.

On the other hand  there is a short root $\gamma\in \D_\eta^\bot$ with $\eta+\gamma $ is a long root.
Applying now the Weyl group of $\D_\eta^\bot$ on $\eta+\gamma$ we get that $\eta+\ä_1 $ is a long root.
In case of $C_n$ this argument does not apply since $\gamma$ spans the $A_1$ factor of $\D_\eta^\bot$.
\end{description}
Hence we have verified $\D\subset \W$ in the cases A, B and C.
It remains to show that in the situation where
$\la\ä,\eta\ra=0$, $\Delta=C_n$ and neither case A nor case B applies, it holds that
$\D\subset \W$ or
the representation of $C_n$ is the one with highest weight
$\w_{2k}$ with $k\ge2$.

We suppose that $\D\not\subset \W$. Hence no long root can be a weight.

First of all we show that under these conditions $\ä$ has to be a short root.
This is true because
$\frac{2\la \eta\pm\ä,\eta\ra}{\|\eta\|^2}=2$ implies $\W\ni \eta\pm\ä-\eta=\pm\ä$. Hence $\ä$ has to be short.

Secondly we note that neither $\eta+\ä$ nor $\eta- \ä$ can be a root because it would be a long root and a weight.
This implies  $n\ge 4$.

In a third step we show that there is no long root $\bb$ such that $\eta+\ä+\bb\in \W$ and
$\eta+\ä-\bb\in \W$. We consider the number
\beq\label{pmbeta}
\frac{2\la \eta+\ä\pm\bb,\ä\ra}{\|\ä\|^2}=2 \pm \frac{2\la \ä,\beta\ra}{\|\ä\|^2}.
\eeq
If $\la\ä,\bb\ra=0$ we have that $\eta+\ä\pm\bb-\ä=\eta\pm\beta\in \W$. But this was excluded (First step or
case B).
Hence we suppose that $\frac{2\la \ä,\beta\ra}{\|\ä\|^2}=2$. We still have that $\eta+\bb\in \Omega$.
We consider the number
$\frac{2\la \eta+\bb,\eta\ra}{\|\eta\|^2}=2 \pm \frac{2\la \eta,\beta\ra}{\|\eta\|^2}\ge 0.$ If this is not zero
we have that $\W\ni\eta+\beta-\eta=\bb$ which was excluded. Hence
$\frac{2\la \eta,\beta\ra}{\|\eta\|^2}=-2$.
But this together with $\frac{2\la \ä,\beta\ra}{\|\ä\|^2}=2$ is a contradiction since the long roots of $C_n$
are of the form $\pm 2e_i$ and the short ones of the form $\pm(e_i\pm e_j)$.
Hence if there is a root such that $\eta+\ä+\bb\in \W$ and $\eta+\ä-\bb\in \W$, it has to be a short one.

If there is no such $\bb$ then $\eta+\ä$ is extremal. Considering the fundamental weights of $C_n$
this gives easily that the highest weight of the representation is $\w_4$.

Finally we suppose that there is such a short root $\bb$. Since $\bb$ is short equation (\ref{pmbeta})
implies $\eta\pm\bb\in \W$. Since we have excluded case A and B it must hold $\la\eta,\bb\ra=0$ and neither
$\eta+\bb$ nor $\eta-\bb$ is a root. On the other hand the same holds for $\ä$ and $\beta$ since any other
would imply that $\ä\pm\bb$ is a long root which was excluded or a short root $\gamma$ orthogonal to $\eta$
and with $\eta\pm\gamma\in\W$. This way we go on attaining that any extremal wight is the sum
of orthogonal short roots whose pairwise sum is no long root. But this is nothing else than the fact that the
highest weight of the representation is $\w_{2k}$ for $k\ge 2$.

\end{description}All in all we have shown the proposition.
\eprf

\bs\label{s2satz0}
Let $\lag \subset \laso(N, \ccc ) $
be an irreducible representation of real type of a  complex simple Lie algebra different from
$\mf{sl}(2,\ccc)$ and satisfying (SII).
If  $\ 0\in\W$ then there is a root $\ä$
such that for the
extremal weight from property (SII) holds
$\Lambda = a \alpha$ or the representation is  congruent to one of the following:
\bnum
\item $\D=C_4$ with highest weight $\w_4$.
\item $\D=D_n$ with highest weight $2\w_1$.
\enum
\es

\bprf
Let $\LL$ and $\ä$ be the extremal weight and the root from property (SII). We suppose that $\LL$ is
not the multiple of a root.

First of all we consider the case where $0\in\W_\ä$. By the previous lemma this is true in the
following cases:
\bnum
\item[(a)] $\D\not=C_n$, because in this case $\Delta\subset \W$.
\item[(b)] $\D=C_n$ but the highest weight of the representation is not equal to $\w_{2k}$ with $k\ge 2$, because
this again implies $\D\subset \W$.
\item[(c)] $\D=C_n$ and $\alpha$ is a short root, because for representations with $0\in\W$ holds that
the short roots are weights.
\enum
For $0\in\W_\ä$  property (SII) gives
$0=\LL -\ä -\bb$ or $0=-\LL+\bb$. The second case was excluded thus we
have to consider the first case. Suppose that
$\LL=\ä+\bb$ where $\ä+\bb\not\sim \gamma\in\Delta$. In particular $\ä+\bb$ is not a root which implies
that $\la\ä,\bb\ra\ge 0$. We consider three cases.

\begin{description}
\item{{\em Case 1: $\D=G_2$.}} In this case $\ä+\beta\not\sim\gamma\in\D$ implies 
$\la\ä,\beta\ra>0$ and $ \ä$ and $\bb$ must have different length. Thus we can chose a long root
$\gamma$ not proportional neither to $\ä$ nor to $\bb$ and such that $\la\ä,\gamma\ra<0$
and $\la\bb,\gamma\ra<0$ which implies 
$\gamma\in\W_\ä$ as well as $\gamma\in \W_\bb$. (SII) implies then
$\gamma-\bb\in\D$ or $\gamma-\ä\in\D$ or $\gamma+\ä+\bb\in\D$. The first two cases are not possible
because of lemma \ref{wurzellemma0}.
For the third case we suppose that $\ä$ is the long root and consider
$\frac{2\la\gamma+\bb,\ä\ra}{\|\ä\|^2}=0$ because $\ä$ is long and both terms have opposite sign.
Hence $\gamma+\ä+\bb$ can not be a root.

\item{{\em Case 2: $\D\not=G_2$ and $\la\ä,\bb\ra>0$.}} This implies $\ä-\bb\in\Delta$.
We consider the number
$k:= \frac{2\la\ä,\ä+\bb\ra}{\|\ä\|^2}=2+ \frac{2\la\ä,\bb\ra}{\|\ä\|^2}\ge 3$. Since $G_2$ was excluded
we have that $k\in\{3,4\}$.
Hence $\ä+\bb-k\ä=\bb-(k-1)\ä\in\W_\ä$. Then property (SII) implies 
$\bb-(k-1)\ä=-\ä-\bb+\gamma$ with $\gamma \in \D_0$, i.e.
$2\bb-(k-2)\ä\in\Delta$. At first this implies $k=3$ and thus
$\frac{2\la\ä,\bb\ra}{\|\ä\|^2}=1$. Secondly we must have $\frac{2\la\ä,\bb\ra}{\|\bb\|^2}=2$, therefore
$\|\ä\|^2=2\|\bb\|^2$, i.e. $\ä$ as well as $2\bb-\ä$ are long roots and $\bb$ and $\bb-\ä$ are short ones.

This implies
$\frac{2\la\bb-\ä,\ä+\bb\ra}{\|\bb-\ä\|^2}= \frac{2(\|\bb\|^2-\|\ä\|^2)}{\|\bb\|^2}=-2$.
Hence $\ä+\bb+2(\bb-\ä)=3\bb-\ä\in\W$ and since
$\frac{2\la\ä,\ä-3\bb\ra}{\|\ä\|^2}=2-3=-1$ holds
$\ä-3\bb\in\W_\ä$. (SII) then gives
$\ä-3\bb=\bb-\gamma$ or $\ä-3\bb=-\bb-\ä+\gamma$ with $\gamma\in\D_0$. But none of these
equations can be true.
\item{{\em Case 3: $\la\ä,\bb\ra=0$ and $\D\not=G_2$.}}
Since $\ä+\bb\not\sim\gamma \in\D$ the rank of $\D$ has to be greater than $3$ or it is $\D=D_n$ and
$\LL=2e_i$, i.e. $\LL=2w_1$. In the second case we are ready and we exclude this representation in the following.
We can suppose $rk \D\ge 4$.
In this situation we prove the following lemma.
\blem
Let $rk\D\ge 4$ and let $\LL=\ä+\bb$ be an extremal weight of a representation satisfying property
(SII) for $(\LL,-\LL+\ä,\ä)$ with $\bb\in\D$ satisfying $\la\ä,\bb\ra=0$ and
$\ä+\bb\not\sim\gamma\in\D$. Then $\D$ is a root system
with roots of the same length or $\D=C_n$ and $\ä$ and $\bb$ are two short roots.
\label{lemma1}
\elem
\bprf
Suppose that $\D$ has roots of different length.

First we assume that $\bb$ is a long root. We consider the root system $\D_\ä^\bot$ which contains $\bb$.
We notice that $\bb$ lies not in an $A_1$ factor of $\D_\ä^\bot$ because otherwise $\ä+\bb$ would be the multiple
of a root.
Since $\bb$ is long we find a short root $\gamma\in\D_\ä^\bot$ such that
$\frac{2\la\bb,\gamma\ra}{\|\gamma\|^2}=-2$. Hence $\ä+\bb+2\gamma\in\W$ and
--- since
$\frac{2\la\ä,\ä+\bb+2\gamma\ra}{\|\ä\|^2}=2$ --- it is $-\ä-\bb-2\gamma\in\W_\ä$.
But this contradicts property (SII).

Now we suppose that $\ä$ is a long root. Here we consider the root system $\D_\bb^\bot$ containing $\ä$.
Again $\ä$ lies not in an $A_1$ factor of $\D_\bb^\bot$ because otherwise $\ä+\bb$ would be the multiple of a root.
Since $\ä$ is long we find a short root $\gamma\in\D_\bb^\bot$ such that
$\frac{2\la\ä,\gamma\ra}{\|\gamma\|^2}=-2$. Hence $\ä+\bb+2\gamma\in\W$.
Now we have that $\frac{2\la\ä,\gamma\ra}{\|\ä\|^2}=-1$ and therefore
$\frac{2\la\ä,\ä+\bb+2\gamma\ra}{\|\ä\|^2}=2-1=1$. Thus $-\ä-\bb-2\gamma\in\W_\ä$.
Again this contradicts (SII).

If $\ä$ and $\bb$ are short and orthogonal and the root system is not $C_n$, i.e. it is $B_n$ or $F_4$,
then the sum of two orthogonal short roots is the multiple of a root.
\eprf

Now we prove a second
\blem\label{lemma2}
The assumptions of the previous lemma imply that there is no $\gamma\in\D$
such that
\beq
\la\ä,\gamma\ra=0\mbox{ and }\frac{2\la\bb,\gamma\ra}{\|\gamma\|^2}=1.
\label{gamma}\eeq
\elem
\bprf
Lets suppose that there is a $\gamma\in \D$ such that
$\la\ä,\gamma\ra=0$ and $\frac{2\la\bb,\gamma\ra}{\|\gamma\|^2}=1$. In case of $C_n$ $\gamma $ is a short root.
We note that both together imply that neither $\ä+\gamma$ nor $\ä-\gamma$ is a root.
But $\gamma-\bb$ is a root, in case of $C_n$ a short one. Furthermore $\LL-\gamma\in \W$
Hence
\[
\frac{2\la\LL-\gamma,\gamma-\bb\ra}{\|\gamma-\bb\|^2}=
\frac{2\la\ä+\bb-\gamma,\gamma-\bb\ra}{\|\gamma-\bb\|^2}
=-2
.\]
Hence $\LL-\gamma+2(\gamma-\bb)=\ä-\bb+\gamma\in \W$. Now
$\frac{2\la\ä-\bb+\gamma,\ä\ra}{\|\ä\|^2}=2$, i.e.
$-\ä+\bb-\gamma\in \W_\ä$.
(SII) implies now that $-\ä+\bb-\gamma=\bb+\delta$ or $-\ä+\bb-\gamma=-\ä-\bb+\delta$ for $\delta\in \D_0$.
But both options are not possible since $\ä+\gamma$ is not a root and because $\gamma$ is short.
\eprf
We conclude that lemma \ref{lemma1} left us with representations of
$A_n$, $D_n$, $E_6,\ E_7,\ E_8$ or $C_n$ where $\LL$ is the sum of two orthogonal (short) roots but not a root.

Now one easily verifies that lemma \ref{lemma2} implies  $n\le 4$ and $\D\not=A_4$.
Hence the remaining representations are $2\w_1$, $2\w_3$ and $2\w_4$ of $D_4$, which are congruent
to each other, and $w_4$ of $C_4$.
\end{description}
To finish the proof we have to consider the representation of highest weight $\w_{2k}$ (with $k\ge 2$)
of $C_n$ supposing
$\ä$ is a long root. $0\in\W$ implies that the short roots are weights.
Let $\bb$ be a short root with $\la\ä,\bb\ra<0$, i.e. $\beta\in\W_\alpha$. (SII) then gives
$\bb=\w_{2k}-\ä+\delta$ or $\bb=\w_{2k}-\delta$ for  a $\delta \in\Delta_0$.
Analyzing roots and fundamental weights of $C_n$ we get that (SII) implies $k=2$ and
$\alpha=2e_i$ for $1\le i\le 4$. But for $n> 4$  lemma \ref{lemma2} applies analogously.
The remaining representation is $\w_4$ of $C_4$.
\eprf

\bfolg
Let $\lag\subset\laso(N,\ccc)$ be an orthogonal algebra of real type different from $\mf{sl}(2,\ccc)$
and satisfying (SII).
If $0\in\W$, in particular if
$\D =G_2,F_4$ or $E_8$ then it is the complexification of a Riemannian holonomy representation
with the exception of $G_2$ in corollary \ref{wurzelgewichtfolge}.
\efolg

\bprf
If $\LL$ is the multiple of a root then we are in the situation of corollary \ref{wurzelgewichtfolge}.
For
$D_n$ the remaining representations are those which appear in corollary \ref{planarfolge}.
The representation of highest weight $\w_4$ of $C_4$ is the complexification of the holonomy representation
of the Riemannian symmetric space of type $EI$, i.e. of $E_6/Sp(4)$ resp. $E_{6(6)}/Sp(4)$.

Furthermore analyzing the roots and fundamental representations of the exceptional algebras we notice that
every representation of $G_2$, $F_4$ and $E_8$ contains zero as weight.
\eprf

\subsection{Representations with the property (SII) where zero is no weight}

First we need a
\blem
Let $0\not\in \Omega$. Then there is a weight $\lam\in \W$, such that for every root holds
$\left|\frac{2\la\lam,\ä\ra}{\|\ä\|^2}\right|\le 1$.
\elem
\bprf
Let $\lam$ be a weight and $\ä$ a root such that
$\frac{2\la \lam,\ä\ra}{\|\ä\|^2}=:k\ge 2$. If
$k$ is even we have that $0\not=\lam-\frac{k}{2}\ä\in\W$. But for this weight holds
 $\frac{2\la\lam-\frac{k}{2}\ä,\ä\ra}{\|\ä\|^2}=k-k=0$.
If $k$ is odd we have that $0\not=\lam-\frac{k-1}{2}\ä\in\W$ and
 $\frac{2\la\lam-\frac{k-1}{2}\ä,\ä\ra}{\|\ä\|^2}=1$.
 \eprf

\bs\label{s2satz1}
Let $\lag \subset \laso(N, \ccc ) $
be an irreducible representation of real type of a  complex simple Lie algebra different from
$\mf{sl}(2,\ccc)$, with $0\not\in\W$ and satisfying (SII). Then
$\left|\frac{2\la \LL,\bb\ra}{\|\bb\|^2}\right|\le 3$ for all roots
$\bb\in\D$.
\es

\bprf
Let $\ä\in \D$ with the property (SII), i.e.
$\Omega_\alpha\subset \left\{ \Lambda-\alpha+\beta\ |\ \beta\in \Delta_0\right\}\cup
\left\{ -\Lambda+\beta\ |\ \beta\in \Delta_0\right\}$.

By the previous lemma there is a $\lam\in\Omega$ such that
$\left|\frac{2\la \lam,\bb\ra}{\|\bb\|^2}\right|\le 1$ for all roots $\bb\in\D$. Applying the Weyl group one
can choose $\lam$ such that $\la\lam,\ä\ra<0$.

$\la\lam,\ä\ra<0$ implies  $\lam\in\Omega_\ä$.
Hence (SII) gives 
$\lam=\LL-\ä-\gamma$ or $\lam=-\LL+\gamma$ with $\gamma\in \D_0$. The second case gives
for every $\bb\in\D$
\[
\left|\frac{2\la \LL,\bb\ra}{\|\bb\|^2}\right|\le
\left|\frac{2\la \lam,\bb\ra}{\|\bb\|^2}\right|+
\left|\frac{2\la \gamma,\bb\ra}{\|\bb\|^2}\right|\le 3,
\]
because we have excluded $G_2$.

Thus we have to consider the first case
$\LL=\lam+\ä+\gamma$ with $\gamma\in\D_0$ and have to verify that
\beq\label{kleiner3}
\left|\frac{2\la \LL,\bb\ra}{\|\bb\|^2}\right|=
\left|\frac{2\la \lam,\bb\ra}{\|\bb\|^2}+
\frac{2\la \ä,\bb\ra}{\|\bb\|^2}+
\frac{2\la \gamma,\bb\ra}{\|\bb\|^2}\right|\le 3.
\eeq
 for all roots
$\bb\in\D$.

For $\beta=\pm\ä$ this is satisfied:
\[
\frac{2\la \LL,\ä\ra}{\|\ä\|^2}=
\pm\frac{2\la \lam,\ä\ra}{\|\ä\|^2}\pm 2+
\frac{2\la \gamma,\ää\ra}{\|\ä\|^2}
=\mp 1\pm 2 +\frac{2\la \gamma,\ää\ra}{\|\ä\|^2}\le 3.
\]
Now we have to show (\ref{kleiner3}) for all $\bb\in\D$ with $\bb\not\sim \ä$. For this we
consider three cases.
\begin{description}
\item{{\em Case 1: All roots have the same length.}}
This implies $\left|\frac{2\la \gamma,\bb\ra}{\|\bb\|^2}\right|\le 1$ for all roots which are not
proportional to each other.
Thus we have (\ref{kleiner3}) for all $\bb\not\sim\gamma $:
\[
\left|\frac{2\la \LL,\bb\ra}{\|\bb\|^2}\right|\le
\left|\frac{2\la \lam,\bb\ra}{\|\bb\|^2}\right|+
\left|\frac{2\la \ä,\bb\ra}{\|\bb\|^2}\right|+
\left|\frac{2\la \gamma,\bb\ra}{\|\bb\|^2}\right|\le 3.
\]
For $\beta=\pm\gamma$ we have
\[
\frac{2\la \LL,\gamma\ra}{\|\gamma\|^2}=
\frac{2\la \lam,\gamma\ra}{\|\gamma\|^2}
+\frac{2\la \ä,\gamma\ra}{\|\gamma\|^2}
+2.
\]
This has absolute value $\ge 4  $ only if $\la\lam,\gamma\ra >0$ and $\la \ä,\gamma\ra>0$. This implies
that $\ä-\gamma$ is a root. But for this one holds
$\frac{2\la \lam,\gamma-\ä\ra}{\|\gamma-\ä\|^2}= \frac{2\la \lam,\gamma\ra}{\|\gamma-\ä\|^2}
-\frac{2\la \lam,\ä\ra}{\|\gamma-\ä\|^2}=2$
since all roots have the same length. This is a contradiction to the choice of $\lam$.
\item{{\em Case 2: There are long and short roots and $\bb$ is a long root.}}
This implies again $\left|\frac{2\la \gamma,\bb\ra}{\|\bb\|^2}\right|\le 1$ for all $\bb$ which are not
proportional to $\gamma$. This implies (\ref{kleiner3}) in this case.

For $\bb=\pm \gamma$ we argue as above, recalling that $\gamma-\alpha$ and $\ä$ have to be short roots in this
case.
Hence
$\frac{2\la \lam,\gamma-\ä\ra}{\|\gamma-\ä\|^2}= \frac{2\la \lam,\gamma\ra}{\|\gamma-\ä\|^2}
-\frac{2\la \lam,\ä\ra}{\|\gamma-\ä\|^2}\ge \frac{2\la \LL,\gamma\ra}{\|\gamma\|^2}
-\frac{2\la \lam,\ä\ra}{\|\ä\|^2} \ge 2$ which is a contradiction.
\item{{\em Case 3: There are long and short roots and $\bb$ is a short root.}}

First we consider the case where $\bb=\pm\gamma$. Again  (\ref{kleiner3}) is not satisfied
only if $\la\lam,\gamma\ra $ and $\la \ä,\gamma\ra$ are non zero and have the same sign, lets say $+$.

If $\ä$ is a short root too, then because of $\la \ä,\gamma\ra\not=0$
 lemma \ref{wurzellemma0} gives that $\ä-\gamma$ is also a short root.
Hence
$\frac{2\la \lam,\gamma-\ä\ra}{\|\gamma-\ä\|^2}= \frac{2\la \lam,\gamma\ra}{\|\gamma-\ä\|^2}
-\frac{2\la \lam,\ä\ra}{\|\gamma-\ä\|^2}= \frac{2\la \lam,\gamma\ra}{\|\gamma\|^2}
-\frac{2\la \lam,\ä\ra}{\|\ä\|^2} = 2$ yields a contradiction.

If $\ä$ is a long root, then $\gamma-\ä$ has to be a short one and we get again a contradiction:
$\frac{2\la \lam,\gamma-\ä\ra}{\|\gamma-\ä\|^2}= \frac{2\la \lam,\gamma\ra}{\|\gamma-\ä\|^2}
-\frac{2\la \lam,\ä\ra}{\|\gamma-\ä\|^2}\ge \frac{2\la \lam,\gamma\ra}{\|\gamma\|^2}
-\frac{2\la \lam,\ä\ra}{\|\ä\|^2} \ge 2$.

Now suppose that $\beta\not\sim \gamma$.
Then $\frac{2\la \LL,\bb\ra}{\|\bb\|^2}=
\frac{2\la \lam,\bb\ra}{\|\bb\|^2}+
\frac{2\la \ä,\bb\ra}{\|\bb\|^2}+
\frac{2\la \gamma,\bb\ra}{\|\bb\|^2}$ has absolute value $\ge 4$ only if all three right hand side terms have
the same sign --- lets say they are positive ---
and at least one of the last two terms has absolute value greater than one, i.e. $\gamma$ or
$\alpha$ is a long root. If $\alpha $ is a long root then $\alpha-\beta $ is a short one and
arguing as above gives the contradiction. If $\ä$ is a short root then $\la\ä,\bb\ra>0$ implies by
lemma \ref{wurzellemma0} that $\bb-\alpha$ is a short root. Again we have a contradiction:
$\frac{2\la \lam,\bb-\ä\ra}{\|\bb-\ä\|^2}= \frac{2\la \lam,\bb\ra}{\|\bb-\ä\|^2}
-\frac{2\la \lam,\ä\ra}{\|\bb-\ä\|^2}= \frac{2\la \lam,\bb\ra}{\|\bb\|^2}
-\frac{2\la \lam,\ä\ra}{\|\ä\|^2} = 2$.
\end{description}
\eprf

\bs \label{s2satz2}
Under the same assumptions as in the previous proposition holds that
$\left|\frac{2\la\LL,\eta\ra}{\|\eta\|^2}\right|\le 2$ for all long roots $\eta$.
\es

\bprf
Let $\LL$ and $\ä$ be the extremal weight and the root from property (SII).
We suppose that there is a long root $\eta$ with
\beq\label{indirekt}\frac{2\la\LL,\eta\ra}{\|\eta\|^2}=-3\eeq and
derive a contradiction considering different cases.
\begin{description}
\item{{\em Case 1: All roots have the same length.}}
By applying the Weyl group we find an extremal weight $\LL^\prime$ such that
$a:=\frac{2\la\LL^\prime,\ä\ra}{\|\ä\|^2}=-3$.

First we find a root $\bb$ with
\[ \frac{2\la\ä,\beta\ra}{\|\bb\|^2}=1\mbox{ and }
\frac{2\la\LL^\prime,\beta\ra}{\|\bb\|^2}\le-2.\]
This is obvious: We find a $\bb$ such that $\frac{2\la\ä,\beta\ra}{\|\bb\|^2}=1$. If
$\frac{2\la\LL^\prime,\beta\ra}{\|\bb\|^2}\ge-1$ then we consider the root $\ä-\bb$. It satisfies
$ \frac{2\la\ä,\ä-\beta\ra}{\|\ä-\bb\|^2}=1$ and we have
\[\frac{2\la\LL^\prime,\ä-\beta\ra}{\|\ä-\bb\|^2}=-3- \frac{2\la\LL^\prime,\beta\ra}{\|\ä-\bb\|^2}\le -2.\]

Hence we have
$\LL^\prime+k\beta\in \W$ for $0\le k\le 2$ and $\LL^\prime+k\ä\in \W$ for $0\le k\le 3$. Furthermore
\[\frac{2\la\LL^\prime+l\beta,\ä\ra}{\|\ä\|^2}=-3- \frac{2\la\LL^\prime,\ä\ra}{\|\ä\|^2}=-3+l.\]
But this gives
\[\LL^\prime + k \ä +l\bb\in \W_\ä \mbox{ for } 0\le k\le 2, 0\le k+l\le 2.\]
Among others (SII) implies the existence of $\gamma_i$ and $\delta_i$ from $\D_0$ for $i=0,1, 2$ such that
that the following alternatives must hold
\begin{align}
\LL^\prime+\ä&&=&&\LL+\gamma_0&&\mbox{ or }&&\LL^\prime&&=&&-\LL+\delta_0 \label{nulltens}\\
\LL^\prime+3\ä&&=&&\LL+\gamma_1&&\mbox{ or } &&\LL^\prime+2\ä&&=&&-\LL+\delta_1 \label{zweitens}\\
\LL^\prime+\ä+2\bb&&=&&\LL+\gamma_2&&\mbox{ or }&&\LL^\prime+2\bb&&=&&-\LL+\delta_2. \label{fuenftens}
\end{align}
First we suppose that the first alternative of (\ref{nulltens}) holds, i.e
$\LL^\prime+\ä=\LL+\gamma_0$. Since $a=-3$ and both $\LL$ and
$\LL^\prime$ are extremal we have that $\ä\not=-\gamma_0$. Hence the first case of (\ref{zweitens}) can not be
true and we have $\LL^\prime + 2\ä=-\LL+\delta_1$.
We consider now (\ref{fuenftens}): The left side of (\ref{nulltens}) gives
$\LL^\prime+2\bb+\ä=\LL+\gamma_0+2\bb$. If the left side of (\ref{fuenftens}) holds, we would have
$\gamma_0=-\beta$. Hence $\LL+\bb\in \W$ and on the other hand
$\W\ni \LL^\prime+\ä=\LL-\beta$ which contradicts the extremality of $\LL$.
Thus the right hand side of (\ref{fuenftens}) must be satisfied.
From $\LL^\prime + 2\ä=-\LL+\delta_1$ follows $\LL^\prime +2\bb=-\LL+\delta_1 +2 (\bb-\ä)$ and therefore
$\delta_1=-(\bb-\ä)$. Again we have $-\LL+(\bb-\ä)\in \W$ and  $-\LL-(\bb-\ä)\in \W$ which contradicts the
extremality of $\LL$.

If one starts with the right hand side of (\ref{nulltens}) we can proceed analogously and get a
contradiction in the case
where all roots have the same length.
\item{{\em Case 2. The roots have different length and $\alpha$ is a short root.}}
On one hand we find a short root $\bb$ which is orthogonal to $\ä$ and $\ä+\bb$ is a long root, and
on the other we can find an extremal weight $\LL^\prime$ such that
\[\frac{2\la\LL^\prime,\ä+\bb\ra}{\|\ä+\bb\|^2}=-3.\]
Since $\ä\bot\bb$ we have
\[
-3\ =\  \frac{2(\la\LL^\prime,\ä\ra+\la\LL^\prime,\bb\ra)}{\|\ä\|^2+\|\bb\|^2}\\
\ =\  \frac{1}{2} \left(
\frac{2\la\LL^\prime,\ä\ra}{\|\ä\|^2}+
\frac{2\la\LL^\prime,\bb\ra}{\|\bb\|^2} \right).
\]
Because of the previous proposition we get
\[\frac{2(\la\LL^\prime,\ä\ra}{\|\ä\|^2}\ =\
\frac{2(\la\LL^\prime,\bb\ra}{\|\bb\|^2} \ = \ -3.\]
Hence $\LL^\prime +k\ä+l\bb\in\W$ for $0\le k,l\le 3$ and therefore
$\LL^\prime +k\ä+l\bb\in\W_\ä$ for $0\le k\le 2 $ and $0\le l\le 3 $.
(SII) implies the following alternatives
\begin{align}
\LL^\prime+\ä&&=&&\LL+\gamma_0&&\mbox{ or }&&\LL^\prime&&=&&-\LL+\delta_0 \label{1nulltens}\\
\LL^\prime+\ä+3\bb&&=&&\LL+\gamma_1&&\mbox{ or } &&\LL^\prime+3\bb&&=&&-\LL+\delta_1 \label{1erstens}\\
\LL^\prime+2\ä+3\bb&&=&&\LL+\gamma_2&&\mbox{ or } &&\LL^\prime+\ä+3\bb&&=&&-\LL+\delta_2 \label{1zweitens}\\
\LL^\prime+3\ä+2\bb&&=&&\LL+\gamma_3&&\mbox{ or } &&\LL^\prime+2(\ä+\bb)&&=&&-\LL+\delta_3 \label{1drittens}\\
\LL^\prime+3\ä+3\bb&&=&&\LL+\gamma_4&&\mbox{ or } &&\LL^\prime+2\ä+3\bb&&=&&-\LL+\delta_4. \label{1viertens}
\end{align}
If the left hand side of the first alternative is valid then the left hand sides of the remaining four
can not be satisfied:
For (\ref{1erstens}) we would have $3\bb=\gamma_1-\gamma_0$ which is not possible.
(\ref{1zweitens}) would imply $3\bb+\ä=\gamma_2-\gamma_0$ which is by lemma \ref{wurzellemma3} a contradiction since
$\ä\not=-\bb$ and $\gamma_0\not=-\ä$.
(\ref{1drittens}) would imply $2(\ä+\bb)=\gamma_3-\gamma_0$. Since $\ä+\bb$ is a long root this would give
$\gamma_0=-\ä+\bb$ and $\gamma_3=\ä+\bb$ which is a contradiction to the extremality of $\LL$.
(\ref{1viertens}) would give $2\ä+3\bb=\gamma_4-\gamma_0$ which
also is not possible.

Thus for the last four equations the right hand side must hold.
Taking everything together we would get
$\ä=\delta_2-\delta_1=\delta_4-\delta_2$ and $\bb= \delta_4-\delta_3$. This gives $2 \ä= \delta_4-\delta_1$ and
thus
\[\frac{2\la\delta_4,\ä\ra}{\|\ä\|^2}- \frac{2\la\delta_1,\ä\ra}{\|\ä\|^2}= \frac{4 \|\ä\|^2}{\|\ä\|^2}=4.\]
The extremality of $\LL$ prevents that $\ä=\delta_4=-\delta_1$. Hence
$\delta_1$ and $\delta_4$ are long roots, in particular
\[
\frac{2\la\delta_4,\ä\ra}{\|\ä\|^2}=- \frac{2\la\delta_1,\ä\ra}{\|\ä\|^2}=2.\]
For $\bb$ again $\bb= \delta_4=-\delta_3$ can not hold by the extremality of $\LL$ and we have
\[0=
\frac{2\la\bb,\ä\ra}{\|\ä\|^2}= \frac{2\la\delta_4,\ä\ra}{\|\ä\|^2}- \frac{2\la\delta_3,\ä\ra}{\|\ä\|^2}
=2-\frac{2\la\delta_3,\ä\ra}{\|\ä\|^2}\]
which forces $\delta_3$ to be a long root too. Now we have a contradiction because the short root
$\bb$ is the sum of two long roots. This is impossible.

If we start with the right hand side of the first alternative one proceeds analogously.

\item{{\em Case 3. The roots have different length and $\alpha$ is a long root.}}
In this case we find an extremal weight $\LL^\prime$ such that
$
\frac{2\la\LL^\prime,\ä\ra}{\|\ä\|^2}=-3$. Now we can write $\ä=\ä_1+\ä_2$ with $\ä_1\bot \ä_2$ two short roots.
As above we get
\beq
\frac{2(\la\LL^\prime,\ä_1\ra}{\|\ä_1\|^2}\ =\
\frac{2(\la\LL^\prime,\ä_2\ra}{\|\ä_2\|^2} \ = \ -3.\label{extremal}\eeq
Again this implies $\LL^\prime +k\ä+l\bb\in\W$ for $0\le k,l\le 3$ and therefore
$\LL^\prime +k\ä+l\bb\in\W_\ä$ for $0\le k,l\le 2$.
Now (SII) implies the existence of $\gamma_i$ and $\delta_i$ from $\D_0$ for $i=0,\ldots, 8$ such that
that the following alternatives must hold
\begin{align}
&&(L)&&&&&&&&(R)&&\nonumber\\
\LL^\prime+\ä_1+\ä_2&&=&&\LL+\gamma_0&&\mbox{ or }&&\LL^\prime&&=&&-\LL+\delta_0 \label{2nulltens}\\
\LL^\prime+2\ä_1+\ä_2&&=&&\LL+\gamma_1&&\mbox{ or }&&\LL^\prime+\ä_1&&=&&-\LL+\delta_1 \label{2erstens}\\
\LL^\prime+3\ä_1+\ä_2&&=&&\LL+\gamma_2&&\mbox{ or }&&\LL^\prime+2\ä_1&&=&&-\LL+\delta_2 \label{2zweitens}\\
\LL^\prime+ \ä_1+2\ä_2&&=&&\LL+\gamma_3&&\mbox{ or }&&\LL^\prime+\ä_2&&=&&-\LL+\delta_3 \label{2drittens}\\
\LL^\prime+\ä_1+3\ä_2&&=&&\LL+\gamma_4&&\mbox{ or }&&\LL^\prime+2\ä_2&&=&&-\LL+\delta_4 \label{2viertens}\\
\LL^\prime+2\ä_1+2\ä_2&&=&&\LL+\gamma_5&&\mbox{ or }&&\LL^\prime+\ä_1+\ä_2&&=&&-\LL+\delta_5 \label{2fuenftens}\\
\LL^\prime+2\ä_1+3\ä_2&&=&&\LL+\gamma_6&&\mbox{ or }&&\LL^\prime+\ä_1+2\ä_2&&=&&-\LL+\delta_6 \label{2sechstens}\\
\LL^\prime+3\ä_1+2\ä_2&&=&&\LL+\gamma_7&&\mbox{ or }&&\LL^\prime+2\ä_1+\ä_2&&=&&-\LL+\delta_7 \label{2siebentens}\\
\LL^\prime+3\ä_1+3\ä_2&&=&&\LL+\gamma_8&&\mbox{ or }&&\LL^\prime+2\ä_1+2\ä_2&&=&&-\LL+\delta_8. \label{2achtens}
\end{align}
In the following we denote the left hand side formulas with an .L and the right hand side formulas with an .R.
Again we suppose that (\ref{2nulltens}.L) is satisfied, i.e.
$\LL^\prime+\ä_1+\ä_2=\LL+\gamma_0$. Then (\ref{extremal}) and the extremality of $\LL$
implies that $\gamma_0$ does not equal to $\ä_i$.

Now
(\ref{2achtens}.L) would imply that $2(\ä_1+\ä_2)=2\ä=\gamma_8-\gamma_0$.
Since $\ä$ is a long root this is not possible and we have (\ref{2achtens}.R), i.e.
$\LL^\prime+2\ä_1+2\ä_2=-\LL+\delta_8$.

Thinking for a moment gives that (\ref{2siebentens}.L) implies $\gamma_0=-\ä_1$ and
(\ref{2sechstens}.L) implies $\gamma_0=-\ä_2$. On the other hand
(\ref{2viertens}.L) implies $\gamma_0\not=-\ä_1$ and (\ref{2zweitens}.L) implies $\gamma_0\not=-\ä_2$.
Hence  (\ref{2siebentens}.L) entails (\ref{2sechstens}.R) and (\ref{2viertens}.R), as well as
(\ref{2sechstens}.L) entails (\ref{2siebentens}.R) and (\ref{2zweitens}.R).

Now we suppose that (\ref{2siebentens}.L) is satisfied. Then we have (\ref{2achtens}.R),
(\ref{2sechstens}.R) and (\ref{2viertens}.R), i.e.
\[\ä_2=\delta_8-\delta_7\mbox{ and } 2\ä_1=\delta_8-\delta_4.\]
Again because of the extremality of $\LL$ these roots are not proportional.
It implies $2\ä_1-\ä_2=\delta_7-\delta_4$. Now  $\ä_1\bot\ä_2$ and $\delta_7\not=\ä_1,\not=\ä_1-\ä_2$
(Extremality of $\LL$) gives a contradiction.

In the same way we argue supposing that (\ref{2sechstens}.L) holds.

Hence we have shown that neither (\ref{2sechstens}.L) nor (\ref{2siebentens}.L) can be satisfied.
Thus we have (\ref{2sechstens}.R) and (\ref{2siebentens}.R). These together with (\ref{2achtens}.L)
are no contradiction, but if one supposes one of the remaining (\ref{2erstens}.R),
(\ref{2drittens}.R) or (\ref{2fuenftens}.R) we get a contradiction. Hence
(\ref{2erstens}.L),
(\ref{2drittens}.L) and (\ref{2fuenftens}.L) must be valid. But from these together with
(\ref{2nulltens}.L) we derive as above a contradiction.

If we start with the right hand side of the first alternative one proceeds analogously.
\end{description}
All in all we have shown, that the assumption of a long root with (\ref{indirekt}) leads to a contradiction.
\eprf

Now we are in a position that we can use results of \cite{schwachhoefer2} explicitly. First we will cite them.
\bs\cite{schwachhoefer2} \label{s2schwachhoefer}
Let $\lag \subset \laso(N, \ccc ) $
be an irreducible representation of real type of a  complex simple Lie algebra different from
$\mf{sl}(2,\ccc)$. Then it holds:
\bnum
\item If there is an extremal spanning $(\LL_1,\LL_2,\ä)$ triple then
there is no weight $\lam$ for which exists a pair of orthogonal long roots
$\eta_1$ and $\eta_2$ such that
$\left|\frac{2\la\lam,\eta_i\ra}{\|\eta_i\|^2}\right|= 2$.
\item
If furthermore all roots have the same length, then there is no weight $\lam$ for which
exists a triple of orthogonal roots
$\eta_1\bot\eta_2\bot\eta_3\bot\eta_1$ such that
$\left|\frac{2\la\lam,\eta_1\ra}{\|\eta_1\|^2}\right|= 2$ and
$\left|\frac{2\la\lam,\eta_2\ra}{\|\eta_2\|^2}\right|= \left|\frac{2\la\lam,\eta_3\ra}{\|\eta_3\|^2}\right|= 1$.
\enum
\es

We will now show that existence of such a pair or triple of roots implies that (SII) defines
an extremal spanning pair.

\bs \label{s2satz3}
Let $\lag \subset \laso(N, \ccc ) $
be an irreducible representation of real type of a  complex simple Lie algebra different from
$\mf{sl}(2,\ccc)$, with $0\not\in\W$ and satisfying (SII). Then it holds:
If there is a pair of orthogonal long roots
$\eta_1$ and $\eta_2$ such that
$\left|\frac{2\la\LL,\eta_i\ra}{\|\eta_i\|^2}\right|= 2$ for the extremal weight $\LL$ from the property $(SII)$,
then $\LL-\ä$ is an extremal weight, i.e. (SII) defines an extremal spanning triple.
\es

\bprf
Again we argue indirectly considering three different cases for the root $\ä$ from the property (SII)
\begin{description}
\item{{\em Case 1: All roots have the same length or $\ä$ is a long root.}}
Again by applying the Weyl group
the indirect assumption implies that there
is an extremal weight $\LL^\prime$ and a root long $\bb$ orthogonal to $\ä$ such that
$\frac{2\la\LL^\prime,\ä\ra}{\|\ä\|^2}=\frac{2\la\LL^\prime,\bb\ra}{\|\bb\|^2}=-2$.

This gives $\LL^\prime +k\ä+l\bb\in \W$ for $0\le k,l\le 2$
and hence
$
\LL^\prime + k \ä +l\bb\in \W_\ä$ for $ 0\le k\le 1, 0\le l\le 2.$

Among others (SII) implies the existence of $\gamma_i$ and $\delta_i$ from $\D_0$ for $i=0,\ldots, 3$ such that
that the following alternatives must hold
\begin{align}
&&(L)&&&&&&&&(R)&&\nonumber\\
\LL^\prime+\ä&&=&&\LL+\gamma_0&&\mbox{ or }&&\LL^\prime&&=&&-\LL+\delta_0 \label{3nulltens}\\
\LL^\prime+2\ä&&=&&\LL+\gamma_1&&\mbox{ or }&&\LL^\prime+\ä&&=&&-\LL+\delta_1 \label{3erstens}\\
\LL^\prime+\ä+2\bb&&=&&\LL+\gamma_2&&\mbox{ or }&&\LL^\prime+2\bb&&=&&-\LL+\delta_2 \label{3zweitens}\\
\LL^\prime+ä_1+2\ä+2\bb&&=&&\LL+\gamma_3&&\mbox{ or }&&\LL^\prime+\ä+2\bb&&=&&-\LL+\delta_3 .\label{3drittens}
\end{align}
Supposing again (\ref{3nulltens}.L) we conclude that (\ref{3zweitens}.L) and (\ref{3drittens}.L) can not hold
because $\bb$ is long and the extremality of $\LL$. Hence
(\ref{3zweitens}.R) and (\ref{3drittens}.R) must be satisfied. Again the extremality of $\LL$ prevents
that (\ref{3erstens}.R) can be valid. Hence we have (\ref{3erstens}.L).

Now (\ref{3nulltens}.L) gives
\[\frac{2\la\LL,\ä\ra}{\|\ä\|^2}= \frac{2\la\LL^\prime,\ä\ra}{\|\ä\|^2}+2 -\frac{2\la\gamma_0,\ä\ra}{\|\ä\|^2}
=-\frac{2\la\gamma_0,\ä\ra}{\|\ä\|^2}\]
by assumption.

On the other hand (\ref{3nulltens}.L) together with (\ref{3erstens}.L) and
(\ref{3zweitens}.R) and (\ref{3drittens}.R) implies that
 $\ä=\gamma_1-\gamma_0=\delta_3-\delta_2$. We note that
$\gamma_0$ can not be equal to $0$ and $\gamma_1$ not equal to $\ä$.

If $\gamma_0=-\ä$ and $\gamma_1=0$ then $\Lambda =\LL^\prime+2\ä$. Then
(\ref{3zweitens}.R) and (\ref{3drittens}.R) imply
\barr{rcccl}
\la\delta_2,\ä\ra&=& 2\la\LL^\prime,\ä\ra +2\|\ä\|^2&=&0\;\mbox{ and}\\
\la\delta_3,\ä\ra&=& 2\la\LL^\prime,\ä\ra +3\|\ä\|^2&=&\|\ä\|^2.
\earr
Since $\ä$ is long this entails $\delta_2=0$ and $\delta_3=\ä$.
Taking now (\ref{3nulltens}.L) and (\ref{3zweitens}.R) together we get
$\LL=\ä-\bb$. But this forces $0\in\W$ which was excluded.

Thus we have $\ä=\gamma_1-\gamma_0$ with non-proportionality. But this implies, since $\ä$ is long, that
$\frac{2\la\gamma_0,\ä\ra}{\|\ä\|^2}=-1$ and hence $\frac{2\la\LL,\ä\ra}{\|\ä\|^2}= 1$. But this means that
$\LL-\ä$ is an extremal weight.

\item{{\em Case 2: There are roots with different length and $\ä$ is a short root.}}
By assumption there is  a short root $\gamma$ such that
$\gamma\bot \ä$ and $\eta:=\ä+\gamma$ is a long root and an extremal weight $\LL^\prime$ and
a long root $\bb$ such that
$\frac{2\la\LL^\prime,\eta\ra}{\|\eta\|^2}=\frac{2\la\LL^\prime,\bb\ra}{\|\bb\|^2}=-2$.
Analogously to the previous theorem the orthogonality of $\ä$ and $\gamma$ gives
\[-2=\einhalb\left(
\frac{2\la\LL^\prime,\ä\ra}{\|\ä\|^2}+\frac{2\la\LL^\prime,\gamma\ra}{\|\gamma\|^2}\right).\] Hence we have
to consider three cases:
\bnum
\item[(a)] $\frac{2\la\LL^\prime,\ä\ra}{\|\ä\|^2}=\frac{2\la\LL^\prime,\gamma\ra}{\|\gamma\|^2}-2$,
\item[(b)] $\frac{2\la\LL^\prime,\ä\ra}{\|\ä\|^2}=-3$ and $\frac{2\la\LL^\prime,\gamma\ra}{\|\gamma\|^2}-1$,
\item[(c)] $\frac{2\la\LL^\prime,\ä\ra}{\|\ä\|^2}=-1$ and $\frac{2\la\LL^\prime,\gamma\ra}{\|\gamma\|^2}-3$.
\enum
Then an easy calculation shows that $\la\ä,\beta\ra=\la\gamma,\bb\ra=0$ in each case.

We shall consider the case (a),(b) and (c) separately.
\begin{description}
\item{{\em Case (a):}}
Here we can proceed completely analogously to the first case 1. We have that
$
\LL^\prime + k \ä +l\bb\in \W_\ä$ for $ 0\le k\le 1, 0\le l\le 2$ leading to the same
set of equations (\ref{3nulltens}) --- (\ref{3drittens}) and the same implications since $\bb$ is long again.
The proportional case is excluded as above and we get that
$\ä=\gamma_1-\gamma_0$ non proportional. At least one has to be a short root and
$\la \gamma_0,\ä\ra<0$ and $\la \gamma_1,\ä\ra>0$.
On the other hand we have
$\frac{2\la\LL,\ä\ra}{\|\ä\|^2}=-\frac{2\la\gamma_0,\ä\ra}{\|\ä\|^2}$
and $\frac{2\la\LL,\ä\ra}{\|\ä\|^2}=-\frac{2\la\gamma_1,\ä\ra}{\|\ä\|^2}+2$ by (\ref{3nulltens}.L) and
(\ref{3erstens}.L). But this implies that both are short and
$\frac{2\la\LL,\ä\ra}{\|\ä\|^2}=1$ which is the proposition.

\item{{\em Case (b):}}
$\frac{2\la\LL^\prime,\ä\ra}{\|\ä\|^2}=-3$ implies 
$\LL^\prime + k \ä +l\bb\in \W_\ä$ for $0\le k,l\le 2$. (SII) then implies
\begin{align}
&&(L)&&&&&&&&(R)&&\nonumber\\
\LL^\prime+\ä&&=&&\LL+\gamma_0&&\mbox{ or }&&\LL^\prime&&=&&-\LL+\delta_0 \label{4nulltens}\\
\LL^\prime+2\ä&&=&&\LL+\gamma_1&&\mbox{ or }&&\LL^\prime+\ä&&=&&-\LL+\delta_1 \label{4erstens}\\
\LL^\prime+3\ä&&=&&\LL+\gamma_2&&\mbox{ or }&&\LL^\prime+2\ä&&=&&-\LL+\delta_2 \label{4zweitens}\\
\LL^\prime+2\ä+2\bb&&=&&\LL+\gamma_3&&\mbox{ or }&&\LL^\prime+\ä+2\bb&&=&&-\LL+\delta_3 \label{4drittens}\\
\LL^\prime+3\ä+2\bb&&=&&\LL+\gamma_4&&\mbox{ or }&&\LL^\prime+2\ä+2\bb&&=&&-\LL+\delta_3. \label{4viertens}
\end{align}
Supposing (\ref{4nulltens}.L) excludes (\ref{4drittens}.L) and (\ref{4viertens}.L) because $\bb$ is long. Hence
(\ref{4drittens}.R) and (\ref{4viertens}.R) are valid and exclude (\ref{4erstens}.R) and (\ref{4zweitens}.L).
Hence (\ref{4erstens}.L) and (\ref{4zweitens}.L) are satisfied. This gives
$\ä=\gamma_2-\gamma_1=\gamma_1-\gamma_0$ with  $\gamma_0$ different from $0$ and $-\ä$,
$\gamma_1$ different from $0$ and $\ä$ and $\gamma_2$ different from $\pm\ä$.
Hence $\ä+\pm\delta_1$ is a long root with $\ä\bot\delta_1$. But this gives
$\frac{2\la\LL,\ä\ra}{\|\ä\|^2}= \frac{2\la\LL^\prime,\ä\ra}{\|\ä\|^2}+4=1$, i.e. $\LL-\ä$ is an extremal weight.

\item{{\em Case (c):}}
Here we have that $\LL^\prime +k\gamma +l\bb\in\W_\ä$ for $0\le k\le 3$ and $0\le l\le 2$
since $\frac{2\la \LL^\prime +k\gamma +l\bb,\ää\ra}{\|\ä\|^2}=-1$.
The equations implied by (SII) lead easily to a contradiction:
\begin{align}
&&(L)&&&&&&&&(R)&&\nonumber\\
\LL^\prime+\ä&&=&&\LL+\gamma_0&&\mbox{ or }&&\LL^\prime&&=&&-\LL+\delta_0 \label{5nulltens}\\
\LL^\prime+\ä+3\gamma&&=&&\LL+\gamma_1&&\mbox{ or }&&\LL^\prime+3\gamma&&=&&-\LL+\delta_1 \label{5erstens}\\
\makebox[2.5cm][r]{$\LL^\prime+\ä+2\beta+3\gamma$}&&=&&\LL+\gamma_2&&\mbox{ or }&&\LL^\prime+2\bb+3\gamma&&=&&-\LL+\delta_2. \label{5zweitens}
\end{align}
Supposing (\ref{5nulltens}.L) excludes (\ref{5erstens}.L) and (\ref{5zweitens}.L). Hence
(\ref{5erstens}.R) and (\ref{5zweitens}.R) are valid but contradict to each other because $\bb$ is long.
\end{description}
\end{description}
\eprf

\bs \label{s2satz3a}
Let $\lag \subset \laso(N, \ccc ) $
be an irreducible representation of real type of a  complex simple Lie algebra different from
$\mf{sl}(2,\ccc)$, with $0\not\in\W$ and satisfying (SII).
If furthermore all roots have the same length, and if there is a triple of orthogonal roots
$\eta_1\bot\eta_2\bot\eta_3\bot\eta_1$ such that
$\left|\frac{2\la\LL,\eta_1\ra}{\|\eta_1\|^2}\right|= 2$ and
$\left|\frac{2\la\LL,\eta_2\ra}{\|\eta_2\|^2}\right|= \left|\frac{2\la\LL,\eta_3\ra}{\|\eta_3\|^2}\right|= 1$
then holds one of the cases
\bnum
\item $\LL-\ä$ is an extremal weight, i.e. (SII) defines an extremal spanning triple, or
\item $\LL=\ä+\einhalb(\bb+\gamma)$ with roots $\ä\bot\bb\bot\gamma\bot\ä$.
\enum
\es

\bprf
Let $\ä$ be the root determined by (SII).
The assumption implies that there is an extremal weight $\LL^\prime$ and roots $\bb$ and $\gamma$ such that
\[
\frac{2\la\LL^\prime,\ä\ra}{\|\ä\|^2}= -2\;\mbox{ and }\;
\frac{2\la\LL,\bb\ra}{\|\bb\|^2}= \frac{2\la\LL,\gamma\ra}{\|\gamma\|^2}= -1.\]
Then
$\LL^\prime +k\ä+l\bb+m\gamma\in \W_ä$ for $k,l,m=0,1$. Hence
(SII) implies again
\begin{align}
&&(L)&&&&&&&&(R)&&\nonumber\\
\LL^\prime+\ä&&=&&\LL+\gamma_0&&\mbox{ or }&&\LL^\prime&&=&&-\LL+\delta_0 \label{6nulltens}\\
\LL^\prime+2\ä&&=&&\LL+\gamma_1&&\mbox{ or }&&\LL^\prime+\ä&&=&&-\LL+\delta_1 \label{6erstens}\\
\LL^\prime+\ä+\bb&&=&&\LL+\gamma_2&&\mbox{ or }&&\LL^\prime+\bb&&=&&-\LL+\delta_2 \label{6zweitens}\\
\LL^\prime+ 2\ä+\bb&&=&&\LL+\gamma_3&&\mbox{ or }&&\LL^\prime+\ä+\bb&&=&&-\LL+\delta_3 \label{6drittens}\\
\LL^\prime+\ä+\gamma&&=&&\LL+\gamma_4&&\mbox{ or }&&\LL^\prime+\gamma&&=&&-\LL+\delta_4 \label{6viertens}\\
\LL^\prime+2\ä+\gamma&&=&&\LL+\gamma_5&&\mbox{ or }&&\LL^\prime+\ä+\gamma&&=&&-\LL+\delta_5 \label{6fuenftens}\\
\LL^\prime+\ä+\bb+\gamma&&=&&\LL+\gamma_6&&\mbox{ or }&&\LL^\prime+\bb+\gamma&&=&&-\LL+\delta_6 \label{6sechstens}\\
\LL^\prime+2\ä+\bb+\gamma&&=&&\LL+\gamma_7&&\mbox{ or }&&\LL^\prime+\ä+\bb+\gamma&&=&&-\LL+\delta_7. \label{6siebentens}
\end{align}
Supposing (\ref{6nulltens}.L) excludes (\ref{6siebentens}.R) because of the orthogonality of the roots.
Thus it must hold (\ref{6siebentens}.R). Now we consider two cases:
\begin{description}
\item{{\em Case 1: $\la\gamma_0,\bb\ra=\la\gamma_0,\gamma\ra=0$.}}
This excludes (\ref{6zweitens}.L), (\ref{6viertens}.L) and (\ref{6sechstens}.L) and implies therefore
(\ref{6zweitens}.R), (\ref{6viertens}.R) and (\ref{6sechstens}.R). The latter together with
(\ref{6siebentens}.R) gives $\ä=\delta_7-\delta_6$.

Since $\delta_7\not=0$ this implies
$\la\ä,\delta_7\ra >0$.
On the other hand
(\ref{6siebentens}.R) and the assumption gives
$\frac{2\la\LL,\ä\ra}{\|\ä\|^2}=\frac{2\la\delta_7,\ä\ra}{\|\ä\|^2}>0$. If $\ä\not=\delta_7$ we are done.

But $\delta_7=\ä$ implies
$\LL^\prime+\bb+\gamma=-\LL=-\LL^\prime-\ä-\gamma_0$ and hence
$-2=2-2-\frac{2\la\ä,\gamma_0\ra}{\|\ä\|^2}$, i.e. $\gamma_0=-\ä$. Taking everything together we
get $2\LL=2\ä-(\bb+\gamma)$. 

\item{{\em Case 2: $\la\gamma_0,\bb\ra$ or $\la\gamma_0,\gamma\ra$ not equal to zero.}}
This implies $\gamma_0\not=\pm \ä$ and thus
$\frac{2\la\LL,\ä\ra}{\|\ä\|^2}=-\frac{2\la\gamma_0,\ä\ra}{\|\ä\|^2}=\pm 1$ or zero.
Now (\ref{6erstens}.L) would
imply $\ä=\gamma_1-\gamma_0$, i.e. $\la\ä,\gamma_0\ra< 0$. This would be the proposition.

Hence we suppose
(\ref{6erstens}.R). This together with the starting point
(\ref{6nulltens}.L)  gives
\be
\LL&=&\einhalb\left(\delta_1-\gamma_0\right)\;\;\mbox{ and}\\
\LL^\prime&=&-\ä +\einhalb (\delta_1+\gamma_0).
\ee
The second equation implies using the  assumption that $\la\ä,\delta_1+\gamma_0\ra=0$. For the length of
both extremal weights then holds
\be
\|\LL\|^2&=& \frac{1}{4}\left( \|\delta_1\|^2+\|\gamma_0\|^2 - 2\la\delta_1,\gamma_0\ra\right)\\
\|\LL^\prime\|^2&=& \|\ä\|^2 - \underbrace{\la\ä,\delta_1+\gamma_0\ra}_{=0}+
\frac{1}{4}\left( \|\delta_1\|^2+\|\gamma_0\|^2 + 2\la\delta_1,\gamma_0\ra\right).
\ee
This gives $0=\|\ä\|^2+\la\delta_1,\gamma_0\ra$. Since all roots have the same length this implies $
\delta_1=-\gamma_0$. Hence $\LL$ is a root. But this was excluded.
\end{description}
\eprf

Now using the proposition \ref{s2schwachhoefer} of Schwachh\"ofer we get a corollary.

\bfolg \label{sf}
Let $\lag \subset \laso(N, \ccc ) $
be an irreducible representation of real type of a  complex simple Lie algebra different from
$\mf{sl}(2,\ccc)$, with $0\not\in\W$ and satisfying (SII). Then it holds:
\bnum
\item
\label{sf1}
There is no a pair of orthogonal long roots
$\eta_1$ and $\eta_2$ such that
$\left|\frac{2\la\LL,\eta_i\ra}{\|\eta_i\|^2}\right|= 2$ for the extremal weight $\LL$ from the property $(SII)$.
\item\label{sf2}
If furthermore all roots have the same length, and if there is a triple of orthogonal roots
$\eta_1\bot\eta_2\bot\eta_3\bot\eta_1$ such that
$\left|\frac{2\la\LL,\eta_1\ra}{\|\eta_1\|^2}\right|= 2$ and
$\left|\frac{2\la\LL,\eta_2\ra}{\|\eta_2\|^2}\right|= \left|\frac{2\la\LL,\eta_3\ra}{\|\eta_3\|^2}\right|= 1$
then $\LL=\ä+\einhalb(\bb+\gamma)$ with roots $\ä\bot\bb\bot\gamma\bot\ä$.
\enum
\efolg

Before we apply this corollary we have to deal with the remaining exception in the second point.

\blem\label{ausnahmelemma}
If the representation of a simple Lie algebra with roots of the same length
has an extremal weight $\LL$ such that
$\LL=\ä+\einhalb(\bb+\gamma)$ with roots $\ä\bot\bb\bot\gamma\bot\ä$. Then it holds
\bnum
\item
There is no root $\delta$ such that $\la\delta,\bb\ra=0$, $\la\delta,\gamma\ra\not=0$ and $\delta\not\sim\gamma$.
\item
The root system is $D_n$ and the representation has one of the following highest weights:
$\w_3$ for arbitrary $n$, $\w_1+\w_3$ or $\w_1+\w_4$ for $n=4$ and $\w_2$ for $n=3$.
\enum\elem

\bprf
The first point is easy to see: If there is such a $\delta$ then we have
\[\frac{2\la\LL,\delta\ra}{\|\delta\|^2}=
\frac{2\la\ä,\delta\ra}{\|\delta\|^2}+\einhalb
\frac{2\la\gamma,\delta\ra}{\|\delta\|^2}=
\frac{2\la\ä,\delta\ra}{\|\delta\|^2}\pm\einhalb\not\in\mathbb{Z}.\]
This is a contradiction.

Now we consider the different root systems with roots of constant length.
\begin{description}
\item{$A_n$:} Here the assumption means that $\LL=e_i-e_j+\einhalb
\left(e_p-e_q+e_r-e_s\right)$ with all indices different from each other. But then
$\frac{2\la\LL,e_i-e_p\ra}{\|e_i-e_p\|^2}$ is not an integer.
\item{$D_n$:}
If $\ä=e_i\pm e_j$, $\bb=e_p\pm e_q$ and $\gamma=e_r\pm e_s$ with all indices different we get the
same contradiction as in the
$A_n$ case. Thus we are left with two cases.

The first is $\bb+\gamma=e_p+e_q+e_p-e_q=2e_p$ and hence
$\LL=e_i\pm e_j +e_p$. This leads to $\LL=\w_3$ or for $n=3$ to $\LL=\w_2$.

The second is $\ä=e_i+e_j$, $\bb=e_i-e_j$ and $\gamma=e_p\pm e_q$. For $n>4$ we found
a root $e_p+e_s$ which leads to a contradiction by applying the first point. For $n=4$ we have
$\LL=\frac{3}{2}e_i+\einhalb\left(e_j+e_p\pm e_q\right)$. But this yields the remaining representations.
\item{$E_6$:} $E_6 $ has two different types of roots:
\[
e_i\pm e_j\;\;\mbox{ and }\;\;
\einhalb\big( e_8-e_7-e_6 \ \underbrace{\pm e_5\pm e_4\pm e_3\pm e_2\pm e_1}_{\makebox[2cm][c]{even number of minus signs}}\big).
\]
The only possibility for $\bb$ and $\gamma$ for which the first point yields no contradiction
is $\bb=e_i+e_j$ and $\gamma=e_i-e_j$. Hence $ \LL=\ä+e_i$.  $\ä\bot\bb$ and $\ä\bot\gamma
$ implies $\ä=e_p+e_q$. But then $\la\LL, \einhalb( \ldots )\ra\not\in \mathbb{Z}$.
\end{description}
Proceeding analogously for $E_7$ and $E_8$ we prove the second assertion.
\eprf

Now using all these properties we can find the representations without weight zero and
satisfying (SII).

\bs \label{s2satz4}
Let $\lag \subset \laso(N, \ccc ) $
be an irreducible representation of real type of a  complex simple Lie algebra different from
$\mf{sl}(2,\ccc)$, with $0\not\in\W$ and satisfying (SII).
Then the roots system and the highest weight of the representation is is one of the following
(modulo congruence):
\begin{description}
\item{$A_n$:} $\w_4$ for $n=7$.
\item{$B_n$:} $\w_n$ for $n=3,4,7$.
\item{$D_n$:} $\w_1$, $2\w_1$ for arbitrary $n$ and $\w_8$ for $n=8$.
\end{description}
\es

\bprf
We apply proposition \ref{s2satz2} and corollary \ref{sf} to the remaining representations with $0\not\in\W$,
i.e. representations of $A_n,\ B_n,\ C_n,\ D_n,\  E_6$ and $E_7$.
Therefore we use a fundamental system such that the extremal weight $\LL$ determined by (SII)
is the highest weight. It can be written in the fundamental representations $\LL=\sumk m_k\w_k$ with
$m_k\in\mathbb{N}\cup\{0\}$.
\begin{description}
\item{$A_n$:}
Proposition \ref{s2satz2} gives for the largest root
\[2\ \ge\ \frac{2\la\LL,e_1-e_{n+1}\ra}{\|e_1-e_{n+1}\|^2}\ =\ \sumk m_k\la\w_k,e_1-e_{n+1}\ra\ =\  \sumk m_k.\]
Since the representation has to be self dual we have that $m_i=m_{n+1-i}$.

First we consider the case that $\LL=\w_i+\w_{n+1-i}$. For $n>2$ we get in case $i>1$ that
$\la\LL,e_2-e_n\ra=2$. But $(e_2-e_n)\bot(e_1-e_{n+1})$ gives a contradiction to
\ref{sf1} of corollary \ref{sf}.  For $n\ge 2$ it has to be
\[\LL=\w_1+\w_n=2e_1 + e_2 +\ldots e_n=e_1-e_{n+1}\]
recalling that for $A_n$ holds that $e_1=-(e_2+\ldots +e_{n+1})$. Thus the representation is the adjoint one with
$0\in \W$.

Now we consider the case that $n+1$ is even and $\LL=2\w_{\frac{n+1}{2}}$. This representation is orthogonal
but again we have $\la\LL,e_2-e_n\ra=2$ for $n>2$. But this is impossible because of point \ref{sf1}
of corollary \ref{sf}.

For $n+1$ even we have to study the case $\LL=\w_{\frac{n+1}{2}}$. This representation is orthogonal if
$\frac{n+1}{2}$ is even. The weights of this representation are given by
$\pm e_{k_1}\pm \ldots \pm e_{k_{\frac{n+1}{2}}}$ where the $\pm$'s are meant to be independent of each other.

We will show that (SII) implies  $n\le 7$.

Hence suppose that there is a root $\ä$ such that (SII) with $\LL$.  We have to consider two cases for $\ä$.
The first is that $\ä=e_i-e_j$ with $1\le i\le\frac{n+1}{2}<j\le n+1$. W.l.o.g. we take
$\ä=e_{\frac{n+1}{2}}-e_{\frac{n+1}{2}+1}$ and consider the  weight
\[\lam:=e_1+\ldots e_{\frac{n+1}{2}-3}+e_{\frac{n+1}{2}+1}+e_{\frac{n+1}{2}+2}+e_{\frac{n+1}{2}+3}.\]
$\la\lam,\ä\ra<0$ implies $\lam\in\W_\ä$. Then
$\lam-(\LL-\ä)\in \Delta_0$ or $\lam+\LL\in\Delta_0$.
We check the first alternative:
$\LL-\ä=e_1+\ldots e_{\frac{n+1}{2}-1}+e_{\frac{n+1}{2}+1}$ implies
\[
\lam-(\LL-\ä)
=e_{\frac{n+1}{2}-3}+e_{\frac{n+1}{2}-2}+e_{\frac{n+1}{2}+2}+e_{\frac{n+1}{2}+3}.\]
But this is not a root.

For the second alternative we get, recalling that $-e_1=e_2+\ldots + e_{n+1}$,
\[
\lam+\LL
=e_1+\ldots +e_{\frac{n+1}{2}-3}-e_{\frac{n+1}{2}+4}-\ldots -e_{n+1}.\]
This is not a root if $\frac{n+1}{2}>4$, i.e. $n>7$.

For the second type of root $\ä=e_i-e_j$ with $1\le i< j\le \frac{n+1}{2}$ and
$\frac{n+1}{2}<i<j\le n+1$ one derives analogously that $n\le 5$.

Hence for $\LL=\w_{\frac{n+1}{2}}$ the property (SII) can only be fulfilled if $n\le 7$. These representations
are orthogonal for $n=7$ and $n=3$. $A_3 $ is isomorphic to $D_3$ and the representation with highest weight $\w_2$ of $A_3$ is equivalent
to the one with $\w_1$ of $D_3$.

\item{$B_n$:}
Again proposition \ref{s2satz2} gives for the largest root
\[2\ \ge\ \frac{2\la\LL,e_1+e_2\ra}{\|e_1+e_2\|^2}\ =\ \sumk m_k\la\w_k,e_1+e_2\ra\ =\  m_1+2m_2+\ldots 2m_{n-1}+m_n.\]
The only representations with $0\not\in \W$ are these with $\LL=\w_1+\w_n$ and the spin representation $\LL=\w_n$. There is no
possibility to apply the first point of corollary \ref{sf}. But we verify that for $\LL=\w_1+\w_n$
(SII) implies  $n\le 2$ and for the spin representation $\LL=\w_n$ (SII) implies $n\le 7$.

\begin{description}
\item{{\em The spin representations:}} For these we show that
(SII) implies $n\le 7$
The spin representation of highest weight
$\LL=\einhalb(e_1+\ldots +e_n)$ has weights
$\W=\left\{\einhalb(\ve_1 e_1 +\ldots +\ve_n e_n)|\ve_i=\pm 1 \right\}$. We have to consider three types for the
root $\ä$: $\ä=e_i$, $\ä=e_i+e_j$ and $\ä=e_i-e_j$.

For the first we can assume w.l.o.g. that $\ä=e_1$. Then
$\W_\ä=\{\einhalb(-e_1+\ve_2 e_2 +\ldots +\ve_n e_n)|\ve=\pm 1\}$. It is
$\LL-\ä=\einhalb(-e_1+e_2+\ldots +e_n)$. Hence for $\lam\in \W_\ä$ we have
\be
\LL-\ä-\lam&=&\einhalb((1-\ve_2)e_2 + \ldots + (1-\ve_n)e_n  )\;\mbox{ and }\\
\LL+\lam&=&\einhalb((1+\ve_2)e_2 + \ldots +(1+\ve_n)e_n  )
\ee
If (SII) is satisfied at least one of these expression has to be a root. But if $n\ge 7$
 we can choose $(\ve_2, \ldots \ve_n)$ such that non of them is a root.

The second type of root shall be w.l.o.g. $\ä=e_1-e_2$. In this case
$\W_\ä=\{\einhalb(-e_1+e_2+\ve_3 e_3 +\ldots +\ve_n e_n)|\ve_i=\pm 1\}$ and
$\LL-\ä=\einhalb(-e_1+2e_2+e_3+\ldots +e_n)$. Hence for $\lam\in\W_\ä$
\be
\LL-\ä-\lam&=&\einhalb(e_2+(1-\ve_3)e_3 + \ldots +(1-\ve_n)e_n  )\;\mbox{ and }\\
\LL+\lam&=&\einhalb(2e_2+(1+\ve_3)e_3 + \ldots +(1+\ve_n)e_n  )
\ee
We can choose $\lam$ such that none of them is a roots if $n\ge 4$.

Now we consider the last type of root, $\ä=e_1+e_2$.
$\W_\ä=\{\einhalb(-e_1-e_2+\ve_3 e_3 +\ve_n e_n)|\ve_i=\pm 1\}$ and
$\LL-\ä=\einhalb(-e_1-e_2+e_3+\ldots +\ldots +e_n)$. Hence for $\lam\in\W_\ä$
\be
\LL-\ä-\lam&=&\einhalb((1-\ve_3)e_3 + \ldots +(1-\ve_n)e_n  )\;\mbox{ and }\\
\LL+\lam&=&\einhalb((1+\ve_3)e_3 + \ldots +(1+\ve_n)e_n  )
\ee
We can choose $\lam$ such that none of them is a roots if $n\ge 8$.
Hence if (SII) is satisfied it has to be $n\le 7$ and for $n=7$ the pair of property (SII) is
of the shape $(\LL, e_1+e_2)$.

Now for $n=2$, $n=5$ and $n=6$ the spin representations are symplectic but not orthogonal.

\item{{\em The representations of $\LL=\w_1+\w_n=\frac{3}{2}e_1 +\einhalb(e_2+\ldots +e_n)$.}}
 Then the weights are given by
$\einhalb( a e_{k_1} + \ve_2 e_{k_2} + \ldots +\ve_n e_{k_n})$ with $a\in\{\pm1,\pm3\}$ and $\ve_i=\pm 1$.
For these one shows analogously that (SII) implies  $n\le 2$. For $n=2$ this representation is symplectic.
\end{description}
\item{$C_n$:}
For the largest root we get
\[2\ \ge\ \frac{2\la\LL,2e_1\ra}{\|2e_1\|^2}\ =\ \sumk m_k\la\w_k,e_1\ra\ =\  \sumk m_k.\]
In case that one $m_i=2$ and all others zero we have that $0\in \W$. Hence we suppose that
$\LL=\w_i+\w_j$ for $i\not=j$. If $i>1$ we get for the root $2e_2$ which is orthogonal to $2e_1$
that
$\frac{2\la\LL,2e_2}{\|2e_2\|^2}=2$. Thus by \ref{sf1} of corollary \ref{sf} we have $i=1$.
But $\LL=\w_1 +\w_i$ is only orthogonal if $i$ is odd, but if $i$ is odd we have that $0\in \W$.

Hence we have to deal with the case $\LL=\w_i$. This is orthogonal if $i$ is even, but in this case $0\in \W$.

\item{$D_n$:}
Here we get for the largest root
\[2\ \ge\ \frac{2\la\LL,e_1+e_2\ra}{\|e_1+e_2\|^2}\ =\ \sumk m_k\la\w_k,e_1+e_2\ra\ =\
m_1 +2m_2+\ldots +2m_{n-2}+m_{n-1}+m_{n} .\]
First we consider the representation where this number is equal to $2$.

For the representations $2\w_n$ and $2\w_{n-1}$ we have that $0\in \W$.

For the representations
$\LL=\w_1+\w_n$ and $\LL=\w_1+\w_{n-1}$ we get that $n=4$ or there is no triple as in the second point of
proposition \ref{sf}. Thus suppose in this case $n>4$. We have
that$\la\LL,e_1+e_2\ra=2$ and for the orthogonal roots $\la\LL,e_1-e_2\ra=\la\LL ,e_3\pm e_4\ra=1$.
But this contradicts
proposition \ref{sf},\ref{sf1}.

For $\LL=\w_{n-1}+\w_n=e_1+\ldots +e_{n-1}$ we have that $0\not\in \W$ implies
$n-1$ even. The first point of corollary \ref{sf} then gives for $n>4$ that
$2=\la\LL,e_3+e_4\ra$ which is impossible. Hence $n\le 4$. Then
$1=\la \LL,e_3\pm e_4\ra$ and the second point of corollary \ref{sf} imply $n\le3$.

Now suppose that $\LL=\w_i$ for $2\le i\le n-2$. We apply the first point of corollary \ref{sf}.
If $n\ge 4$ we get that $\la\w_i,e_3+e_4\ra=2$ for $i\ge 4$ but this was excluded. Hence $i\le 3$.

In the case $n=3$ we have that only $\w_2$ is an orthogonal representation. But for this holds that $0\in \W$.

Thus, to get the assertion of the proposition we have to show that
\bnum
\item For the spin representations $\LL=\w_{n-1}$ and $\LL=\w_n$ (SII) implies $n\le 8$
\item $\LL=\w_3$ does not satisfy (SII),
\item $\LL=\w_1+\w_3$ and $\w_1+\w_4$ for $n=4$ do not satisfy (SII).
\enum

\begin{description}
\item{{\em The spin representations:}} For these we show that
(SII) implies first $n\le 8$. Because we are interested in the representations
modulo congruence it suffices to consider the spin representation of highest weight
$\LL=\einhalb(e_1+\ldots +e_n)$ with weights
$\W=\left\{\einhalb(\ve_1 e_1 +\ve_n e_n)|\ve_i=\pm 1 \mbox{ and $\ve_i=-1$ for an even number}\right\}$.

Analogously as for $B_n$ we get for two types of
roots  $\ä=e_i+e_j$ and $\ä=e_i-e_j$ that (SII) implies $n\le 8$ (We have to admit one dimension higher
because of the sign restriction of the weights).

Now for $n$ odd the spin representation is not self dual, and for $n=6$ not orthogonal. For $n=4$ it is
congruent to $\w_1$.

\item{$\LL=\w_3=e_1+e_2+e_3$:} Here it is
$\W=\{(\ve_1e_{k_1}+\ve_2e_{k_2}+\ve_3e_{k_3}|\ve_i=\pm 1\}\cup\{\pm e_i\}$.
For $n=3$ and $n=4$ this is a spin representation. Hence suppose $n\ge 5$.

For $\ä=e_1+e_2$ we get $\LL-\ä=e_3$. Set $\lam:=-e_1+e_4+e_5\in\W_\ä$. Hence
$\LL-\ä-\lam=e_3+e_1-e_4-e_5$ and $\LL+\lam=e_2+e_3+e_4+e_5$. None is a root, i.e. $\w_3$ for $n\ge5$ does not
satisfy (SII).

For $\ä=e_1-e_2$ we get the same.

\item{{\em $\LL=\w_1+\w_3$ and $\w_1+\w_4$ for $n=4$.}} These are congruent to each other
and as above it can be shown that they do not satisfy (SII).
\end{description}

\item{$E_6$ and $E_7$:} For these we refer to \cite{schwachhoefer2}. There is shown that under the conclusions
of proposition \ref{s2satz2} and \ref{s2schwachhoefer} --- which is our situation because of lemma
\ref{ausnahmelemma} --- the only remaining representations are the standard representations of
$E_6$ and  $E_7$. But the first is not
self dual and the latter symplectic but not orthogonal.
\end{description}
\eprf

We get the following
\bfolg
Let $\lag\subset\laso(N,\ccc)$ be an orthogonal algebra of real type different from $\mf{sl}(2,\ccc)$.
If $0\not\in\W$ and (SII) is satisfied, then  it is the complexification of a Riemannian holonomy representation
or the spin representation of $\laso(15,\ccc)$.
\efolg

\bprf
We give the Riemannian manifolds the complexified holonomy representation of which is
one of the representations of proposition \ref{s2satz4}.

The representation with highest weight $\w_4$ of $A_7$ is the complexified holonomy representation
of the symmetric space of type $EV$, i.e. of $E_7/SU(8)$ resp. $E_{7(7)}/SU(8)$.

The spin representations of $B_n$ for $n=3,4$ are the holonomy representations of a non-symmetric
$Spin(7)$--manifold and of the symmetric space of type $FII$, i.e. of
$F_4/Spin(9)$ resp. $F_{4(-20)}/Spin(9)$. For $n=7$ we have an exception.

For $D_n$ first we have  the standard representation, i.e. the complexified holonomy representation
of a generic manifold. The representation with highest weight $2\w_1$ is
the complexified holonomy representation of the symmetric space of type $AI$, i.e. of
$SU(2n)/SO(2n,\rr)$, resp. $Sl(2n,\rr)/SO(2n,\rr)$.
The remaining representation of $Spin(16)$ is the complexified holonomy representation
of the symmetric space of type $EVIII$, i.e. of $E_8/Spin(16)$, resp.
$E_{8(8)}/Spin(16)$.
\eprf

\subsection{Consequences for simple weak-Berger algebras of real type}

Before we conclude the result we need a lemma to exclude both exceptions.

\blem
The spin representation of $B_7$  and the representation of $G_2$ with two times a short root
as highest weight are not weak-Berger.
\elem
\bprf
1.) Suppose that the spin representation of $B_7$ is weak-Berger.
We have shown that it does not satisfy the property (SI). Hence it obeys (SII). Let $(\LL,\ä)$ be the pair
of (SII). We choose a fundamental system such that $\LL=\w_7$ is the highest weight.
In the proof of proposition \ref{s2satz4} we have shown that in this case $\ä=e_i+e_j$.

Let now $Q_\phi$ be the weight element from $\bhg$ and $u_\LL\in V_\LL$ such that
$Q_\phi(u_\LL)=A_{e_i+e_j}\in\lag_{e_i+e_j}$.
Since $Q_\phi(u_\LL)\in \lag_{\phi+\LL}$ this implies that
$\phi=e_i+e_j-\LL$ is a weight of $\bhg$. Hence
$\phi=-\einhalb(e_1+\ldots +e_{i-1}-e_i+e_{i+1}+\ldots +e_{j-1}-e_j+e_{j+1}+\ldots+e_7)$ is also an extremal
weight of
$V$ and we can consider a weight vector $u_{-\phi}\in V_{-\phi}$.
For this we get $Q_\phi(u_{-\phi})\in\mf{t}$. In case it does not vanish it would define a
planar spanning triple $(\phi,-\phi,\left(Q_\phi(u_{-\phi})\right)^\bot)$, i.e. (SI) would be satisfied.
But this was not possible, and thus $Q_\phi(u_{-\phi})=0$.

On the other hand we  have that $0\not=Q_\phi(u_\LL)u_{-\phi}\in V_\LL$ and thus there is a $v\in V_{-\LL}$
such that $H(Q_\phi(u_\LL)u_{-\phi},v)\not=0$. Now the Bianchi identity gives
\[0=
H(Q_\phi(u_\LL)u_{-\phi},v)
+\underbrace{H(Q_\phi(u_{-\phi})v,u_\LL))}_{=0}
+
H(Q_\phi(v)u_\LL, u_{-\phi},v)
.\]
Hence $0\not=Q_\phi(v)\in \lag_{\phi-\LL}$. But
$\phi-\LL=-(e_1+\ldots +e_{i-1}+e_{i+1}+\ldots +e_{j-1}+e_{j+1}+\ldots+e_7)$ is not a root, hence
$\lag_{\phi-\LL}=\{0\}$. This is a contradiction.

2.) Suppose that the representation of $G_2$ with two times a short root as highest weight is weak-Berger.
We will argue analogously as for $B_n$.
\\[.2cm]
\begin{minipage}[b]{8cm}{
In the picture we see the weight lattice of this representation
(the arrows represent the roots).
Obviously there is no planar spanning triple, because there
is no hypersurface which contains all but two extremal weight
(see also proof of proposition \ref{wurzelgewicht}).}

The weak-Berger property implies that there is a pair $(\LL,\ä)$
such that (SII) is satisfied. We choose a fundamental system such that
$\LL=2\eta$ is the maximal weight.
\end{minipage}
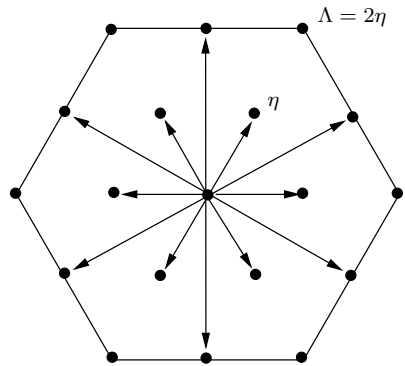
\hfill
\begin{minipage}[b]{5cm}{
\begin{picture}(0,0)%
\includegraphics{g2.pstex}%
\end{picture}%
\setlength{\unitlength}{4144sp}%
\begingroup\makeatletter\ifx\SetFigFont\undefined%
\gdef\SetFigFont#1#2#3#4#5{%
  \reset@font\fontsize{#1}{#2pt}%
  \fontfamily{#3}\fontseries{#4}\fontshape{#5}%
  \selectfont}%
\fi\endgroup%
\begin{picture}(2360,2165)(299,-1643)
\put(1856,-90){\makebox(0,0)[lb]{\smash{\SetFigFont{8}{9.6}{\rmdefault}{\mddefault}{\updefault}$\eta$}}}
\put(2157,417){\makebox(0,0)[lb]{\smash{\SetFigFont{8}{9.6}{\rmdefault}{\mddefault}{\updefault}$\Lambda=2\eta$}}}
\end{picture}}
\end{minipage}\\
Using the realization of $G_2$ from the appendix of \cite{knapp96} we have that
$\eta=e_3-e_2$. Now we have to determine the roots for which (SII) is satisfied.\\
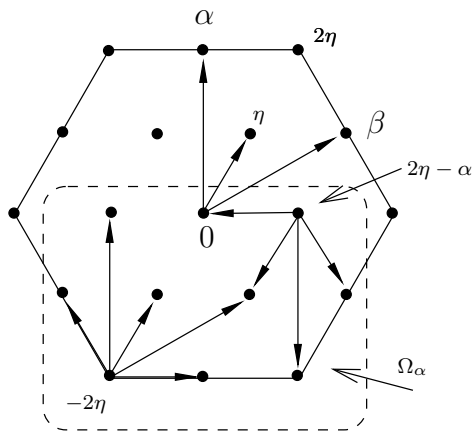
\begin{minipage}[t]{5cm}{
\begin{picture}(0,0)%
\includegraphics{2g2.pstex}%
\end{picture}%
\setlength{\unitlength}{4144sp}%
\begingroup\makeatletter\ifx\SetFigFont\undefined%
\gdef\SetFigFont#1#2#3#4#5{%
  \reset@font\fontsize{#1}{#2pt}%
  \fontfamily{#3}\fontseries{#4}\fontshape{#5}%
  \selectfont}%
\fi\endgroup%
\begin{picture}(2437,2601)(267,-1950)
\put(2611,-1591){\makebox(0,0)[lb]{\smash{\SetFigFont{8}{9.6}{\rmdefault}{\mddefault}{\updefault}$\Omega_\alpha$}}}
\put(2108,387){\makebox(0,0)[lb]{\smash{\SetFigFont{8}{9.6}{\rmdefault}{\mddefault}{\updefault}$2\eta$}}}
\put(2668,-409){\makebox(0,0)[lb]{\smash{\SetFigFont{8}{9.6}{\rmdefault}{\mddefault}{\updefault}$2\eta-\alpha$}}}
\put(2108,387){\makebox(0,0)[lb]{\smash{\SetFigFont{8}{9.6}{\rmdefault}{\mddefault}{\updefault}$2\eta$}}}
\put(624,-1805){\makebox(0,0)[lb]{\smash{\SetFigFont{8}{9.6}{\rmdefault}{\mddefault}{\updefault}$-2\eta$}}}
\put(1424,-838){\makebox(0,0)[lb]{\smash{\SetFigFont{12}{14.4}{\rmdefault}{\mddefault}{\updefault}$0$}}}
\put(1744,-79){\makebox(0,0)[lb]{\smash{\SetFigFont{8}{9.6}{\rmdefault}{\mddefault}{\updefault}$\eta$}}}
\put(1395,495){\makebox(0,0)[lb]{\smash{\SetFigFont{12}{14.4}{\rmdefault}{\mddefault}{\updefault}$\alpha$}}}
\put(2425,-155){\makebox(0,0)[lb]{\smash{\SetFigFont{12}{14.4}{\rmdefault}{\mddefault}{\updefault}$\beta$}}}
\end{picture}}
\end{minipage}
\hfill
\begin{minipage}[b]{8cm}{
In the picture one can see that the long roots $\ä$ and $\bb$ satisfy (SII).
(We illustrate the situation in detail only for $\ä$.)
Contemplate the picture for a moment one
sees that there are no short roots and no other long root
for which (SII) can be valid.

Now $\ä$ and $\bb$ are the only roots with
$\la\LL,\ä\ra>0$ and $\la\LL,\bb\ra>0$. Hence
$\ä=2e_3-e_1-e_2$ and $\bb=-2e_2+e_1+e_3$.}
\end{minipage}\\
We consider the case where $(\LL,\ä)$ satisfies (SII). There is a weight element $Q_\phi$ from $\bhg$
such that
$Q_\phi(u_\LL)=A_{2e_3-e_1-e_2}$, i.e. $\phi=
2e_3-e_1-e_2-\LL=e_2-e_1$. But this is a short root and therefore a weight.
Thus we consider $u_{-\phi}\in =V_{-\phi}$. Then
$Q_\phi(u_{-\phi})\in\mf{t}$. Since there is no planar spanning triple
it has to be zero. As above the Bianchi identity gives that
$\phi-\LL $ has to be a root. But
 $\phi-\LL=e_2-e_1-2e_3+2e_2= 3e_2-2e_3-e_1$ is no root.

 For $\bb$ one proceeds analogously.
 \eprf

Now we can draw the conclusions from the previous sections. If a Lie algebra acts irreducible of real type
the it is semi-simple and obeys the properties (SI) or (SII). The simple Lie algebras with (SI) or (SII)
we have listed above. Thus we get

\btheo
Let $\lag\subset \laso(N,\rr)$ be a irreducible weak-Berger algebra of real type.
Then it is the holonomy representation of a Riemannian manifold.

The conclusion holds in particular if $\lag$ is simple, of real type and the irreducible
component of the $\lason$-projection of an indecomposable, non-irreducible simply connected Lorentzian manifold.
\etheo

\bbem {\bf Quaternionic symmetric spaces.}\label{quatsym1}
With the result of course we have covered all simple irreducible acting Riemannian holonomy groups of real type.

If one considers a quaternionic symmetric space
$G/Sp(1)\cdot H$ with $H\subset Sp(n)$ then of course $\mf{sp}(1) \oplus \mf{h}\subset \laso(4n,\rr) $
is a real Berger algebra of real type
and thus its complexification is a complex Berger algebra of real type.
Then the restriction of this representation to $\mf{h}$ is of  quaternionic, i.e. of non-real type,  its complexification
decomposes into two irreducible components $\ccc^{2n}\oplus\overline{\ccc^{2n}}$.
For this situation in \cite{schwachhoefer2}
is proved that  $\mf{h}^\ccc_{\big|\ccc^{2n}}$ is a complex Berger algebra. This result does not collide
with our list because this representation  is not of real type and hence not orthogonal.
$ \mf{h}\subset \laso(4n,\rr)$ is not a real Berger algebra.
\ebem

\section{Weak-Berger algebras of  non-real type}

In this section we will classify weak-Berger algebras of non-real type,
and we will show that these are Berger algebras. For the classification we will use
the classification of first prolongations of irreducible complex Lie algebras.
We will show that the complexification of the space ${\cal B}_h(\lag_0)$ is isomorphic
to the first prolongation of the complexified Lie algebra.

\bigskip

In this section $\lag_0$ is a real Lie algebra and $E$ a $\lag_0$-module of non-real type, i.e.
$E^\ccc$ is not irreducible. Thus the situation is a little bit more puzzling then in the real case.

Since $\lag_0\subset \laso(E,h)$ with $h$ positive definite, $\lag_0$ is compact.
For a compact real Lie algebra with module of
non-real type
the corresponding complex representation of non-real type is not orthogonal but
unitary (See appendix \ref{real representations}, in particular proposition \ref{main1}).
But if we
switch to the complexified algebra the $(\lag^\ccc,V)$ irreducible remains, but it can no
longer
be unitary of course. We have to handle this situation.

With the same notations as in appendix \ref{real representations} the complex representations space
$W=E^\mathbb{C}$ splits into the irreducible modules $W=V \oplus \overline{V}$
under ${\frak g}_0$.  This splitting is of course ${\frak g}_0^\mathbb{C}$
invariant.

Now we define the complex Lie algebra
\begin{equation}
\label{defg} {\frak g}\ :=\   \left\{ A_{|V}\left|\ A\in {\frak g}_0^\mathbb{C}\subset
\mathfrak{so}(W=V \oplus \overline{V}, H)\right.\right\}\subset \mathfrak{gl}(V).
\end{equation}
Here $H$ denotes again $h^\ccc$.
Since the symmetric bilinear form we start with is positive definite the appendix \ref{real representations}
gives two important results (see proposition \ref{main1}):
\bnum
\item \label{unitaer} Since $\lag_0$ is compact
there is a positive definite hermitian form $\theta^h$ on $V$ which is the
the restriction of the sesqui-linear extension of $h$ on $V$, such that
$({\frak g}_0)_{|V}\subset{\frak u}(V, \theta^h)$.
\item\label{nichtorthogonal} $\lag$ is not orthogonal, in particular $H_{\big|V\times V}=0$. This is the case
since modules of non-real type are symplectic if they are self-dual. Thus they can not be orthogonal.
\enum
In ${\frak g}_0^\mathbb{C}$ as well as in ${\frak g}$ we have a conjugation
$\overline{{ }^{\left.\ \right.}}$ with respect to $ {\frak g}_0$ and $({\frak
g}_0)_{|V}$ respectively.

Since an $A\in {\frak g}_0$ acts on $V\oplus \overline{V}$ by $A(v+
\overline{w})= Av + \overline{Aw}$ we have for $iA\in {\frak g}_0^\mathbb{C}$
that
\[ iA(v+ \overline{w})\ =\  i(Av + \overline{Aw})\ =\ (iAv+ \overline{-iAw}).\]
So we write the action of $A\in {\frak g}_0^\mathbb{C}$ with the help of the
conjugation in ${\frak g}$ as follows
\begin{equation}
\label{g0c}
A(v+ \overline{w})\ =\  Av + \overline{\overline{A} w}.
\end{equation}

This gives the following Lie algebra isomorphism
\[ \begin{array}{rcrcl}
\varphi&:& {\frak g}_0^\mathbb{C} &\simeq & {\frak g}\\
&&A&\mapsto& A_{|V}.
\end{array}\]
This is clearly a Lie algebra homomorphism. It is
 injective because  for $A_{|V}=B_{|V}$ holds that $A(v+ \overline{w}) =
Av + \overline{\overline{A}w}=Bv + \overline{\overline{B}w}=B(v+\overline{w}) $
for all $v,w\in V$, i.e. $A=B$.

By definition it is surjective and $\varphi^{-1}$ is given by
\begin{equation}
\label{phi-1}
 \varphi^{-1}(A)\ :\ v + \overline{w} \ \longmapsto\ Av +
\overline{\overline{A}w}\makebox[2cm][c]{for all} A\in {\frak g}.
\end{equation}
These notations are needed to show the relation to the first prolongation.

\subsection{The first prolongation of a Lie algebra of non-real type}
Now we  define the first prolongation of an arbitrary Lie algebra
${\frak g}\subset \mathfrak{gl}(V)$.

\begin{de}
The ${\frak g}$-module
\begin{align}
{\frak g}^{(1)}
&\ :=\ \{ Q \in V^* \otimes {\frak g}\ |\  Q(u)v=Q(v)u\}.
\label{hg}\\
\intertext{is called {\bf first prolongation} of  ${\frak g}\subset \mathfrak{gl}(V)$.
Furthermore we set}
\tilde{ {\frak g}}&\  :=\  span\{Q(u)\in {\frak g}\ |\ Q\in
{\frak g}^{(1)}, u\in V\}\subset {\frak g}, \nonumber\\
\intertext{and if in ${\frak g}$ a conjugation $\overline{{ }^{\left.\
\right.}}$ is given:}
{\frak g}^{[1,1]}&\ :=\  \{ R\in \overline{V}^* \otimes {\frak g}^{(1)}\ |\
\overline{R(\overline{u},v)}=- R(\overline{v},u)\},\nonumber\\
\tilde{\tilde{{\frak g}}}&\ :=\ span \{ R(\overline{u},v)\ |\ R\in {\frak
g}^{[1,1]}, \overline{u}\in \overline{V}, v\in V\}\subset {\frak g}\nonumber.
\end{align}
\end{de}

We will now describe the spaces ${\cal B}_{H}({\frak g}_0^\mathbb{C})$ and
${\cal K}({\frak g}_0^\mathbb{C})$ --- which are essential for the Berger and
the weak-Berger property --- with the help of the first prolongation of $
{\frak g}$.

In the setting of the above notations we can now prove
the following.

\begin{satz} \label{bhg-isom}
Let $E$ be a non-real type module of ${\frak g}_0$, orthogonal with respect to a positive
definite scalar product $h$,  and $E^\mathbb{C}=V \oplus \overline{V}$ the
corresponding $ {\frak g}_0^\mathbb{C}$ invariant decomposition, ${\frak g}$
defined as in (\ref{defg}).
Then there is an isomorphism
\[\begin{array}{rcrcl}
\phi&:& {\cal B}_{H}({\frak g}_0^\mathbb{C}) & \simeq &  {\frak g}^{(1)} \\
&&Q&\mapsto & Q_{|V\times V}.
\end{array}\]
\end{satz}

\begin{proof}
For the prove we will use the
$\lag_0$--invariant hermitian form $\theta$ on $V$ which is given by
$\theta(u,v)=h^\ccc(u,\overline{v})$,
where $\overline{{ }^{\left.\ \right.}}$ is the conjugation in $E^\ccc=V\+\overline{V}$
with respect to $E$.
The linearity of $\phi$ mapping is clear. we have to show the following:

1.) The definition of $\phi$  is correct, i.e. for $Q\in {\cal B}_{H}({\frak
g}_0^\mathbb{C})$ it is $Q_{|V\times V}\in {\frak g}^{(1)}$.
We have for every $u,v,w \in V$ and $H= h^\mathbb{C}$
that
\begin{eqnarray*}
\theta(Q(u)v,w)&=& h^\mathbb{C}(Q(u)v,\overline{w})
\\
&=&
-h^\mathbb{C}(Q(v)\overline{w},u)\;\; -
\underbrace{h^\mathbb{C}(Q(\overline{w})u,v)}_{\makebox[3cm][c]{
\begin{minipage}{4.5cm}{\scriptsize \centerline{$=0$} since  $h^\mathbb{C}_{V\times V}=0$
(proposition \ref{main1})}
\end{minipage}}}
\\
&\stackrel{\mbox{{\scriptsize $h^\mathbb{C}$ invariant}}}{=}&
h^\mathbb{C}(Q(v)u, \overline{w})
\\
&=&
\theta(Q(v)u,w),
\end{eqnarray*}
i.e. $Q(u)v=Q(v)u$ which means that $Q_{|V\times V} \in
{\frak g}^{(1)}$.

2.) The homomorphism $\phi$ is injective.  Let $Q_1$ and $Q_2$ be in $ {\cal
B}_{H}({\frak g}_0^\mathbb{C})$ with $(Q_1)_{|V\times V} =
(Q_2)_{|V\times V}$.
Then it is
\begin{enumerate}
\item[a)] $(Q_1)_{|\overline{V}\times \overline{V}} = (Q_2)_{|\overline{V}\times
\overline{V}}$, since $Q_1(\overline{u})\overline{v} = \overline{ Q_1(u)v}= \overline{ Q_2(u)v}
= Q_2(\overline{u})\overline{v},$
\item[b)] $(Q_1)_{|\overline{V}\times V} = (Q_2)_{|\overline{V}\times V}$, since
\[\begin{array}{rcccl}
\theta(Q_1(\overline{u})v,w)& =&
h^\mathbb{C}(Q_1(\overline{u})v,\overline{w})
&=&
-h^\mathbb{C}(v,Q_1(\overline{u})\overline{w})\ =
\\
=\ h^\mathbb{C}(v,Q_2(\overline{u})\overline{w})
&=&
h^\mathbb{C}(Q_2(\overline{u})v,\overline{w})
&=&
\theta(Q_2(\overline{u})v,w).
\end{array}\]
\item[c)] $(Q_1)_{|V\times \overline{V}} = (Q_2)_{|V\times \overline{V}}$
because of b) with the same argument as in a).
\end{enumerate}

3.) The homomorphism $\phi$ is surjective.  For $Q\in {\frak g}^{(1)}$ we define
$\phi^{-1}$ using $\varphi$:
\begin{eqnarray*}
( \phi ^{-1}Q)(u):= \varphi^{-1}(Q(u))&\mbox{and}&
 (\phi ^{-1}Q)(\overline{u}):= \varphi^{-1}(\overline{Q(u)})\in \mathfrak{gl}
(E^\mathbb{C}),\\
\mbox{i.e. }
( \phi ^{-1}Q)(u,v)=Q(u)v&,&  (\phi ^{-1}Q)(\overline{u})=\overline{Q(u)}v\ ,\\
 (\phi ^{-1}Q)(u, \overline{v})= \overline{\overline{Q(u)}v}&,&
( \phi ^{-1}Q)(\overline{u}, \overline{v})= \overline{Q(u)v}.
\end{eqnarray*}
It is $(\phi ^{-1}Q)(\overline{u}, \overline{v})= \overline{(\phi ^{-
1}Q)(u,v)}$.

Then obviously
$\phi \circ \phi^{-1} = id$, since $\phi \left( \phi^{-1}(Q)\right)=\phi^{-1}(Q)_{|V\times V}=Q$.

Because of the symmetry of $Q$ we have also that $(\phi ^{-1}Q)\in {\cal B}_{H}({\frak g}_0^\mathbb{C})$:
\begin{align*}
\intertext{$\bullet$ For $u,v\in V, \overline{w}\in \overline{V}$:}
H( (\phi^{-1}Q)(u)v,\overline{w})+
H( (\phi^{-1}Q)(v)\overline{w},u)+
\overbrace{H( (\phi^{-1}Q)(\overline{w})u,v)}^{\makebox[2cm][l]{{\scriptsize $=0$
because $H=0$ on $V\times V$}}}
&=
\\
H( (\varphi^{-1}(Q(u))v,\overline{w})+
\underbrace{H( (\varphi^{-1}(Q(v))\overline{w},u)}_{= -H(
(\overline{w},\varphi^{-1}(Q(v))u)}&=
 H( (Q(u)v -  Q(v)u ,\overline{w})\\
&=0.
 \\
\intertext{$\bullet$ For $u\in V,\overline{v}, \overline{w}\in
\overline{V}$:}
\overbrace{H( (\phi^{-
1}Q)(u)\overline{v},\overline{w})}^{=0}+
H( (\phi^{-1}Q)(\overline{v})\overline{w},u)+
H( (\phi^{-1}Q)(\overline{w})u,\overline{v})& =
\\
 H( (\varphi^{-1}(\overline{Q(v)})\overline{w},u)+
\underbrace{H( \varphi^{-1}(\overline{Q(w)})u,\overline{w})}_{= -H(
(u,\varphi^{-1}(\overline{Q(w)})\overline{v})}
&=
 H(\overline{ (Q(v)w -  Q(w)v} ,u)\\&=0.
\end{align*}
Terms with entries only from $V$ or only from $\overline{V}$ are zero.
\end{proof}

Furthermore we show for the space $ {\cal K}({\frak g})$ an analogous result.

\begin{satz} \label{kg-isom}
Let $E$ be an orthogonal  non-real type module of ${\frak g}_0$ and $E^\mathbb{C}=V \oplus
\overline{V}$ the
corresponding $ {\frak g}_0^\mathbb{C}$ invariant decomposition, ${\frak g}$
defined as in (\ref{defg}). Suppose that $\theta:=\theta^h$ is non-degenerate.
Then there is an isomorphism
\[\begin{array}{rcrcl}
\psi&:& {\cal K}({\frak g}_0^\mathbb{C}) & \simeq & {\frak g}^{[1,1]} \\
&&R&\mapsto & R_{|\overline{V}\times V\times V}.
\end{array}\]
\end{satz}

\begin{proof}
The proof is completely analogous to the previous one.

1.) The definition is correct.
We have for $u,v,w\in V$ and $R\in {\cal K}( {\frak g}_0^\mathbb{C}  )$ that
\[
\underbrace{R(u,v)\overline{w}}_{\in \overline{V}}\ =\
\underbrace{R(\overline{w},v)u}_{\in V}-
\underbrace{R(\overline{w},u)v}_{\in V}\ =\ 0.
\]
but this means that $R(\overline{u},.)_{|V\times V}\in {\frak g}^{(1)}$.

Further $R(u,v) \overline{w}=0$ implies $R(u,v)w=0$ because
\[
\theta(R(u,v)w,z)\ =\  h^\mathbb{C}(R(u,v)w, \overline{z})\ =\
-  h^\mathbb{C}(w,R(u,v) \overline{z})\ =\ 0.\]
This  implies  $R( \overline{u}, \overline{v}) \overline{w}=R( \overline{u},
\overline{v}) w=0$ too.

For a $ R \in {\cal K}( {\frak g}_0^\mathbb{C})$ we have due to the skew
symmetry
\[\overline{R(\overline{u},v)}\stackrel{\mbox{ easy calculation }}{=}
 R(u, \overline{v})\stackrel{\mbox{skew-symm.}}{=}-R( \overline{v},u),\]
i.e. the restriction of $R$ on $\overline{V}\times V  \times V$ is in $ {\frak
g}^{[1,1]}$.

2.) The homomorphism $\psi$ is injective.

Let $R_1$ and $R_2$ be in $ {\cal K}({\frak g}_0^\mathbb{C})$ with
$(R_1)_{\overline{V}\times V\times V}=(R_2)_{\overline{V}\times V\times V}$.
Then again via $\theta$ the remaining non zero terms $R_i( \overline{u},
v)\overline{w}$ are determined by $R_i( \overline{u},v)w$ which are equal for
$i=1,2$ and by the skew symmetry of $R$.

3.) The homomorphism $\psi$ is surjective.

We set
\begin{eqnarray*}
 (\psi ^{-1}R) (\overline{u},v))\ :=\  \varphi^{-1}(R( \overline{u},v)&\mbox{ ,
}&
(\psi ^{-1}R) (u,\overline{v})\ :=\    \varphi^{-1}(\overline{R(
\overline{u},v)})\mbox{ and}\\
 (\psi ^{-1}R)(u,v)&:=& (\psi ^{-1}R)(\overline{u}, \overline{v})\ :=\ 0
\end{eqnarray*}
So we have the skew symmetry, i.e. $\psi^{-1}R\in \wedge^2
E^\mathbb{C}\otimes {\frak g}_0^\mathbb{C}$,  because
\[(\psi ^{-1}R) (u,\overline{v})=   \varphi^{-1}(\overline{R( \overline{u},v)})
=-\varphi^{-1}(R( \overline{v},u))
=- (\psi^{-1}R)(\overline{v},u).\]
The Bianchi identity is also satisfied:
\begin{align*}
\intertext{$\bullet$ For $u\in \overline{V},v,w\in V$:}
(\psi^{-1}R)(\overline{u},v)w+
\overbrace{(\psi^{-1}R)(v,w)\overline{u}}^{=0}
+(\psi^{-1}R)(w,\overline{u})v&=\\
\varphi^{-1}\left(R(\overline{u},v)\right)w+
\varphi^{-1}\left(\overline{R(\overline{w},u)}\right)v
&=\\
\varphi^{-1}\left(R(\overline{u},v)\right)w-
\varphi^{-1}\left(R(\overline{u},w)\right)v
&=\\
R(\overline{u},v)w-R(\overline{u},w)v
&=0
\\
\intertext{$\bullet$ For $\overline{u}, \overline{v}\in \overline{V}, w\in
V$:}
\overbrace{(\psi^{-1}R)(\overline{u},\overline{v})w}^{=0}+
(\psi^{-1}R)(\overline{v},w)\overline{u}
+(\psi^{-1}R)(w,\overline{u})\overline{v}
&=\\
\varphi^{-1}\left(R(\overline{v},w)\right)\overline{u}+
\varphi^{-1}\left(\overline{R(\overline{w},u)}\right)\overline{v}
&=\\
-\varphi^{-1}\left(\overline{R(\overline{w},v)}\right)\overline{u}+
\varphi^{-1}\left(\overline{R(\overline{w},u)}\right)\overline{v}
&=\\-\overline{R(\overline{w},v)u} +
\overline{R(\overline{w},u)v}
&=0
\end{align*}
Terms with entries only from $V$ or only from $\overline{V}$ are zero.
\end{proof}

In contrary to the previous proof, in this proof we only supposed the fact that $\theta^h$ is non-degenerate
and not that $h^\mathbb{C}_{|V\times V}=0$. If we assume $h$ to be positive definite, then both facts
are satisfied.

\subsection{Consequences for Berger and weak-Berger algebras}

Both propositions give three important corollaries.

\begin{folg}\label{wichtig}
Let $ {\frak h}_0\subset {\frak g}_0\subset {\mathfrak
so}(E^\mathbb{C}, H)$ be  subalgebras of non-real type, ${\frak h}$ and ${\frak
g}$ defined as above.  If
\[  {\frak h}^{(1)}=   {\frak g}^{(1)},\]
then $( {\frak h}_0^\mathbb{C})_{H}=({\frak
g}_0^\mathbb{C})_{H}$.  I.e. if in ${\frak g}$ exists a proper subalgebra
which has the same first prolongation and  a compact real form in ${\frak g}_0$
of non-real type, then ${\frak g}_0^{\mathbb{C}}$ and therefore ${\frak g}_0$ can not
be weak-Berger algebras.
\end{folg}

\begin{proof}
Because of $Q\in {\cal B}_{H}({\frak h}_0^\mathbb{C})\simeq {\frak
h}^{(1)}={\frak g}^{(1)}\simeq {\cal
B}_{H}({\frak g}_0^\mathbb{C})$ we have
$Q(u)\in ({\frak g}_0^\mathbb{C})_{H}$ if and only if $Q(u)\in ({\frak
h}_0^\mathbb{C})_{H}$. \end{proof}

\begin{folg}\label{wichtiger}
Let ${\frak g}_0\subset  \mathfrak{so}(E^\mathbb{C}, H) $ be a Lie algebra
of non-real type, and ${\frak g}$ defined as above. Then
\begin{enumerate}
\item $({\frak g}_0^\mathbb{C})_H = {\frak g}_0^\mathbb{C}$ (i.e. ${\frak
g}_0^\mathbb{C}$ is a weak-Berger-algebra) if and only if $ {\frak g}= \tilde{
{\frak g}}$.
\item $\underline{{\frak g}_0^\mathbb{C}}= {\frak g}_0^\mathbb{C}$ (i.e. ${\frak
g}_0^\mathbb{C}$ is a Berger-algebra) if and only if ${\frak g}=\tilde{\tilde{
{\frak g}}}$.
\end{enumerate}
\end{folg}
\begin{proof}
1.) First we show the sufficiency: Let $A\in {\frak g}_0^\mathbb{C}$ be arbitrary.  The
assumption ${\frak g} = \tilde{ {\frak g}}$ gives w.l.o.g. that $\varphi(A)=
Q(u)$ with $Q\in {\frak g}^{(1)}$ and $u\in V$.  But then we have
\[(\phi^{-1}Q)(u)\stackrel{\mbox{ per def. }}
{=} \varphi^{-1}(Q(u))\ =\ \varphi^{-1}( \varphi(A))\ =\ A,\]
with $(\phi^{-1}Q)\in {\cal B}_H( {\frak g}_0^\mathbb{C})$, i.e. $A\in
({\frak g}_0^\mathbb{C})_H$.

Now we show the necessity: If $A\in {\frak g}$, then the assumption $ {\frak
g}_0^\mathbb{C}= ( {\frak g}_0^\mathbb{C})_H$ gives w.l.o.g. that $
\varphi^{-1}(A)=\hat{Q}(u+ \overline{v})$ with $\hat{Q}\in  {\cal
B}_H({\frak g}_0^\mathbb{C})$, $u\in V $ and $\overline{v}\in
\overline{V}$. But by the isomorphism of the proposition \ref{bhg-isom} there is
a $Q\in {\frak g}^{(1)}$ such that
\[
 \varphi^{-1}(A)=\hat{Q}(u+ \overline{v})=( \phi ^{-1} Q)(u+\overline{v})
= \varphi^{-1}(Q(u)) + \varphi^{-1}(\overline{Q(v)}).\]
But this means that
\[A\ = \ \underbrace{Q(u)}_{\in \tilde{{\frak g}}} +
\underbrace{\overline{Q(v)}}_{\in \tilde{{\frak g}}}\ \in \ \tilde{ {\frak
g}},\]
i.e. $ {\frak g}\subset \tilde{ {\frak g}}$.

2.) Both directions are proved completely analogous to 1.)

Suppose that ${\frak g}=\tilde{\tilde{g}}$. Then for $A\in {\frak g}_0^\mathbb{C}$ one has that
$\varphi(A)= R( \overline{u},v)$ and
\[(\psi^{-1}R)( \overline{u},v)\ =\ \varphi^{-1}(R( \overline{u},v)\ =\ A.\]
On the other hand we have for $A\in {\frak g}$ that $\varphi^{-1}(A)=
\hat{R}(z+\overline{u}, v+\overline{w})$.
This gives
\[\begin{array}{rcccl}
 \varphi^{-1}(A)&=&\hat{R}(z,\overline{w})+ \hat{R}(\overline{u},v)&=&
( \psi ^{-1} R)(z,\overline{w})+=( \psi ^{-1} R)(\overline{u},v)\\
&=&
\varphi^{-1}(\overline{R(\overline{z},w)}) + \varphi^{-1}(R(\overline{u},v))&&.
\end{array}\]
and therefore $A\in \tilde{\tilde{{\frak g}}}$. \end{proof}

As a result of the previous and this section we have to investigate complex
irreducible representations of complex Lie algebras with non-vanishing first
prolongation. Fortunately these are classified by Cartan \cite{cartan09},
Kobayashi and Nagano \cite{ko-na65} in a rather short list. In the next section
we will present this list and check for the entries with the help of the previous
corollaries whether they are Berger or weak-Berger algebras.

\subsection{Lie algebras  with non-trivial first prolongation and the result}

There are only a few complex Lie algebras ${\frak g}$ contained irreducibly in
$\mathfrak{gl}(V)$which have non
vanishing first prolongation. The classification is due to \cite{cartan09} and
\cite{ko-na65}. We will cite them following \cite{schwachhoefer1} in two tables.

\begin{center}
\begin{tabular}{r|c|c|cl|c}
\multicolumn{6}{l}{{\bf Table 1} Complex Lie-groups and algebras with  ${\frak
g}^{(1)}\not=0 $ and ${\frak g}^{(1)}\not=V^*$:}\\
\multicolumn{6}{l}{ }\\
&$G$&${\frak g}$& \multicolumn{2}{|c|}{$V$} &${\frak g}^{(1)} $\\
\cline{1-6}
 1.& $Sl(n, \mathbb{C})$&$\mathfrak{sl}(n, \mathbb{C})$&$
\mathbb{C}^n,$&$n\ge2$&$(V \otimes \odot^2  V^*)_0$
\\
2.& $Gl(n, \mathbb{C})$&$\mathfrak{gl}(n, \mathbb{C})$&$ \mathbb{C}^n$,&$n\ge
1$&$V \otimes \odot^2 V^*$
\\
3.&$Sp(n, \mathbb{C})$& $\mathfrak{sp}(n, \mathbb{C})$&$ \mathbb{C}^{2n}$,&$n\ge
2$&$ \odot^3 V^*$
\\
4.& $\mathbb{C}^* \times Sp(n, \mathbb{C})$&$\mathbb{C} \oplus\mathfrak{sp}(n,
\mathbb{C})$,&$
\mathbb{C}^{2n}$,&$n\ge 2$&$\odot^3 V^*$
\end{tabular}
\end{center}
\vspace{.3cm}
\begin{center}\label{table2}
\begin{tabular}{r|c|c|cl}
\multicolumn{5}{l}{{\bf Table 2} Complex Lie-groups and algebras with first
prolongation ${\frak
g}^{(1)}=V^*$:}\\
\multicolumn{5}{l}{ }\\
&$G$&${\frak g}$&\multicolumn{2}{|c}{$ V$} \\ \cline{1-5}
1.&$CO(n, \mathbb{C})$&$ \mathfrak{co}(n, \mathbb{C})$&$ \mathbb{C}^{n}$,&$n\ge
3$
\\
 2.&$Gl(n, \mathbb{C})$&$ \mathfrak{gl}(n, \mathbb{C})$&$ \odot^2
\mathbb{C}^n$,&$n\ge 2$
\\
3.&$Gl(n, \mathbb{C})$&$ \mathfrak{gl}(n, \mathbb{C})$&$
\wedge^2\mathbb{C}^n$,&$n\ge 5$
\\
4.&$Gl(n, \mathbb{C})\cdot Gl(m, \mathbb{C})$&$
\mathfrak{sl}(\mathfrak{gl}(n, \mathbb{C}) \oplus \mathfrak{gl}(m,
\mathbb{C}))$&$
\mathbb{C}^n \otimes \mathbb{C}^m$,&$m,n\ge2$
\\
5.&$\mathbb{C}^*\cdot Spin(10, \mathbb{C})$&$ \mathbb{C}
\oplus\mathfrak{spin}(10, \mathbb{C})$&$\Delta^+_{10}\simeq\mathbb{C}^{16}$&
\\
6.&$\mathbb{C}^* \cdot E_6$&$ \mathbb{C} \oplus {\frak e}_6$&$
\mathbb{C}^{27}$&
\end{tabular}
\end{center}

We have to make two remarks about the second table:

The fourth Lie algebra is defined as
\begin{eqnarray*}
\mathfrak{sl}(\mathfrak{gl}(n, \mathbb{C}) \oplus \mathfrak{gl}(m,
\mathbb{C}))&=&\{(X,Y)\in \mathfrak{gl}(n, \mathbb{C}) \oplus \mathfrak{gl}(m,
\mathbb{C})|tr\ X + tr\  Y=0\}
\\
&=&
\left(\mathfrak{gl}(n, \mathbb{C}) \oplus \mathfrak{gl}(m,
\mathbb{C})\right)\cap \mathfrak{sl}(n+m, \mathbb{C}).
\end{eqnarray*}
The identification with the Lie algebra of the group is given as follows
\begin{eqnarray*}
\mathfrak{sl}(\mathfrak{gl}(n, \mathbb{C}) \oplus \mathfrak{gl}(m,
\mathbb{C}))&\simeq&LA(Gl(n,\mathbb{C})\cdot GL(m,\mathbb{C}))\subset
\mathfrak{gl}(n\cdot m, \mathbb{C})\\
(A,B)&\longmapsto& (x \otimes u\mapsto Ax \otimes u - x \otimes Bu).
\end{eqnarray*}
In entry 5. $\Delta^+_{10}$ denotes the irreducible $Spin(10, \mathbb{C})$
spinor module.  The representation in 6. is one of the two $27$-dimensional,
irreducible ${\frak e}_6$ representations, which are conjugate to each other as representations
of the compact real form of ${\frak e}_6$.

\paragraph{The algebras of table 1}

The first three entries of table 1 are all complexifications of Riemannian
holonomy algebras $ \mathfrak{su} (n)$, $ \mathfrak{u} (n)$  acting on
$\mathbb{R}^{2n}$ and
$\mathfrak{sp}(n)$ acting on $\mathbb{R}^{4n}$ and therefore Berger algebras.

The fourth has the compact real form $i \mathbb{R} \oplus \mathfrak{sp}(n)\simeq
\mathfrak{so}(2) \oplus \mathfrak{sp}(n)$
acting irreducible on $ \mathbb{R}^{4n}$ where $i\ id$ corresponds to the
element $J\in \mathfrak{u} (2n)$.  Since the representation of $\mathfrak{sp}
(n) $ on $ \mathbb{R}^{4n}$ is of non-real type we are in the situation of corollary
\ref{wichtig}, because $( \mathbb{C} Id \oplus \mathfrak{sp}(n,
\mathbb{C}))^{(1)}=
\mathfrak{sp}(n, \mathbb{C})^{(1)}$. Hence $\mathbb{C} \oplus {\mathfrak
sp}(2n, \mathbb{C})$ is not a weak-Berger algebra.

\paragraph{The algebras of table 2}

If one looks at the unique (up to inner automorphisms) compact real form and the reellification of the Lie
algebras and representations in table 2
one sees that they correspond to the holonomy representation of
Riemannian symmetric spaces which are K\"ahlerian. This gives the following
proposition.

\begin{samepage}
\begin{satz}\label{unbewiesen}
The compact real forms of the algebras in  table 2 and the reellification of the
representations are equivalent to the holonomy representations of the following
Riemannian, K\"ahlerian symmetric spaces (see \cite{helgason78}):
\begin{center}
\begin{tabular}{r|l|c|c|l}
&Type&non-compact&compact&dim.\\
\cline{1-5}
 1.& $BD \ I$&$SO_0(2,n)\big/SO(2)\times SO(n)$&$
SO(2+n)\big/SO(2)\times SO(n)$&$2n$\\
2.&$C\ I$&$Sp(n, \mathbb{R})\big/U(n)$&$Sp(n)\big/U(n)$&$n(n+1)$\\
3.&$D\ III$&$SO^*(2n)\big/U(n)$&$SO(2n)\big/U(n)$&$n(n-1)$\\
4.&$A\ III$&$SU(n,m)\big/U(n)\cdot U(m)$&$SU(n+m)\big/U(n)\cdot U(m)$&$2nm$\\
5.&$E\ III$&$\left( {\frak e}_{6(-14)}, \mathfrak{so}(2) \oplus
\mathfrak{so}(10) \right)$&$\left( {\frak e}_{6(-78)}, \mathfrak{so}(2) \oplus
\mathfrak{so}(10) \right)$&$32$\\
6.&$E\ VII$&$\left( {\frak e}_{7(-25)}, \mathfrak{so}(2) \oplus \mathfrak{e}_6
\right)$&$\left( {\frak e}_{7(-133)}, \mathfrak{so}(2) \oplus \mathfrak{e}_6
\right)$&$54$\\
\multicolumn{5}{l}{ }\\
\multicolumn{5}{l}{{\bf Table 3} Riemannian, K\"ahlerian symmetric spaces corresponding to
table 2 }
\end{tabular}
\end{center}
\end{satz}
\end{samepage}

So we obtain that all algebras corresponding to table 2 are Berger algebras and
therefore also weak-Berger algebras.

\btheo
Let $\lag$ be a Lie algebra and $E$ an irreducible $\lag$--module of non-real type.
If $\lag\subset \laso(E,h)$ is a weak-Berger algebra then it is a Berger algebra.
\etheo

\paragraph{Consequences for Lorentzian holonomy}
All in all we have shown, that every real Lie algebra ${\frak g}_0$ of non-real type,
i.e. contained in $\mathfrak{u} (n)$,  that can be weak-Berger  is a Berger
algebra. Further each of these Lie algebras is the holonomy algebra of a Riemannian
manifold, the remaining entries of table 1 of non-symmetric ones, and the entries of
table 2 of symmetric ones.

Before we  we apply this to the irreducible components of the
$\mathfrak{so}(n)$-projection of the holonomy algebra of an indecomposable
Lorentzian manifold with light like invariant subspace, we prove a lemma to get
the result in full generality.

\begin{lem}
Let ${\frak g}\subset {\frak u}(n)\subset \mathfrak{so}(2n)$ be a Lie algebra
with the decomposition property of theorem \ref{theoI}, ie. there
exists decompositions of $\mathbb{R}^{2n}$ into orthogonal subspaces and of
${\frak g}$ into ideals
\[
\mathbb{R}^{2n} = E_0 \oplus E_1 \oplus \ldots \oplus E_r\ \mbox{ and }\
{\frak g} = {\frak g}_1 \oplus \ldots \oplus {\frak g}_r\]
where ${\frak g}$ acts trivial on $E_0$, ${\frak g}_i$ acts irreducible on $E_i$
and ${\frak g}_i (E_j)=\{0\}$ for $i\not=j$.  Then ${\frak g}\subset {\frak
u}(n)$ implies  $dim\  E_i= 2 k_i$ and
${\frak g}_i \subset {\frak u}(k_i)$ for $i=1, \ldots ,r$.
\end{lem}

\begin{proof}
 Let $\rr^{2n}= \mathbb{C}^n$ and $ \theta$ be the positive definite hermitian form
on $\mathbb{C}^n$. Let $E_i $ be an invariant subspace on which ${\frak g}$ acts
irreducible.  If $E_i= V^i_\mathbb{R}$ for a complex vector space $V^i$, then we
can restrict $\theta$ to $V^i$. Because $\theta$ is positive definite it is
non-degenerate on $V^i$ --- since $\theta(v,v)>0$ for $v \not=0$ ---  we get
that $ {\frak g}_i\subset {\frak u}(V^i, \theta)$, i.e. ${\frak g}\subset {\frak
u}(k_i)$.

Hence we have to consider a subspace $E_i$ which is not the reellification of a
complex vector space. Let $J$ be the complex structure on $\mathbb{R}^{2n}$.
We consider the real vector space $JE_i$, which is invariant under $ {\frak g}$,
since $J$ commutes with ${\frak g}$. Then the space $JE_i \cap E_i$ is contained
in $E_i$ as well as in $JE_i$ and invariant under $ {\frak g}$. Because ${\frak
g}$ acts irreducible on $E_i$ we get two cases. The first is $E_i \cap JE_i =
E_i = JE_i$, but this was excluded since $E_i$ was not a reellification.  The second
is $E_i\cap JE_i=\{0\}$. So we have  two invariant irreducible subspaces on
which $ {\frak g}$ acts simultaneously, i.e. $A(x, Jy)=(Ax, AJy)$, but this is
not possible because of the Borel-Lichnerowicz decomposition property from
theorem \ref{theoI}.
\end{proof}

\begin{theo}\label{theoun}
Let $(M,h)$ be an indecomposable $n+2$-dimensional  Lorentzian manifold with
light like holonomy-invariant subspace.  Set ${\frak g}:= pr_{\mathfrak{so}(n)}
\mathfrak{hol}_p(M,h)$  and suppose ${\frak g}\subset \mathfrak{u} (n)$. Then
${\frak g}$ is the holonomy algebra of a Riemannian manifold.
\end{theo}

\begin{proof} ${\frak g}\subset \mathfrak{u} (n)$ is a weak-Berger algebra.
Then all the ${\frak g}_i$ of the decomposition of theorem \ref{theoI} are
unitary because of the lemma and weak-Berger because of corollary
\ref{zerlegung}. Hence they are weak-Berger of non-real type. Then ${\frak g}_i$
corresponds to a compact real form of the entries of table 1 or 2. But these are
all Riemannian holonomy algebras, and therefore $\mathfrak{g}$ is a Riemannian holonomy algebra. \end{proof}

\bbem {\bf Quaternionic symmetric spaces.}\label{quatsym2}
Again we have to make a remark about quaternionic symmetric spaces (see remark \ref{quatsym1}).
If $G/Sp(1)\cdot H$ with $H\subset Sp(n)$ is a quaternionic symmetric space then the
corresponding complex irreducible representation of $H$ is of quaternionic, i.e. of non real type, and it is
Berger \cite{schwachhoefer2}. But the real representation of $H$, i.e. the reellification of the complex one,
is not. Thats why it does not occur in the above list. The place of
$Sp(1)\cdot H$ would be in a list of real
semisimple, but non-simple, weak-Berger algebras of real type.
\ebem

\begin{appendix}
\section{Representations of real Lie algebras}
\label{real representations}

In this appendix we will collect and illustrate some standard facts about representations of real Lie algebras.

Because of the theorem \ref{theoI} and proposition \ref{zerlegung} we are interested
in irreducible real representations of real Lie algebras which are orthogonal.

First we will recall some facts about irreducible complex representations of real Lie
algebras, in particular orthogonal or unitary ones.

Then we will use the results of E. Cartan (\cite{cartan1914}, see also
\cite{goto78}, pp.363 and \cite{iwahori59}), in order to reduce the study of
real representations to that of complex ones.

Throughout the whole section ${\frak g}$ is a real Lie algebra.

\subsection{Preliminaries}

First of all we recall the Schur-lemma.

\begin{satz}[Schur-lemma]
Let $\kappa_1$, $\kappa_2$ be irreducible representations of ${\frak g}$ on
$\mathbb{K}$-vector spaces  $V_1$ and $V_2$. Let $f\in Hom_{\frak g}(V_1, V_2)$
be an invariant homomorphism, i.e.
\[ f \circ \kappa_1(A) \ = \ \kappa_2(A) \circ f\makebox[4cm][r]{for all $A\in
{\frak g}$.}\]
Then holds
\begin{enumerate}
\item $f$ is zero or an isomorphism, i.e. $V_1 \not\simeq V_2$ implies
$Hom_{\frak g}(V_1, V_2)=0$.
\item If $V_1=V_2=:V$ and if $f$ has an eigenvalue $\lambda\in \mathbb{K}$, then
 $f= \lambda\ id_V$. I.e. if $\mathbb{K}= \mathbb{C}$ and $V_1=V_2$ we have
always $f=\lambda\ id$ with $\lambda \in \mathbb{C}$.
\end{enumerate}
\end{satz}

For invariant bilinear forms, i.e. forms $\beta$ which satisfy
\begin{equation}
\label{invarianz}
\beta( \kappa(A)u,v) + \beta(u, \kappa(A)v)=0
\makebox[3cm][r]{for all $A\in {\frak g}$}
\end{equation}
this gives the following consequence.

\begin{folg}
Let $\kappa$ be an irreducible representation of $ {\frak g} $ on a
$\mathbb{K}$--vector space $V$ and $ \beta$ be the invariant bilinear form. Then $ \beta$ is
zero or non-degenerate.

If $\mathbb{K}= \mathbb{C}$, then the space of invariant bilinear forms is zero
or one-dimensional. It is one-dimensional if and only if  $V\simeq_\kappa V^*$.
Then it is generated by a symmetric or an anti-symmetric bilinear form.
\end{folg}

This consequence is obvious by applying the Schur-lemma to the endomorphism
of $V$, which is induced by two invariant bilinear forms.

For complex representations and invariant sesqui-linear forms, i.e. forms
$\theta$ with
\[\theta (\lambda u,v)= \lambda \theta (u,v)\makebox[2cm][c]{ and }
\theta (u,\lambda v)= \overline{\lambda} \theta (u,v),\]
one has an analogous result.
\begin{folg}
Let $\kappa $ be an irreducible representation of $ {\frak g}$ on a $
\mathbb{C}$-vector space $V$. Every invariant sesqui-linear form is zero or
non-degenerate, and the space of invariant sesqui-linear-forms is zero or one-dimensional.
It is one dimensional if and only if
$\overline{V}\simeq_\kappa V^*$. In this case it is generated by a hermitian or
an anti-hermitian form, and the spaces of invariant hermitian and invariant
anti-hermitian forms are one-dimensional real subspaces, identified by the multiplication with $i$.
\end{folg}

In these corollaries we refer to the dual and the conjugate representations,
which are defined as follows:
\begin{eqnarray*}
(\kappa^*(A) \alpha)v&=& - \alpha (\kappa(A)v)\\
\overline{\kappa}(A)\overline{v}&=&\overline{\kappa(A)v}.
\end{eqnarray*}

\begin{de}
Let $\kappa$ be an arbitrary representation of a Lie algebra ${\frak g}$ on a
$\mathbb{K}$-vector space $V$.
\begin{enumerate}
\item Then $\kappa$ is called {\bf self-dual} if there is an invariant
isomorphism between $V$ and $V^*$. This is equivalent to the existence of an
invariant bilinear form $\beta$.
\item If $ \mathbb{K}=\mathbb{C}$, then $\kappa$ is called {\bf self-conjugate}
if there is an invariant isomorphism from $V$ to $\overline{V}$, i.e. there
exists an anti-linear bijective mapping $J:V\longrightarrow V$ which is invariant,
i.e. $J \circ \kappa(A)= \kappa(A) \circ J\ \mbox{ for all }A\in {\frak g}$.
\end{enumerate}
\end{de}

It is evident that the existence of an invariant hermitian form or a self-conjugate
representation is only possible for {\bf real Lie algebras}.

\subsection{Irreducible complex representations of real Lie algebras}

\begin{de}
Let $\kappa$ be an {\bf irreducible}  complex representation of a real Lie algebra
${\frak g}$ on $V$. $\kappa$ is called
\begin{description}
\item[of real type] if $\kappa$ is self-conjugate with $J^2=1$,
\item[of quaternionic type] if $\kappa$ is self-conjugate with $J^2=-1$ and
\item[of complex type] if $\kappa$ is not self-conjugate.\end{description}
\end{de}

From the Schur-lemma it is clear that every complex irreducible
representation is either real, complex or quaternionic: If $\kappa$ is
self-conjugate, then $J^2$ is a linear automorphism of $\kappa$ so that $J^2=
\lambda \ id$.  Furthermore $\lambda$ must be real because of
\[\lambda Jv= J^2 Jv = J J^2v = J \lambda v = \overline{\lambda} Jv.\]
Dividing $J$ by $\sqrt{|\lambda|}$ one gets $\lambda=\pm1$.

Now it holds
\begin{satz} \label{realtype}
Let $\kappa$ be a complex irreducible representation of a real Lie algebra
${\frak g}$.  Then $\kappa$ is of real type if and only if the reellification
$\kappa_\mathbb{R}$ is reducible.
\end{satz}

\begin{proof}

($\Longrightarrow$) Let $\kappa $ be of real type, i.e. there is an anti-linear,
invariant automorphism of the complex representation space $V$ with $J^2= $id.
Then $J$ is $\mathbb{R}$-linear, and $V_{\mathbb{R}}$ splits into invariant,
real vector spaces
\begin{eqnarray*}
V_{\pm }&=& \{ v\in V\ |\ Jv=\pm v\  \}\\
V_{\mathbb{R}}&=& V_+\oplus V_-.
\end{eqnarray*}
So $\kappa_\mathbb{R}$ is reducible.

($\Longleftarrow$) Let $W$ be a real, $\kappa_\mathbb{R}$-invariant subspace of
$ V_\mathbb{R}$.  On $V_\mathbb{R}$ the multiplication with $i$ gives an
$\mathbb{R}$-automorphism, which defines two subspaces of $V_\mathbb{R}$: $W\cap
iW$ and $W+iW$.  Then both are complex vector spaces in an obvious way, such
that they are  complex subspaces of $V$.  Since $W$ is $\kappa_\mathbb{R}$
invariant, both are $\kappa$ invariant. Since $\kappa$ is irreducible, it
remains the case that $W\cap iW=\{0\}$ and $W\oplus $i$W=V$.  But this means
that $V=W^\mathbb{C}$ such that $W$ defines a conjugation $J$ in $V$ with the
desired properties.
\end{proof}

\paragraph{Orthogonal and unitary representations}

\begin{satz}\label{dualtype}
Let $\kappa$ be an irreducible representation of ${\frak g}$.
If $\kappa$ is of complex type, then it can not be both, unitary and
self-dual. If $\kappa$ is not of complex type, then it is unitary if and only if
it
is self-dual. In particular one has for real and quaternionic representations ($J$
denotes the automorphism):
\begin{enumerate}
\item If $\kappa$ is of real type, then it is orthogonal if and only if it is
unitary
with respect to $\theta$ for which holds $J^*\theta=\overline{\theta}$. It is
symplectic if and only if it is unitary with respect to $\theta$ satisfying
$J^*\theta=-\overline{\theta}$.
\item If $\kappa$ is of quaternionic type, then it is orthogonal if and only if
it is
unitary with respect to $\theta$ with $J^*\theta=-
\overline{\theta}$. It is symplectic if and
only if it is unitary with respect to $\theta$ satisfying
$J^*\theta=\overline{\theta}$.
\end{enumerate}
\end{satz}

\begin{proof}
Unitary is equivalent to $V^*\simeq_\kappa \overline{V}$ and therefore self-dual is the same as
$V\simeq_\kappa \overline{V}$.  This gives the proposition.  For the remaining
single points we get:

1.)  Let $\kappa$ be of real type with respect to a real structure $J$. By this
$J$ one gets from an invariant bilinear form $\beta$ an invariant sesqui-linear
form $\beta(.,J.)$ which is the complex multiple of an invariant hermitian form
$\theta$ and vice versa. Then one gets for $\beta$ symmetric/anti-symmetric:
\begin{eqnarray*}
J^*\theta(u,v)&=& \theta (Ju,Jv) \ =\ \lambda\beta(Ju, J^2v)\
\stackrel{J^2=id}{=}\  \lambda\beta(Ju, v)
\\
&=&
\pm \lambda\beta(v Ju) \ =\ \pm\theta (v,u)\ =\ \pm\overline{\theta (u,v)}.
\end{eqnarray*}
2.) analogous with $J^2=-id$. \end{proof}

\begin{folg}\label{+unitary}
If $\kappa$ is positive definite unitary, then it is
\begin{enumerate}
\item of real type if and only if it is orthogonal,
\item of complex type if and only it is not self-dual,
\item of quaternionic type if and only if it is symplectic.
\end{enumerate}
\end{folg}

\begin{proof} If $\theta $ is positive definite it can not be $J^*\theta = -
\theta$. \end{proof}

\subsection{Irreducible real representations}

For a real irreducible representation $\rho$ of a real Lie algebra $ {\frak g}$
on a real vector space $E$  two cases are possible:
$\rho^{\mathbb{C}}$ is irreducible or reducible.  We will describe these cases
due to results of E. Cartan (\cite{cartan1914}, see also
\cite{goto78}, pp.363 and \cite{iwahori59}), in order to reduce the study of
real representations to that of complex ones.

\subsubsection{Representations of real type}

\begin{satz}\label{realtype1}
Let ${\frak g} $ be a real Lie algebra and $\rho$ a representation of $ {\frak
g}$ on a real vector space $E$ such that $ \rho^{\mathbb{C}}$ is {\bf
irreducible} on $ E^\mathbb{C}$. Then the complex representation $
\rho^{\mathbb{C}}$ is of real type.

If otherwise  $\kappa$ is a complex representation of $ {\frak g} $ of {\bf real
type} on $V$, then $\kappa$ is the complexification of a real irreducible
representation of ${\frak g}$.
\end{satz}

\begin{proof}

1.) We show the existence of a $  \rho^{\mathbb{C}}$-invariant anti-linear
isomorphism $J$ with $J^2=id$.
If we denote by $J$ the conjugation in $ E^\mathbb{C}$ with respect to $E$,
then it is $J^2=1$ and we have
\[J \left(  \rho^{\mathbb{C}}(A)(u+iv) \right) =  \rho^{\mathbb{C}}(A)(u)- i
\rho^{\mathbb{C}}(A)(v)
=   \rho^{\mathbb{C}}(A)\left( J(u+iv) \right)\]
i.e.
$J$ is $ \rho^{\mathbb{C}}$-invariant.

2.) In the proof of proposition \ref{realtype} we had already shown that for
complex representations of real type holds that $V=W^\mathbb{C}$. \end{proof}
So the following definition makes sense.
\begin{de} Irreducible real representations with irreducible complexification
and irreducible complex representations with reducible reellification (i.e. of
real type) are
called representations of
{\bf real type}.
\end{de}

We have the following correspondence:
\begin{eqnarray}
\left\{
\mbox{real representation of real type}\right\}_{/\sim}&\leftrightarrow&
\left\{
\mbox{complex representations of real type}
\right\}_{/\sim}\label{type1}\\\nonumber
\rho&\mapsto& \rho^\mathbb{C}\\\nonumber
(\kappa_\mathbb{R})_{|\mbox{maximal invariant
subspace}}&\leftarrowtail& \kappa.
\end{eqnarray}
Here $\sim$ denotes the equivalence of representations.

\subsubsection{Representations of non-real type}

The situation in this case is described by the following
\begin{satz}
Let ${\frak g} $ be a real Lie algebra and
 $\rho$ be an irreducible representation of $ {\frak g}$ on a real vector space
$E$ such that $ \rho^{\mathbb{C}}$ is {\bf reducible} on $ E^\mathbb{C}$.
\begin{enumerate}
\item  If $V\subset E^\mathbb{C}$ is any invariant subspace of $
\rho^{\mathbb{C}}$. Then holds
\[E^\mathbb{C}= V \oplus \overline{V},\]
where $\overline{{ }^{\left.\ \right.}}$ is the conjugation in $E^\mathbb{C}$
with respect to $E$.  $V$ and $\overline{V}$ are irreducible and unique as
maximal invariant proper subspaces. The representations on $V$ and
$\overline{V}$ are conjugate to each other.
\item The irreducible representations of $ {\frak g}$ on $V$ and on
$\overline{V}$ are of complex or of quaternionic type, its reellifications are
equivalent to $\rho$.
\end{enumerate}
If otherwise $\kappa$ is a complex irreducible representation of complex or
quaternionic type, then $\kappa$ is the restriction on the maximal invariant proper subspace of
the complexification of $\kappa_\mathbb{R}$.

If $\kappa$ is of complex type, then exists and $\kappa_\rr$--invariant complex structure $J$ on $V_\rr$.
$\kappa$ is of quaternionic type if and only if there exists and $\kappa_\rr$--invariant quaternionic
structure $(I,J,K)$ on $V_\rr$.
\end{satz}

\begin{proof}

1.) Let $V\subset E^\mathbb{C}$ be any invariant, proper  subspace of $
\rho^{\mathbb{C}}$.  Lets denote by $\overline{{ }^{\left.\ \right.}}$ the
conjugation in
$E^\mathbb{C}$ with respect to $E$.

We consider $W:=V + \overline{V}$.  Now it is $\overline{W}=W$ which is
equivalent to $W= F^\mathbb{C}$, where $F=W\cap E$ is a real subspace of $E$.
Since $W$ is invariant under $ \rho^{\mathbb{C}}$, $F$ is invariant under
$\rho$. Now $\rho$ is irreducible and therefore $F=E$, i.e. $V +
\overline{V}= E^\mathbb{C}$.
Analogously one shows that $V \cap \overline{V}=\{0\}$, so that one gets
\[V \oplus \overline{V} =E^\mathbb{C}.\]
It remains to show that $V $ is irreducible:  This is clear since every
invariant subspace $U\subset V$ is invariant in $ E^\mathbb{C}$, but then holds
that
$U \oplus \overline{U} = E^\mathbb{C}$ which implies $U=V$. $\overline{V}$ is
irreducible too.

Hence we have two irreducible representations of ${\frak g}$, one on $V$ and one on
$ \overline{V}$, which are conjugate to each other:
\[ \rho^{\mathbb{C}}(A) \overline{v} = \overline{ \rho^{\mathbb{C}}v}.\]
So we will denote it by $\kappa$ and $\overline{\kappa}$.

2.) In order to show that $\kappa$ and $ \overline{\kappa}$ are of complex or
quaternionic type, we verify that $ \kappa_\mathbb{R}$ and $
\overline{\kappa}_\mathbb{R}$ are irreducible.

For this we show that $\kappa_\mathbb{R}$ and $\overline{\kappa}_\mathbb{R}$
are isomorphic to $\rho$.
The isomorphism between $V$ and $E$ is given by
\[\begin{array}{rcrcl}
\psi&:& V_\mathbb{R}&\longrightarrow& E\\
&&v&\longmapsto & \frac{1}{2}(v+\overline{v}).
\end{array}\]
This is obviously an isomorphism of real vector spaces. (Of course this is also
an isomorphism between $\overline{V}_\mathbb{R}$ and $E$.) It is also invariant
since
\[\psi \circ \kappa_\mathbb{R}(A)(x+iy)
=
\psi (\rho(A)x + i\rho(A)y)= \rho(A)x = \rho(A)(\psi (x+iy)
\]
for all $x+iy\in V_\mathbb{R}$.

The existence of the complex and the quaternionic structure on $V_\rr$ is clear.
\end{proof}
Again on defines:

\begin{de} Irreducible real representations with reducible complexification
and irreducible complex representations with irreducible reellification are
called representations of {\bf non-real type} (of {\bf complex} or {\bf quaternionic type} respectively).
\end{de}

Again we have the correspondence
\begin{eqnarray}\label{type2}
\left\{
\begin{array}{l}\mbox{real representations}\\
\mbox{of non-real type}
\end{array}\right\}_{/\sim}&\leftrightarrow&
\left\{
\begin{array}{l}\mbox{complex representations}\\
\mbox{of non-real type}
\end{array}\right\}_{/\approx}\\
\nonumber
\rho&\mapsto & \rho^\mathbb{C}_{|\mbox{maximal invariant
subspace}}\\\nonumber
\kappa_\mathbb{R}&\leftarrowtail&\kappa.
\end{eqnarray}
Here $\sim$ denotes the equivalence of representation and $\approx$ the
equivalence
\begin{eqnarray*}
\kappa_1\approx \kappa_2 &\Leftrightarrow& \kappa_1 \sim \kappa_2
\mbox{ or } \kappa_1 \sim \overline{\kappa_2}.
\end{eqnarray*}

On  the real space $E\simeq V_\mathbb{R}$
we have the complex structure $J$, i.e. an $\mathbb{R}$--automorphism with
$J^2=-1$ given by the multiplication with $i$: $Jv=iv$.
 $J$ commutes with $\rho$ since
\[\rho(A)(Jv)= \kappa_\mathbb{R} (A)(Jv)= \kappa (A)iv=i \kappa(A)v= J(
\kappa(A)v).\]
One describes the complex vector space $V$ as a subspace in
$E^\mathbb{C}$ as follows.
One extends the complex structure to an automorphism of $E^\mathbb{C}$
also denoted by $J$ and with the property $J^2=-1$. Then one defines
\[ V_{\pm}:= \{v\in E^\mathbb{C}|Jv=\pm i\ v\}\subset E^\mathbb{C}\]
and gets $E^\mathbb{C}= V_+ \oplus V_-$.
Furthermore it is
\begin{equation}
\label{v+form} V_{\pm}= \{x\mp iJx|x\in E \}\mbox{ and therefore
}V_{\pm}=\overline{V}_{\mp}.
\end{equation}
Then one has the following isomorphisms, invariant under the corresponding
representations:
\begin{equation}
\label{isomorphisms}
\begin{array}{rcccccl}
E&\simeq_\mathbb{R}&V&\simeq_\mathbb{C}&V_+&\simeq_\mathbb{R}& V_-
=\overline{V_+}\\
\frac{1}{2}(v+\overline{v})&\leftarrowtail
&v&\mapsto&\frac{1}{2}(v-iJv)&\mapsto&\frac{1}{2}(v+iJv).
\end{array}
\end{equation}

\subsection{Orthogonal real representations}

Let now $\rho$ be a real representation of ${\frak g}$ on $E$ which should
be orthogonal (or symplectic) with respect to a (anti-)symmetric bilinear form
$h$.

On $E^\mathbb{C}$ $h$ defines a (anti-)symmetric bilinear form by bilinear
extension, denoted by $h^\mathbb{C}$ and a (anti-)hermitian form
by conjugate linear extension in the second component, denoted by
$h^\prime$. Both
are invariant under $ \rho^{\mathbb{C}}({\frak g})$. The hermitian form has
the same signature as the symmetric form $h$. The existence of an invariant
{anti-hermitian} form is equivalent to the existence of an invariant hermitian
form.

For the conjugation in $E^\mathbb{C}$ we have the following relations
\[
\begin{array}
{rcccccl}
h^\prime(u,v)&=&h^\mathbb{C}(u,\bar{v})&=&\overline{h^\mathbb{C}
(\overline{u},v)} &=&\overline{h^\prime(v,u)}
\end{array}.\]

\subsubsection{Orthogonal  or symplectic representations of real type}

From these introductory remarks we obtain the following proposition for real type representations
which can be found in \cite{berger55}(for the orthogonal case).

\begin{satz}\label{prop1} \cite{berger55}
Let $\rho$ be a real representation of real type of a real Lie algebra ${\frak
g}$ on a real vector space $E$, orthogonal or symplectic with respect to $h$.
Let $\beta^h$ denote the complex linear and $\theta^h $ the hermitian
extension of $h$ on $V=E^\mathbb{C}$.
Then both are non-degenerate and $\rho^\mathbb{C}$ is orthogonal/symplectic with
respect to $\beta^h$ and
unitary with respect to $\theta^h$. $\theta^h$ has the same index as $h$ in case
$h$ is orthogonal.
\end{satz}

This gives a

\begin{folg}
If $\rho$ is of real type, then the space of invariant bilinear form is one-
dimensional  and generated by a symmetric or an anti-symmetric form.\end{folg}
 The {\it proof} is clear because the irreducibility of $\rho^{\mathbb{C}}$
gives that $h_1^\ccc=h_2^\ccc$, which implies $h_1=h_2$. \hfill $\Box$\\

We will now prove the other direction of proposition \ref{prop1}.

\begin{satz}\label{realtypeorthogonal}
Let ${\frak g}$ be a real Lie algebra and $\kappa$ an irreducible, complex
 representation of real type on $V$, which decomposes  $\kappa_\mathbb{R}$-invariant
into $V=E \oplus i E  $, and set $\rho = (\kappa_\mathbb{R})_{|E}$ the
corresponding irreducible real representation. If $\kappa$ is unitary (and
therefore self-dual), then $\rho$ is self-dual, i.e. orthogonal or symplectic
and we have two cases:
\begin{enumerate}
\item If $ \kappa$ is orthogonal, then  $\rho $ is orthogonal.
\item If $\kappa$ is symplectic, then $\rho $ is symplectic
\end{enumerate}
\end{satz}

\begin{proof}
Let $\kappa$ be unitary with respect to  $\theta$, which  defines two bilinear
mappings on $E$
\begin{eqnarray*}
h_1(x,y)&=& Re \left( \theta(x,y) \right)\ \mbox{ symmetric}\\
h_2(x,y)&=& Im \left( \theta(x,y) \right)\ \mbox{ anti-symmetric.}
\end{eqnarray*}
Both are $\rho$ invariant. If both are degenerate, then both are zero by the
Schur-lemma and so $\theta$ must be zero, which is a
contradiction.

1.) If in addition $\kappa$ is orthogonal, then for $\theta$ holds by
proposition \ref{dualtype} that $J^*\theta =\overline{\theta}$, where $J$ is the
conjugation of $E$ in $ E^\mathbb{C}$. But in this case $h_2$ is zero, because
$E=\{v\in V|Jv=v\}$:
\[
h_2(x,y)=Im\theta(x,y)=Im\theta(Jx,Jy)=Im\overline{\theta(x,y)}=-
Im\theta(x,y)=-h_2(x,y).\]
Hence $h_1$ must be non degenerate and therefore $\rho$ orthogonal.

2.) If $\kappa$ is symplectic one shows analogously with proposition
\ref{dualtype} that $h_1=0$ and therefore $\rho$ symplectic.
\end{proof}
Both results give the following equivalence:
\begin{eqnarray}
\label{typ1equivalenz}
\left\{
\rho\mbox{ real, real type, self-dual}
\right\}_{/\sim}&\leftrightarrow&
\left\{
\begin{array}{l}
 \kappa \mbox{ complex, real type,}\\
 \mbox{self-dual $\hat{=}$ unitary}
 \end{array}
\right\}_{/\sim}\\
\label{typ1equivalenz2}
\left\{
\begin{array}{l}
\rho\mbox{ real, real type,}\\
\mbox{orthogonal/symplectic}
\end{array}
\right\}_{/\sim}&\leftrightarrow&
\left\{
\begin{array}{l}
 \kappa \mbox{ complex, real type,}\\
\mbox{orthogonal/symplectic}\end{array}
\right\}_{/\sim}.
\end{eqnarray}

\subsubsection{Orthogonal representations of non-real type}

For non-real type representations we have the $ \rho^{\mathbb{C}}$-invariant
decomposition $E^\mathbb{C}=V \oplus\overline{V}$.

In a basis, adapted to this  decomposition $h^\mathbb{C}$ and $h^\prime$
are given as follows
\begin{eqnarray*}
h^\mathbb{C}=
\left(
\begin{array}{cc}
	A & B 		\\
B^t	 &  	\overline{A} 	\\
\end{array}
\right)
&\mbox{and}&
h^\prime=
\left(
\begin{array}{cc}
	B & A 		\\
\overline{A}	 &  	B^t	\\
\end{array}
\right)
\end{eqnarray*}
where $A=A^t$ and $B^t=\overline{B}$ are quadratic matrices with the dimension of
$V$.

Now one defines a bilinear and a sesqui-linear form on $V$ resp. on
$\overline{V}$:
\[\begin{array}{rcccl}
\beta^h(u,v)&:=& h^\mathbb{C}(u,v)&=&
h^\prime(u,\overline{v})\makebox[5cm][r]{symmetric/anti-symmetric}\\
\theta^h(u,v)&:=& h^\mathbb{C}(u,\overline{v})&=& h^\prime(u,
v)\makebox[5cm][r]{hermitian/anti-hermitian}
\end{array}\]
for $u,v\in V$ resp. $\overline{V}$. Both are invariant under $\kappa=
\rho^{\mathbb{C}}_{|V}({\frak g})$.

From the Schur-lemma it is clear that at least one of them is non-degenerate,
since
$h^\mathbb{C}$ is non-degenerate.

Using the isomorphisms of (\ref{isomorphisms}) we can give $ \theta^h$ and
$\beta^h$ explicitly:
\begin{eqnarray}
\beta^h(x-iJx, y-iJy)  &=&\label{h-omega}   \frac {1}{4}\left(  h(x,y)-
h(Jx,Jy)
-i \left(h(J  x,y) + h(x,Jy) \right) \right)\\
\theta^h (x-iJx,y-iJy)
&=&\label{h-theta}
\frac {1}{4}\left( h(x,y) + h(Jx,Jy) +i \left(h(x,Jy) - h(Jx,y) \right) \right).
\end{eqnarray}

Again we have the proposition of Berger (for the orthogonal case).

 \begin{satz}\cite{berger55}
Let $\rho $ be a real orthogonal/symplectic representation of non-real type,
i.e.
$(E,\rho)=(V_\mathbb{R},\kappa_\mathbb{R})$.
Then $\kappa$ is invariant under $ \beta^h$ and $ \theta^h$ and at
least one of them is non-degenerate, i.e. $ \kappa$ is orthogonal/symplectic or
unitary/anti-unitary with respect to $ \beta^h$ or $ \theta^h$.

Furthermore holds: If $ {\frak g}$ contains a real sub-algebra ${\frak h}\not=0$
such that ${\frak h}= {\frak p}_\mathbb{R}$ where ${\frak p}$ is a complex
Lie algebra, then $\theta^h=0$ i.e. $\beta^h$ non-degenerate.
\end{satz}

\begin{proof} We only have to prove the second assertion.

By assumption we have a complex Lie structure on ${\frak h}$, i.e. a
automorphism $J$ with $J^2=-1$ and $J \circ ad_X = ad_X \circ J$.
As above for vector spaces we have here a Lie algebra decomposition
\[ {\frak g}^\mathbb{C}\supset {\frak h}^\mathbb{C}= {\frak p}_+ \oplus
{\frak p}_-\mbox{ with }{\frak p}_{\pm}= \{v\in {\frak h}^\mathbb{C}|
Jv=\pm i v\}.\]
Then ${\frak p}\simeq_\mathbb{C} {\frak p}_+$.

Let now $\rho^{\mathbb{C}}$ be extended to ${\frak g}^\mathbb{C}$. Then
because of its linearity $h^\mathbb{C}$ is invariant under $
\rho^{\mathbb{C}}({\frak g}^\mathbb{C})$.  But if we suppose that
$\theta^h$ is invariant under ${\frak g}$ we have for a $H\in {\frak h}$ and
$\kappa=  \rho^{\mathbb{C}}_{|V}$ as above
\begin{eqnarray*}
0&=& \theta^h (\kappa(JH)v,w) + \theta^h (v, \kappa(JH)w)
\\
&\stackrel{p.d.}{=}& h^\mathbb{C}( \kappa(JH)v,\overline{w}) +
h^\mathbb{C}(v, \overline{\kappa(JH)w})
\\
&=&
h^\mathbb{C}(  \rho^{\mathbb{C}}
(JH)v,\overline{w}) + h^\mathbb{C}(v, \overline{  \rho^{\mathbb{C}} (JH)w})
\\
&\stackrel{H\in {\frak p}_+}{=}&
i \left(h^\mathbb{C}(  \rho^{\mathbb{C}}
(H)v,\overline{w}) - h^\mathbb{C}(v, \overline{  \rho^{\mathbb{C}} (H)w})
\right)
\\
&=&
i \left(\theta^h (\kappa(H)v,w) - \theta^h (v, \kappa(H)w) \right)
\\
&\stackrel{\theta^h\mbox{ invariant}}{=}&2i \theta^h (\kappa(H)v,w)
\end{eqnarray*}
for all $H\in {\frak h}$, $v,w\in V$. This means ${\frak h}\subset ker\
\kappa=0$ .
\end{proof}
We also can show the other direction.

\begin{satz}
Let $ {\frak g} $ be a real Lie algebra, $\kappa$ be a complex representation
of non-real type (of complex or quaternionic type), i.e. $\rho=\kappa_\mathbb{R}$ is
irreducible. Then holds:
\begin{enumerate}
\item If $ \kappa$ is unitary with respect to $\theta$ or orthogonal with
respect to $ \beta$, then $\rho$ is orthogonal with respect to $h$ and $
\theta^h= \theta$ or $ \beta^h= \beta$.
\item If $ \kappa$ is anti-unitary with respect to $\theta$ or symplectic with
respect to $ \beta$, then $\rho$ is symplectic with respect to $h$ and $
\theta^h= \theta$ or $ \beta^h= \beta$.
\end{enumerate}
\end{satz}

\begin{proof}
We define a bilinear form on $E=V_\mathbb{R}$ by
\[h(x,y):= Re\ \theta(x-iJx, y-iJy)\makebox[1cm][c]{ or }h(x,y):= Re\ \beta(x-
iJx, y-iJy).\]
This form is invariant and --- since $Re\  i z = -Im\  z$ --- also non-degenerate.
(The difference to real type is that here the arguments in $\theta/
\beta$ run over the whole complex vector space $V$.) $h$ is symmetric if
$\kappa$ is unitary or orthogonal and anti-symmetric if $ \beta$ is anti-symmetric
or anti-unitary.  The fact that the extensions  are equal to $\theta$
resp. $ \beta$ follows from the formulas
(\ref{h-theta}) and (\ref{h-omega}).
\end{proof}
Again we have the following correspondence:
\begin{eqnarray}\label{typ2equivalenz}
\left\{
\rho\mbox{ real, non-real type, orthogonal}
\right\}_{/\sim}&\leftrightarrow&
\left\{\begin{array}{l}
 \kappa \mbox{ complex, non-real type,}\\
\mbox{unitary or orthogonal}
\end{array}
\right\}_{/\approx}\\
\label{typ2equivalenz2}
\left\{
\rho\mbox{ real, non-real type, symplectic}
\right\}_{/\sim}&\leftrightarrow&
\left\{\begin{array}{l}
 \kappa \mbox{ complex, non-real type,}\\
 \mbox{symplectic or anti-unitary}
\end{array}
\right\}_{/\approx}
\end{eqnarray}

The fact that a complex representation is unitary if and only if it is anti-unitary
(the anti-hermitian form is $i\theta$) implies that a real, orthogonal
representation of non-real type with non-degenerate $\theta^h$ on the corresponding
complex representation is also symplectic.  This corresponds to the equality of
real matrix algebras:
\begin{eqnarray*}
 \mathfrak{u} (n)&=& \mathfrak{so}(2n) \cap \mathfrak{sp}(2n)\\
&=& \left\{X\in \mathfrak{gl}(2n)\left|X^t=-X\right.\right\}
\cap
\left\{\left.\left(\begin{array}{cc}A&B\\C&-A^t\end{array}\right) \right|
A,B,C\in \mathfrak{gl}(n), B^t=B, C^t=C \right\}\\
&=&
\left\{\left.\left(\begin{array}{cc}A&B\\-B&A\end{array}\right) \right|
A,B\in \mathfrak{gl}(n), A^t=-A, B^t=B, \right\}.
\end{eqnarray*}
I.e. if a complex representation $\kappa$ of non-real type is unitary, then
$\kappa_\mathbb{R}$ is orthogonal and symplectic.

Furthermore one proves the following
\begin{lem}
Let $h$ be symmetric, $\beta^h$, $\theta^h$ as above and $J$ the complex
structure
on $E$. Then holds
\begin{enumerate}
\item
$ \beta^h=0$ if and only if $h(x,y)= h(Jx,Jy)$ for all $x,y\in E$.
\item
$ \theta^h=0$ if and only if $h(x,y)= -h(Jx,Jy)$ for all $x,y\in E$.
\end{enumerate}\label{jh}
\end{lem}

\begin{proof}
If we write every element of $V=V_+$ in the form (\ref{v+form}) we get the
proposition due to formulas (\ref{h-omega}) and (\ref{h-theta}).
\end{proof}

We will now prove the main result for the case that  $h$ is positive definite.

\begin{satz}\label{main1}
Let $\rho$ be irreducible of non-real type and orthogonal with respect to $h$
where $h$ is {\bf positive definite}.

Then the corresponding complex representation $\kappa$ of non-real type  is unitary,
with respect to a positive definite hermitian form, which is the standard
hermitian form for representations of compact Lie groups/Lie algebras.

$\kappa$ is not orthogonal, i.e. the linear extension $ \beta^h$ of $h$
vanishes on $V\times V$.
\end{satz}

\begin{proof}
We can prove this in two ways.

If $\theta^h$ is degenerate, then it is zero and we have by lemma \ref{jh}
that $h(x,x)=-h(JxJx)$. But this is not possible if $h$ is positive definite. So $
\theta^h$ is non degenerate and by formula (\ref{h-theta}) positive definite,
since $h$ is positive definite.  But the existence of a positive definite
hermitian form entails by corollary \ref{+unitary} for non-real type representations,
i.e. of complex or quaternionic type, that the representation can not be orthogonal. So
$\beta^h=0$.

An easier way to argue is that representations of compact Lie algebras are
unitary with respect to a standard positive definite hermitian form.  This form
is unique and thats why equal to $\theta^h$ and by corollary \ref{+unitary}  the
representation can not be orthogonal. \end{proof}

\end{appendix}

\bibliography{GEOBIB,SPINBIB,HOLBIB,ALGBIB,thomas}
\end{document}